\documentclass{amsart}
\usepackage{geometry}                
\usepackage{graphicx}
\usepackage{overpic}
\usepackage{amsmath}
\usepackage{enumitem}
\usepackage{tikz}
\usepackage{hyperref}
\usetikzlibrary{shapes.geometric}
\usetikzlibrary{3d,calc}
\usepackage{amssymb}
\usepackage{mathtools}
\usepackage{epstopdf,amsmath}
\usepackage{amsthm}
\usepackage{xcolor}
\usepackage{upref}
\numberwithin{equation}{section}

\newtheorem{prop}{Proposition}[section]
\newtheorem{lemma}[prop]{Lemma}
\newtheorem{cor}[prop]{Corollary}
\newtheorem{theorem}[prop]{Theorem}
\newtheorem{conj}[prop]{Conjecture}
\newtheorem{definition}[prop]{Definition}
\newtheorem{claim}[prop]{Claim}
\newtheorem*{lemma*}{Lemma}
\newtheorem*{theorem*}{Theorem}

\DeclareMathOperator{\supp}{supp}

\newtheorem{remark}[prop]{Remark}

\newcommand{\RR}{\mathbb{R}}

\newcommand{\pp}{\mathbin{\!/\mkern-5mu/\!}}

\title{A restriction estimate in $\RR^3$ using brooms}

\author{Hong Wang}

\begin{document}

\begin {abstract} If $f$ is a function supported on the truncated paraboloid in $\mathbb{R}^3$ and $E$ is the corresponding extension operator, then we prove that for all $p> 3+ 3/13$, $\|Ef\|_{L^p(\mathbb{R}^3)} \leq C\|f\|_{L^{\infty}}$. The proof combines Wolff's two ends argument with polynomial partitioning techniques. We also observe some geometric structures in wave packets. 
\end {abstract}

\maketitle
\setcounter{tocdepth}{1}
\tableofcontents

\section{Introduction}

We strengthen a result of Guth \cite{Guth1} about the restriction conjecture (Conjecture \ref{conj: rest}) in $\mathbb{R}^3$. 

\begin{definition}
	Given a complex valued function $f$ supported in the unit disk in $\mathbb{R}^{n-1}$, we define 
	\[Ef(x_1, \dots, x_n): =\int_{|\xi|\leq 1} e^{i(\xi_1x_1+\cdots +\xi_{n-1} x_{n-1} +|\xi|^2 x_n)} f(\xi) d\xi.\]
\end{definition}
 Elias Stein \cite{Stein2} made the restriction conjecture about $Ef$ in the 1960s. 
 
 \begin{conj}[Stein \cite{Stein2}] \label{conj: rest}
 	For any $p>\frac{2n}{n-1}$ and $q\geq \frac{(n+1) q'}{n-1}$, we have 
 	\begin{equation}
 	\|Ef\|_{L^p(\mathbb{R}^n)} \leq C_{p,q} \|f\|_{L^q(\mathbb{R}^{n-1})}. 
 	\end{equation}
 \end{conj}

We refer to a survey of Tao \cite{Tao4} for a good presentation of Stein's conjecture. The conjecture in $\mathbb{R}^2$ was proven by Fefferman \cite{Fefferman} and Zygmund \cite{Zygmund}. The conjecture in $\mathbb{R}^3$ has been studied by various mathematicians but has not been proved for all $p$ and $q$. We are especially interested on the case $q=\infty$, because it implies the estimates for other $q$ with $q>p$ using the factorization theory of Nikishin and Pisier (see \cite{Bourgian}).  In this paper, we study the conjecture in $\mathbb{R}^3$ and for $q=\infty$, which asks if 
\begin{equation}\label{eq: rest}
\|Ef\|_{L^p(\mathbb{R}^3)} \leq C(p, S) \|f\|_{L^{\infty}(\mathbb{R}^2)}
\end{equation}
for all $p>3$. 

In 2002, Tao \cite{Tao3} obtained \eqref{eq: rest} for $p>3+1/3$ using the so-called two ends argument, which was introduced by Wolff \cite{Wolff}. In 2010, this estimate was improved by Bourgain and Guth in \cite{BG} for $p>3+5/17=3/29...$ (see Section 4.8 of \cite{BG}). In 2014, by introducing polynomial partitioning techniques Guth \cite{Guth1} improved  the range of $p$ to $p>2+1/4=3.25$, which is the previous best known estimate. In this paper, we provide a small improvement for $p> 3+3/13 = 3.23...$ based on the polynomial partitioning techniques and the two ends arguments. 

\begin{remark}
	The two ends argument involves the idea of splitting $Ef$ into a global part and a local part and an inductive argument for the local part, which can be found in Tao's paper \cite{Tao3}. The precise way of splitting the function in this paper is different from the one in \cite{Tao3}. But the inductive argument for the local part is general and works for different ways of splitting  the function. 
\end{remark}

\begin{theorem}\label{thm: main}
	The restriction estimate \eqref{eq: rest} holds for all $p> 3+3/13$. 
\end{theorem}

To prove his estimate, Guth \cite{Guth1} splits $Ef$ into a narrow part and a broad part (see Subsection~\ref{subsection: poly}). The narrow part can be controlled using induction on scale. The broad part, which is more difficult, is handled by an estimate for $\|Ef\|_{BL^p(B_R)}$, the $BL^p$--norm of $Ef$ inside a large ball $B_R$ (see Subsection~\ref{subsection: bl}). In this paper, we prove a stronger estimate about $\|Ef\|_{BL^p(B_R)}$ in Theorem~\ref{thm: induction2}.
Obtaining Theorem~\ref{thm: main} from Theorem~\ref{thm: induction2} requires no modification to Guth's argument in \cite{Guth1}, where he used an $\epsilon$--removal argument due to Tao \cite{Tao2}, so we omit the proof of Theorem~\ref{thm: main} from this paper. 


\begin{theorem}\label{thm: induction2}
	If $f$ is supported on the unit disk in $\mathbb{R}^2$, then for any small $\epsilon>0$, there exists a large constant $C_{\epsilon}$ depending only on $\epsilon$ such that for any large enough radius $R$, and for $p\geq 3+3/13$, we have 
	\begin{equation}\label{eq: induction2}
	\|Ef\|_{BL^p(B_R)} \leq C_{\epsilon} R^{\epsilon} \|f\|_{L^2}^{2/p} \underset{d(\theta)=R^{-1/2}}{\max} \|f_{\theta}\|_{L^2_{avg}(\theta)}^{1-2/p},
	\end{equation}
	where $\|f_{\theta}\|_{L^2_{avg}(\theta)}^2 \coloneqq Vol(\theta)^{-1}\|f_{\theta}\|_{L^2(\theta)}^2.$
\end{theorem}
The $BL^p$--norm of $Ef$ is introduced in Guth's paper \cite{Guth1} and \cite{Guth2} to capture the difficult part of $\|Ef\|_{L^p}$. We give its full description in Subsection~\ref{subsection: bl}. One could think of $\|Ef\|_{BL^p(B_R)}$ as $\|Ef\|_{L^p(B_R)}$ for most of this paper and especially in the introduction. On the right-hand side of \eqref{eq: induction2}, the term $ \underset{d(\theta)=R^{-1/2}}{\max} \|f_{\theta}\|_{L^2_{avg}(\theta)} \leq \|f\|_{L^{\infty}}$ plays a similar role as $\|f\|_{L^{\infty}}$ but it works better for induction.

We begin by highlighting some of Guth's arguments to illustrate how our idea grows out of Guth's ideas. The argument of Guth relies on induction on both scale (the radius of the spatial ball $B_R$) and the $L^2$--mass $\|f\|_{L^2}$. His argument consists of an application of the following trichotomy together with induction on scale. As a result of the polynomial partitioning, the function $Ef$ can be split into a cellular, a transverse and a tangential term (see Subsection~\ref{subsection: poly}). The cellular and the transversal contribution is estimated similarly, using the induction hypothesis on mass and radius. The tangential term is estimated directly, without appealing to induction. 

The unconditional estimate for the tangential term remains favorable if the induction hypothesis is changed to accommodate the paper's Theorem~\ref{thm: induction2}, which is reflected in reducing the exponent weight on the $L^2$--mass. However, in this new setup, the estimate for the cellular part is no longer an immediate consequence of induction on radius. Most of the novelty of this paper goes into finding a new way to deal with the cellular contribution. Our argument contains a multi-step iteration. And the scale of the wave packets, the building blocks of $Ef$ (see Subsection~\ref{subsection: wp}), changes throughout the iteration process. Essentially, at each step we split the cellular component into two parts: a local part and a global part. The local part will be estimated using induction on the radius (see Lemma~\ref{lem: telocal}). The global part of $Ef$, in the critical cases,  is mainly concentrated in thin neighborhoods of algebraic varieties intersecting $B_R$ and needs a delicate analysis that forces us to introduce a new geometric object we call the ``broom''.  A large part of the argument in this paper goes into quantifying how brooms interact with thin neighborhoods of algebraic varieties. 

\begin{figure}	
	\begin{tikzpicture}[ scale=0.5]
	\draw[red] (0,0)--(2,1)--(4.5,1)--(2.5,0)--cycle;
	\draw[red](0,0)--(0,-0.3)--(2.5,-0.3)--(4.5,0.7)--(4.5,1);
	\draw[red](2.5,-0.3)--(2.5,0) node[black, below right]{$S$};
	\node(a)[cylinder,rotate=205, shift={(0.8,0.20)},blue, draw, minimum height=55mm, minimum width=5mm]{};
	\node(b)[cylinder,rotate=215, shift={(0.9,0.40)},blue, draw, minimum height=55mm, minimum width=5mm]{};
	\node(c)[cylinder,rotate=225, shift={(1,0.60)},blue, draw, minimum height=55mm, minimum width=5mm]{};
	\draw[blue, dashed] (3.25,1)--(3.25,6)--(-6.75,1)--(-6.75,-9)--(3.25,-4)--cycle;
	\draw[black,thick](3.75,1)--(3.75,1.5)--(3.25,1.5) node[above right]{perpendicular};
	\end{tikzpicture}
	\caption{A broom}
	\label{sbroom}
\end{figure}

To describe a broom we recall the wave packet decomposition of $Ef$ introduced by Bourgain \cite{Bourgian}. The wave packet decomposition says that inside a large ball of radius $R$,  we can decompose $Ef$ into a sum over wave packets $Ef_{\theta,v}$. Each wave packet $Ef_{\theta,v}$ is essentially supported in a tube $T_{\theta,v}$ of length $R$, radius $R^{1/2+\delta}$ for some small $\delta >0$. The axis of $T_{\theta, v}$ points in a direction depending only on $\theta$ and the location of $T_{\theta,v}$ is described by $v$. The absolute value $|Ef_{\theta,v}|$ of a wave packet is approximately a constant function on $T_{\theta,v}$. 

Roughly speaking, a broom is a collection of large tubes $T_{\theta,v}$ satisfying certain properties with respect to  a thin neighborhood of a plane. More precisely, for some $R^{1/2+\delta}< r< R$, let $S$ denote the $r^{1/2}$--neighborhood of a plane $\Sigma$ inside a ball of radius $r$. The dimensions of $S$ are $r\times r\times r^{1/2}$. 

In particular, $S$ is an example of the thin neighborhood of an aforementioned algebraic variety where $Ef$ is concentrated on.  For simplicity, we assume in the introduction that every such thin neighborhood is of the form of $S$. Now we are ready to give a better description of a broom with respect to $S$ as illustrated in Figure~\ref{sbroom}.

 A broom consists of tubes $T_{\theta,v}$ which 
\begin{enumerate}
	\item intersect at a common point in $S$;
	\item lie inside the $R^{1/2+\delta}$--neighborhood of a plane perpendicular to $S$;
	\item form an angle with $S$ of at most $r^{-1/2}$. 
\end{enumerate}

The key ingredient of our proof is an improved estimate about the $L^2$--mass of $Ef$ on a typical $S$ compared to the $L^2$--mass of $Ef$ on a ball $B_r$ containing $S$. We explain using two examples how brooms help to obtain such an estimate. The first example is when  $Ef$ is concentrated on only  one wave packet $Ef_{\theta,v}$ from a broom.  Since $|Ef_{\theta,v}|$ is approximately a constant function on $T_{\theta,v}$, its $L^2$--mass is evenly distributed on $T_{\theta,v}$. So the $L^2$--mass contained in $S$ is only a small proportion because the volume of $S$ is much smaller than the volume of $T_{\theta,v}\cap B_r$.  
\begin{center}
	\includegraphics[scale=0.2]{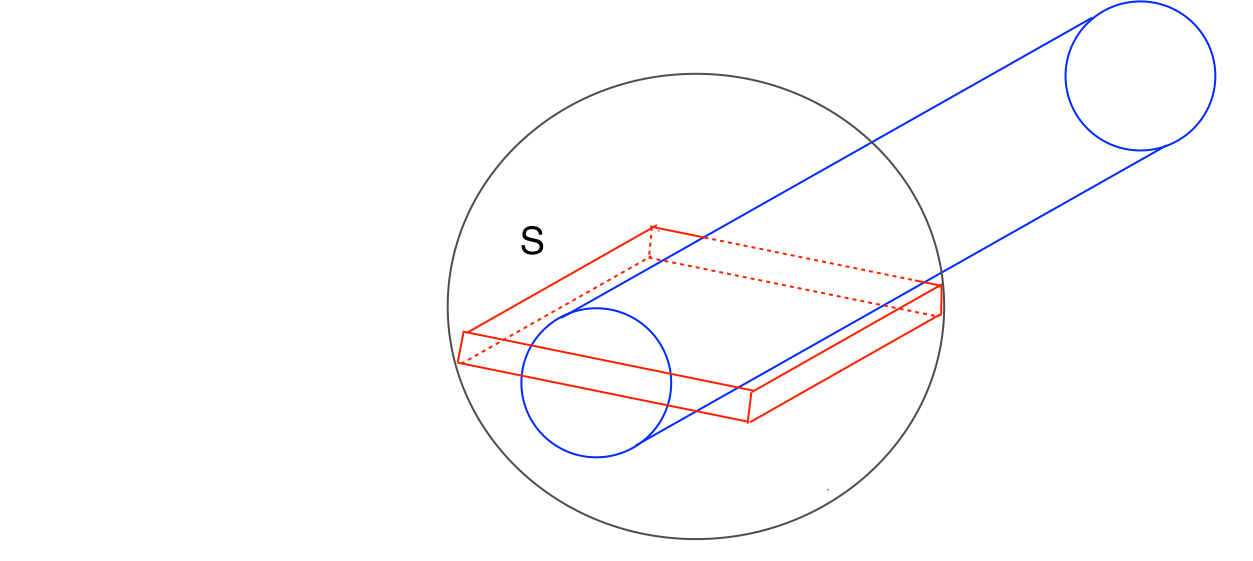}
\end{center}

The second example is when $Ef$ is concentrated on  a large number of wave packets from a broom. Since the tubes in a broom intersect in $S$ on one end and spread out on the other end, a typical $S'$ on the other end intersects the broom in only a few tubes, which implies that the $L^2$--mass of $Ef$ on $S'$ is small because most tubes in the broom miss $S'$. 

 \begin{center}
	\includegraphics[scale=0.25]{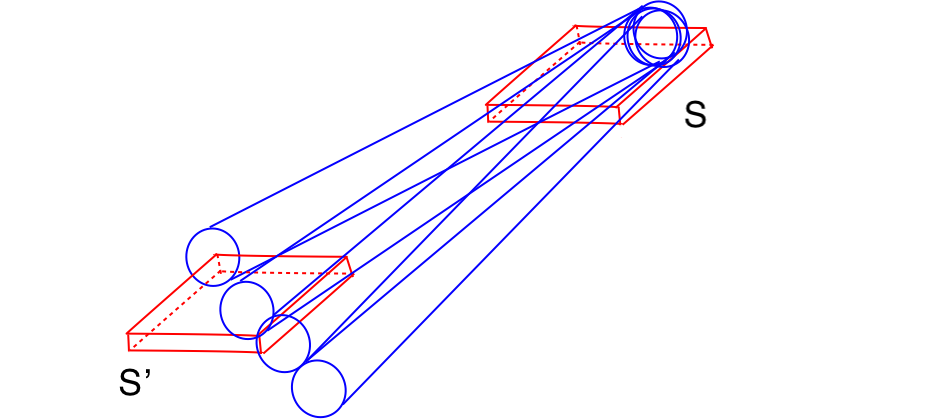}
\end{center}

\textbf{The rest of this paper is organized in the following way.}

In Subsection~\ref{subsection: poly}, we recall Guth's proof for $p> 3.25$ in \cite{Guth1} using polynomial partitioning, which is the starting point of our proof.  In Subsection~\ref{subsection: te}, we recall Wolff's two ends argument and estimate the local part in Lemma~\ref{lem: telocal}.
In Section~\ref{section: poly}, we apply polynomial partitioning iteratively and prove a polynomial structure lemma, Lemma~\ref{lem: recursion}, which says that $Ef$ is mainly concentrated on thin neighborhoods of a collection of algebraic surfaces in the difficult cases, while the easy cases are proved in Section~\ref{section: easy}. In Section~\ref{section: broom}, we define the key geometric objects: brooms and bushes.  In Section~\ref{section: sim}, we give the precise way of splitting $Ef$  into global part and local part using the broom structure, assuming a geometric lemma, Lemma~\ref{lem: plane}, saying that the algebraic surfaces in Lemma~\ref{lem: recursion}  can be viewed as planes. We give the proof of Lemma~\ref{lem: plane} in Section~\ref{section: plane}.  In Section~\ref{section: l2}, we estimate the global part of $Ef$ using the broom structure.   Finally, in Section~\ref{section: proofmain}, we conclude the proof by combining the lemmas from the previous sections.  Section~\ref{section: 38} is devoted to prove a technical lemma, Lemma~\ref{lem: find scale}, in Section~\ref{section: poly}. 

\vspace{10pt}

\textbf{Notation.} If $X$ is a finite set, we use $|X|$ to denote its cardinality.  If $X$ is an infinite,  Lebesgue measurable set,  we use $|X|$ to denote its Lebesgue measure.  We use $d(\theta)$ to denote the diameter of a small disk or the diameter of a small cap on a paraboloid. We use $B(x, r)$ to denote a ball of radius $r$ centered at $x$, and we use $B_r$ to denote a ball of radius $r$. When $A, B>0$, we use $A\lesssim B$ or $A= O(B)$ to denote the estimate $A \leq CB$ where $C$ is an absolute positive constant.  For a parameter $\epsilon$, We use $A \lesssim_{\epsilon} B$ or $A= O_{\epsilon}(B)$ to denote the estimate $A \leq C_{\epsilon} B$ where $C_{\epsilon}>0$ is a constant only depending on $\epsilon$.  We say that two positive quantities $A \sim \lambda $ if $\lambda \lesssim A \lesssim \lambda$.  If $T$ is a compact symmetric convex body with center $C_T$ and $k>0$, then we use $kT$ to denote the rescaling of $T$ by a factor of $k$ with the same center $C_T$.

\subsection*{Acknowledgments}
 I would like to thank my advisor Larry Guth for his guidance and encouragement throughout this project. I would also like to thank Donghao Wang and Ruixiang Zhang for helpful discussion on  Proposition~\ref{lem: plane}. I would like to thank the anonymous referees for the thorough reading and for many writing suggestions that improve the presentation of this paper.   I would like to also thank Susan Ruff and Shuanglin Shao for their writing advice. This research is partially supported by Larry Guth's Simons Investigator Award and partially supported by National Science Foundation (Grant No. DMS-1638352) and the S.S. Chern Foundation for Mathematics Research Fund. 

\section{Preliminaries}
In this section, we first recall the wave packet decomposition and the definition of the $BL^p$--norm. Then we proceed to the proof sketch of Guth's result \cite[Theorem 3.1]{Guth1}. In the end of this section, we review  Wolff's two ends argument.
\subsection{Wave packet decomposition}\label{subsection: wp}
We start with the definition of wave packet decomposition at scale $R$ following the description in \cite[Section 3]{Guth2}.

\begin{definition}\label{def: wp}
	We decompose the unit disk in $\mathbb{R}^2$ into finitely overlapping small disks $\theta$ of radius $R^{-1/2}$. These small disks are referred to as $R^{-1/2}$--caps. Let $\psi_{\theta}$ be a smooth partition of unity adapted to $\{\theta\}$, and write $f=\sum_{\theta}\psi_{\theta} f$ and define $f_{\theta} \coloneqq\psi_{\theta}f$. We cover $\mathbb{R}^2$ by finitely overlapping balls of radius about $R^{\frac{1+\delta}{2}}$, centered at vectors $v\in R^{\frac{1+\delta}{2}} \mathbb{Z}^2$, where $\delta$ is a small number satisfying $\epsilon^7 <\delta \leq \epsilon^3$.  Let $\eta_v$ be a smooth partition of unity adapted to this cover. We can now decompose
	$$f=\sum_{\theta,v} (\eta_v (\psi_{\theta}f)^{\wedge})^{\vee} = \sum_{\theta,v} \eta_v^{\vee}\ast (\psi_{\theta}f).$$
	We choose smooth functions $\tilde{\psi}_{\theta}$ such that $\tilde{\psi}_{\theta}$ is supported on $\theta$, but $\tilde{\psi}_{\theta}=1$ on a $cR^{-1/2}$ neighborhood of the support of $\psi_{\theta}$ for a small constant $c>0$. We define
	\[ f_{\theta,v}\coloneqq \tilde{\psi}_{\theta} [\eta_{v}^{\vee} \ast(\psi_{\theta} f)].\]
	Since $\eta_v^{\vee} (x)$ is rapidly decaying for $|x|\gtrsim R^{\frac{1-\delta}{2}}$, we have 
	\begin{equation}\label{eq: wpR}
	 f=\sum_{(\theta,v): \,  d(\theta)=R^{-1/2}} f_{\theta,v} +\text{RapDec}(R)\|f\|_{L^2}.
	\end{equation}
	Here the notation $d(\theta)$ denotes the diameter of $\theta$, and  $\text{RapDec}(R)$ means that the quantity is bounded by $O_N(R^{-N})$ for any large integer $N>0$. 
\end{definition}
The wave packets $Ef_{\theta,v}$ are the building blocks for $Ef$ and satisfy two useful properties.  The first property is that the functions $f_{\theta,v}$ are approximately orthogonal: for any set $\mathbb{T}$ of $(\theta,v)$, we have (see inequality (3.1) of \cite{Guth2})
\begin{equation}\label{eq: l2orthogonality}
\sum_{(\theta,v)\in \mathbb{T}} \|f_{\theta,v} \|_{L^2}^2 \sim \|\sum_{(\theta,v)\in \mathbb{T}} f\|_{L^2}^2.
\end{equation}
The second property is that on the ball $B(0, R)$ the function $Ef_{\theta,v}$ is essentially supported on a tube $T_{\theta,v}$: 
\[T_{\theta,v}\coloneqq\{ (x', x_3)\in B(0, R), |x'+2x_3\omega_{\theta} +v|\leq R^{1/2+\delta}\},\] 
where $\omega_{\theta}$ is the center of the cap $\theta$. More precisely, we have the following lemma. 

\begin{lemma}[{\cite[Lemma 3.1]{Guth2}}]\label{lem: wpsupp}
	If $x\in B(0,R)\setminus T_{\theta, v}$, then
	\[|Ef_{\theta,v}(x)|\leq \text{RapDec}(R) \|f\|_{L^2}.\]
\end{lemma}

The tube $T_{\theta,v}$ is a cylinder of length $R$ and radius about $R^{1/2+\delta}$, pointing in the direction $G(\omega_{\theta})$ where 
\begin{equation}\label{eq: Gtheta}
 G(\omega) : = \frac{(-2\omega_1, -2\omega_2, 1)}{|(-2\omega_1, -2\omega_2, 1)|}.
\end{equation}
Since the directions in $G(\theta)$ are within a range of $R^{-1/2}$, for a fixed $\theta$ and $v$, any direction in $G(\theta)$ defines essentially the same tube. 

\begin{definition}\label{def: concentrate}
If $\mathbb{T}$ is a set of $(\theta,v)$, we say that $f$ is concentrated on wave packets from $\mathbb{T}$ if $f = \sum_{(\theta,v)\in \mathbb{T}} f_{\theta,v} +\text{RapDec}(R)\|f\|_{L^2}.$
\end{definition}
 Since each pair $(\theta,v)$ corresponds to a unique tube in $B(0,R)$, we also use $\mathbb{T}$ to denote the collection of tubes $T_{\theta,v}$. 
We have discussed  wave packet decomposition over a ball centered at the origin. Sometimes we will also need to do wave packet decomposition over a ball $B(y, \rho)$ centered elsewhere and of a different radius.  We apply a translation to recenter $B(y, \rho)$ at the origin and write $Ef(x)= E\tilde{f}(\tilde{x})$ where $x= y+\tilde{x}$ for some $\tilde{x}\in B(0,\rho)$ and 
\[\tilde{f}(\xi)\coloneqq e^{ i (\langle y', \xi)\rangle + y_3|\xi|^2)} f(\xi) +\text{RapDec}(\rho).\]
We can decompose  $\tilde{f}$  into scale $\rho$ wave packets (with $\rho$ in the place of $R$ in Definition~\ref{def: wp}): 
\begin{equation}\label{eq: wprho}
\tilde{f} = \sum_{(\tilde{\theta}, \tilde{v}): \,d(\tilde{\theta})=\rho^{-1/2}} \tilde{f}_{\tilde{\theta}, \tilde{v}}.
\end{equation}

When $R^{1/2}\leq \rho \leq R$, we recall how two wave packet  decompositions \eqref{eq: wpR} and \eqref{eq: wprho} relate to one another. 

\begin{lemma}[\cite{Guth2} {\cite[Lemma 7.2]{HR}}]\label{lem: slwp}
	If $f$ is concentrated on a collection  $\mathbb{T}$ of scale $R$ wave packets, then $\tilde{f}$ is concentrated on a set $\tilde{\mathbb{T}}$ of scale $\rho$ wave packets with the following property: for every $(\tilde{\theta}, \tilde{v})\in \tilde{\mathbb{T}}$ there exists $(\theta, v)\in \mathbb{T}$ such that 
	\begin{enumerate}
		\item $ \mathit{dist}_H(T_{\tilde{\theta},\tilde{v}} + y, T_{\theta,v}\cap B(y, \rho)) \lesssim R^{1/2+\delta}$;
		\item $\mathit{Angle} \,(G(\omega_{\theta}), G(\omega_{\tilde{\theta}})) \lesssim \rho^{-1/2}$.
	\end{enumerate}
Here $\mathit{dist}_H$ is the Hausdorff distance. 
\end{lemma}
 Lemma~\ref{lem: slwp}  follows from several lemmas  in  \cite[Section 7]{Guth2}, but this version of statement is borrowed from \cite{HR}. 
Lemma~\ref{lem: slwp} says that for every $(\tilde{\theta}, \tilde{v})\in \tilde{\mathbb{T}}$, there exists a parent (in fact, maybe multiple parents) such that the translated tube $T_{\tilde{\theta}, \tilde{v}} +y$ lies close to $T_{\theta,v}$ and points in a similar direction to $T_{\theta,v}$.  When the underlying ball $B(y, \rho)$ is clear from the context, to ease the notation, we write 
$$ f =\sum_{(\tau , w): \, d(\tau)=\rho^{-1/2}} f_{\tau, w}$$
where $\tau =\tilde{\theta}$,  and $w$ corresponds to a unique $\tilde{v}$ in the sense that $T_{\tau,w}= T_{\tilde{\theta}, \tilde{v}}+y$, and  $Ef_{\tau,w}(x) = E\tilde{f}_{\tilde{\theta}, \tilde{v}}(x-y).$  By Lemma~\ref{lem: wpsupp},  in $B(y, \rho)$, $Ef_{\tau, w}$ is essentially supported on the tube $T_{\tau,w}$. If $\mathbb{T}$ is a set of tubes $T_{\tau,w}$, we say that $f$ is concentrated on wave packets from $\mathbb{T}$ if $f=\sum_{T_{\tau,w}\in\mathbb{T}} f_{\tau,w} +\mathit{RapDec}(\rho)\|f\|_{L^2}.$

\begin{lemma}[{\cite[Lemma 3.4]{Guth2}}]\label{lem: l2}
	Suppose that $f$ is concentrated on a set of wave packets from $\mathbb{T}$ and that for every $(\theta,v)\in \mathbb{T}$, $T_{\theta,v}\cap B_r\neq \emptyset$ for some $r\geq R^{1/2+\delta}$, then 
	\[\|Ef\|_{L^2(10B_r)}^2 \sim r\|f\|_{L^2}^2.\]
\end{lemma}

We also recall two lemmas regarding the locally constant property of a function whose Fourier transform is supported on a compact convex set. 
\begin{lemma}[{\cite[Lemma 6.1]{GWZ}}]\label{lem: lc1}
	Let $\Theta$ be a compact symmetric convex set centered  at $C_{\Theta} \in \mathbb{R}^3$. If $\widehat{g}_{\Theta}$ is supported in $\Theta$ and $T_{\Theta}$ is the dual convex  $\Theta^* \coloneqq \{ y: |y\cdot x-C_{\Theta}|\leq 1, \forall x\in \Theta\}$, then there exists a positive function $\eta_{T_{\Theta}}$ satisfying the following properties:
\begin{enumerate}
	\item  $\eta_{T_{\Theta}}$ is 
	essentially supported on $10T_{\Theta}$ and rapidly decaying away from it: for any integer $N>0$, there exists a constant $C_N$ such that $\eta_{T_{\Theta}}(x) \leq C_N(1+n(x, 10 T_{\Theta}))^{-N}$ where $n(x, 10T_{\Theta})$ is the smallest positive integer $n$ such that $x\in 10n T_{\Theta}$, 
	\item $\|\eta_{T_{\Theta}}\|_{L^1}\lesssim 1$, 
	\item \begin{equation}\label{eq: 2.5}
	|g_{\Theta}|\leq \sum_{T\pp T_{\Theta}} c_T \chi_T\leq |g_{\Theta}|\ast \eta_{T_{\Theta}}
	\end{equation}
	where $c_T \coloneqq \max_{x\in T} |g_{\Theta}|(x)$ and the sum $\sum_{T\pp T_{\Theta}}$ is  over a finitely overlapping cover $\{T\}$ of  $\mathbb{R}^3$ with each $T \pp T_{\Theta}$. Here $T\pp T_{\Theta}$ means that $T$ is a  translated copy of  $T_{\Theta}$. 
\end{enumerate}
\end{lemma}

\begin{lemma}[{\cite[Lemma 6.2]{GWZ}}]\label{lem: lc2}
	Let $g_{\Theta}$ and $\eta_{T_{\Theta}}$ be defined as in Lemma~\ref{lem: lc1} and $T\pp T_{\Theta}$.  Then for any integer $N>0$,  there exists a positive function $w_T=1$ on $10T$ and $w_T(x) \leq C_N (1+n(x, T))^{-N}$ such that for any $1\leq p <\infty$, 
	\begin{equation}\label{eq: 2.6}
	\int_T (|g_{\Theta}|\ast \eta_{T_{\Theta}})^p \lesssim_p \int |g_{\Theta}|^p w_T.
	\end{equation}
\end{lemma}

We say that $g_{\Theta}$ is locally constant on  a translated copy $T$  of $ T_{\Theta}$, if \eqref{eq: 2.5} and \eqref{eq: 2.6} hold for $g_{\Theta}$ and $T$.   In particular, $\widehat{Ef_{\theta,v}}$ is supported on $\{(\omega, |\omega|^2) ,\omega\in \theta\}$ in distributional sense, which is contained in a small plate of radius $R^{-1/2}$ and thickness $R^{-1}$, whose dual convex is parallel to $T_{\theta,v}$ (up to an $R^{\delta}$--rescaling).  Hence, we say that $Ef_{\theta,v}$ is locally constant on $T_{\theta,v}$. 

Last, we recall the local $L^2$--orthogonality lemma.  

\begin{lemma}[{\cite[Lemma 6.1]{BD}}]\label{lem: l2loc}
	Let $Q$ be a cube of side length $\geq R^{-1}$ and $\{q\}$ be a finitely overlapping cover of $Q$ by a collection of cubes $q$ of side length $R^{-1}$. Then for each $B_R$ and each function $g_q$ with $\supp\widehat{g}_q \subset q$, we have
	\[\int_{B_R} |\sum_{q\subset Q} g_q|^2 \lesssim \sum_{q\subset Q} \int |g_q|^2 w_{B_R}\]
	where $w_{B_R}$ is a positive function $=1$ on $B_R$ and rapidly decaying outside of $B_R$: $w_{B_R}(x) \leq C_N (1+ n(x, B_R))^{-N}.$
\end{lemma}

\subsection{$BL^p$--norm}\label{subsection: bl}
We recall the $BL^p$--norm $\|Ef\|_{BL^p(B_R)}$ defined in Guth's papers \cite{Guth1, Guth2}.

We decompose the unit disk into finitely overlapping disks $\alpha$ of radius $K^{-1}$, where $1\leq K^{\epsilon^2} \leq R^{\epsilon^{12}/100}$, and $\epsilon$ is the small constant  in Theorem~\ref{thm: induction2}. We decompose $f=\sum_{\alpha} f_{\alpha}$ by a partition of unity with $f_{\alpha}$ supported in $\alpha$. The  wave packets at scale $R$ in $Ef_{\alpha}$ are those $Ef_{\theta,v}$ with $\theta\subset \alpha$. 

Now we are ready to define  the $BL^p$--norm. Let $A$ be a positive integer such that  $1\leq A^{\epsilon} \leq K^{\epsilon^2}$.  For each small ball $B_K$, we define 
\begin{equation}\label{eq: muEf}
\mu_{Ef} (B_K) \coloneqq \underset{V_1, \dots, V_A: \text{ lines in } \mathbb{R}^3}{\min}  \Big(   \underset{\alpha: \mathit{Angle}(G(\alpha), V_{a})\geq K^{-1} \text{ for all } 1\leq a \leq A}{\max} \, \int_{B_K} |Ef_{\alpha}|^p    \Big).
\end{equation}
If a set $U\subset \mathbb{R}^3$ is covered by finitely overlapping balls $B_K$, then we define $\mu_{Ef}(U)$ by summing over the covering: $\mu_{Ef}(U) = \sum_{B_K\subset U} \mu_{Ef}(B_K)$. 
Now we define the $BL^p$--norm on $B_R$ as:
\[ \|Ef\|_{BL^p_A (B_R)}^p \coloneqq \mu_{Ef}(B_R). \]

The $BL^p$--norm is closely related to the bilinear norm, which was studied by Tao in \cite{Tao3}. 
\begin{lemma}\label{lem: brbl} Suppose that $f$ is a function supported in the unit disk in $\mathbb{R}^2$,  $\mu_{Ef}$  is defined as in \eqref{eq: muEf}, and $\alpha_1, \alpha_2$ are caps of radius $K^{-1}$.  Then
	\[ \mu_{Ef}(B_K) \leq \int_{B(0, 2K)} \sum_{\alpha_1, \alpha_2 \text{ nonadjacent} } \int_{B_K} |Ef_{\alpha_1}(x -y) Ef_{\alpha_2}(x)|^{p/2} dx dy.\]
	Here the $K^{-1}$--caps $\alpha_1$ and $\alpha_2$ are said to be nonadjacent if $\mathit{dist}(\alpha_1, \alpha_2)\geq K^{-1}$. 
	\end{lemma}
\begin{proof}
	By Lemma~\ref{lem: lc1} and the definition of the $BL^p$--norm, 
	\begin{align*}
\mu_{Ef}(B_K) & \leq \sum_{\alpha_1, \alpha_2 \text{ nonadjacent }} \int_{B_K} |Ef_{\alpha_1}(x)|^{p/2}  dx\int_{B_K} |Ef_{\alpha_2}(x)|^{p/2} dx\\
& \leq  \int_{B(0, 2K)} \sum_{\alpha_1, \alpha_2 \text{ nonadjacent }} \int_{B_K} |Ef_{\alpha_1}(x-y)|^{p/2} \int_{B_K} |Ef_{\alpha_2}(x)|^{p/2}  dx dy.
	\end{align*}
	The last inequality follows from Fubini's theorem and the fact that $B_K \subset B(x, 2K)$ for any $x\in B_K$. 
\end{proof}

 A  similar statement of Lemma~\ref{lem: brbl} can be found in \cite[Lemma 3.8]{Guth1}. 
 Lemma~\ref{lem: brbl} implies that 
\begin{equation}\label{eq: br}
\|Ef\|_{BL^p_A(B_R)}^p \leq \int_{B(0, 2K)}  \sum_{\alpha_1, \alpha_2 \text{ nonadjacent} } \||Ef_{\alpha_1}(\cdot -y) Ef_{\alpha_2}|^{1/2}\|_{L^p({B_R})}^p dy . 
\end{equation}

One can also view the $BL^p$--norm as approximately an $L^p$--norm with broadness: if $\supp f$ lies inside a small cap of radius $K^{-1}$, then $\|Ef\|_{BL^p_A(B_R)}=0$.  Similar to the usual  $L^p$--norms,  the  triangle inequality and a version of  H\"{o}lder's inequality hold for the  $BL^p$--norm if the parameter   $A$ is allowed to change on the right-hand side.  More precisely, 
\begin{equation}\label{eq: triangle}
\|E(f+g)\|_{BL^p_{A_1+A_2}(U)} \lesssim  \|Ef\|_{BL^p_{A_1}(U)} + \|Ef\|_{BL^p_{A_2}(U)},
\end{equation}
and if $1\leq p, p_1, p_2<\infty$, $0\leq \alpha_1, \alpha_2\leq 1$ obey $\alpha_1+\alpha_2 =1$ and $\frac{1}{p}=\frac{\alpha_1}{p_1} +\frac{\alpha_2}{p_2}$, then 
\begin{equation}\label{eq: holder}
\|Ef\|_{BL^p_{A_1+A_2}(U)} \leq \|Ef\|_{BL^{p_1}_{A_1}(U)}^{\alpha_1} \|Ef\|_{BL^{p_2}_{A_2}(U)}^{\alpha_2}.
\end{equation}
See Lemma 4.1 and Lemma 4.2 in \cite[Section 4]{Guth2} for a proof of inequalities~\eqref{eq: triangle} and \eqref{eq: holder}.

The value of $A$ can change  from line to line whenever we  apply \eqref{eq: triangle} and \eqref{eq: holder} in Guth's argument (see page 13) as well as in the proof of Theorem~\ref{thm: induction2} (see Section~\ref{section: 38}  for more discussion). However, for the key arguments in this paper,   it is convenient to  drop the parameter $A$. 
In what follows, we omit $A$ from our notation unless it plays a role in the proof.

\begin{lemma}\label{lem: tao}
	If $f$ is a function supported in the unit disk in $\mathbb{R}^2$. Then for any $2\leq p\leq 10/3$, 
	\begin{equation}
	\|Ef\|_{BL^p(B_R)}^p \lesssim R^{\frac{5}{2}-\frac{3p}{4}+\epsilon^{10}} \|f\|_{L^2}^p.
	\end{equation}
\end{lemma}
\begin{proof}
	We claim that 
		\begin{equation}\label{eq: tao}
	\|Ef\|_{BL^{10/3}(B_R)} \lesssim R^{\epsilon^{10}} \|f\|_{L^2}.
	\end{equation}
		To see this, we apply  inequality~\eqref{eq: br} and Tao's bilinear restriction theorem \cite[Theorem 1.1]{Tao3},
	\begin{align*}
	\|Ef\|_{BL^{10/3}(B_R)} &\leq  \int_{B(0, 2K)} \sum_{\alpha_1, \alpha_2 \text{ nonadjacent} } \| |Ef_{\alpha_1}(\cdot -y) Ef_{\alpha_2}|^{1/2}\|_{L^{10/3}(B_R)}^{10/3} dy\\ 
	&\leq \int_{B(0, 2K)} \sum_{\alpha_1, \alpha_2 \text{ nonadjacent}} \|f_{\alpha_1}\|_{L^2}^{5/3} \|f_{\alpha_2}\|_{L^2}^{5/3} dy \\
	&\leq \mathit{Vol}(B(0, 2K))\cdot K^{2/3} \|f\|_{L^2}^{10/3}\\
	&\lesssim R^{\epsilon^{10}} \|f\|_{L^2}^{10/3}.
	\end{align*}
	By inequality~\eqref{eq: br} and Lemma~\ref{lem: l2}, we have 
	\begin{equation}\label{eq: l2br}
	\|Ef\|_{BL^2(B_R)}\lesssim R^{1/2+\epsilon^{10}} \|f\|_{L^2}. 
	\end{equation}
Now we combine 	\eqref{eq: holder}, \eqref{eq: tao} and \eqref{eq: l2br} to conclude that
	\begin{equation}\label{eq: 3p4}
	\|Ef\|_{BL^p(B_R)}^p \lesssim R^{\frac{5}{2}-\frac{3p}{4}+\epsilon^{10}} \|f\|_{L^2}^p. 
	\end{equation}
\end{proof}
\begin{remark}
	When $p\geq 10/3$, all restriction estimates proved in this paper follow from \eqref{eq: tao}.  For the rest of the paper, we only consider the range $2\leq p\leq 10/3$. 
\end{remark}
\subsection{Some basic reductions}\label{subsection: reduction}
We sort the wave packets $Ef_{\theta,v}$ according to the size of $\|f_{\theta,v}\|_{L^2}$. 
Given a fixed $B_R$, there are $\lesssim R^4$ wave packets intersecting $B_R$. Since the contribution from those wave packets with $\|f_{\theta,v}\|_{L^2} \leq R^{-10} \|f\|_{L^2}$ satisfies the inequality in Theorem~\ref{thm: induction2}, it suffices to consider the wave packets with $R^{-10} \|f\|_{L^2} \leq \|f_{\theta,v}\|_{L^2}$.  By inequality~\eqref{eq: l2orthogonality}, $\|f_{\theta,v}\|_{L^2} \lesssim \|f\|_{L^2}$.  For each dyadic number $\lambda$ with $R^{-10} \leq \lambda \lesssim 1$, let $f_{\lambda}$ denote the sum over wave packets with $\|f_{\theta,v}\|_{L^2}\sim \lambda$. 

Since there are $\lesssim \log R$ choices of $\lambda$, there exists a $\lambda_0$ such that 
\[ \|Ef_{\lambda_0}\|_{BL^p(B_R)}\gtrsim (\log R)^{-1} \|Ef\|_{BL^p(B_R)}.\]
Let $\mathbb{T}_0$  denote the collection of tubes $T_{\theta,v}$ with  $\|f_{\theta,v}\|_{L^2}\sim \lambda_0$. From now on, we assume that $f$ is 
 concentrated on wave packets from $\mathbb{T}_0$.  
\subsection{Polynomial partitioning}\label{subsection: poly}
In this section, we provide a proof sketch of Guth's theorem. 

\begin{theorem}[{\cite[Theorem 3.1]{Guth1}}]\label{thm: guth}
	If $f$ is supported in the unit disk in $\mathbb{R}^2$, then for any small $\epsilon>0$, there exists a large constant $C_{\epsilon}$ depending only on $\epsilon$ such that for any large radius $R$ and  any $p \geq 13/4$, 
	\[ \|Ef\|_{BL^p(B_R)} \leq C_{\epsilon} R^{\epsilon} \|f\|_{L^2}^{12/13} \underset{d(\theta)=R^{-1/2}}{\max} \|f_{\theta}\|_{L^2_{avg}(\theta)}^{1/13}.\]
\end{theorem}

The idea is to use the zero set of a polynomial to partition $\mu_{Ef}(B_R)= \|Ef\|_{BL^p(B_R)}^p$. After partitioning, we obtain a cellular part, a tangential part and a transversal part of $Ef$. We shall estimate the tangential part directly and use induction on the radius $R$ for the cellular part and the transversal part.

The base case of the induction is when $R\leq C_{\epsilon}'$ for some constant $C_{\epsilon}'$ only depending on $\epsilon$.  By Lemma~\ref{lem: tao},  $\|Ef\|_{BL^p(B_R)}^p \lesssim R^{\frac{5}{2}-\frac{3p}{4}} \|f\|_{L^2}^p \lesssim C_{\epsilon}' \|f\|_{L^2}^p$.  It suffices to choose $C_{\epsilon}$ sufficiently large  depending on $C_{\epsilon}'$.

 To proceed the induction, assume that Theorem~\ref{thm: guth} holds for all $\rho<R/2$.  We apply  \cite[Theorem 5.5]{Guth2} to find the partitioning polynomial. More precisely,   \cite[Theorem 5.5]{Guth2}  says  that for any degree $d\geq 1$, we can find a nonzero polynomial $P$ of degree $\leq d$ such that the complement of its zero set $Z(P)$ in $B_R$ is a union of $\sim d^3$ disjoint cells $U_i'$:
\begin{equation} \label{eq: ZP}
B_R\setminus Z(P)=\bigsqcup_i U_i',
\end{equation}
and the $BL^p$--norm is roughly the same in each cell: 
\begin{equation}\label{eq: ZPcell}
\|Ef\|_{BL^p(U_i')}^p \sim d^{-3} \|Ef\|_{BL^p(B_R)}^p. 
\end{equation}

We choose the degree $d \sim \log R$ to make the induction argument work. 
Furthermore, \cite[Theorem 5.5]{Guth2} says that  $Z(P)$ is a finite union of  smooth algebraic surfaces.

Now we introduce a slightly different treatment of the cells from that used in Guth's argument, which does not change the conclusion of Guth's theorem but will become beneficial when we introduce our new ideas in the next section. 

The cells $U_i'$ may have a variety of shapes; in order to apply induction on scale, we would like to put each of them  inside a smaller ball of radius $\frac{R}{d}$. To do so, it suffices to multiply $P$ by another polynomial $G$ of degree $6d$, and consider the cells cut off by the zero set of $P\cdot G$. More precisely, for $k=1,2,3$, let  $G_k = \Pi_{j=-d}^{d} (x_k -a_k -\frac{Rj}{d})$ such that $\{x_k=a_k\} \cap B_R\neq \emptyset$. The degree of $G_k$ is $\sim d$. Let $G = G_1 \cdot G_2 \cdot G_3$ be the product  and $Q\coloneqq P\cdot G$ be our new partitioning polynomial, then we have a new decomposition of $B_R$, 
\[ B_R\setminus Z(Q) =\bigsqcup_i O_i'.\]
The zero set $Z(Q)$ decomposes $B_R$ into at most $O(d^3)$ cells $O_i'$ by the Milnor-Thom Theorem \cite{Milnor, Thom}. Note that  in passing to the cells $O_i'$, the equidistribution property of the $BL^p$--norm over the cells may be lost; nevertheless, it is more or less kept in the cellular case, as we describe below. 

The general idea is to apply induction on the radius in each cell, and sum up the contributions of $Ef$ over the cells. To efficiently sum up the contributions, we need to understand how each wave packet $Ef_{\theta,v}$ interacts with the cells. By Lemma~\ref{lem: wpsupp}, a wave packet $Ef_{\theta,v}$ has negligible contribution to a cell $O_i'$ if its essential support $T_{\theta,v}$ does not intersect $O_i'$. To analyze how $T_{\theta,v}$ intersects a cell $O_i'$, we need to shrink $O_i'$ further. We define the wall $W$ as the $R^{1/2+\delta}$--neighborhood of $Z(Q)$ in $B_R$ and define the new cells as $O_i \coloneqq  O_i'\setminus W$. 

To summarize, we decompose $B_R= W\sqcup (\sqcup_i  O_i)$ and we have the following inequalities: 
\begin{equation}\label{eq: alcell}
\|Ef\|_{BL^p(B_R)}^p \lesssim \|Ef\|_{BL^p(W)}^p + \sum_i \|Ef\|_{BL^p(O_i)}^p
\end{equation}
and 
\begin{equation}\label{eq: cell}
\|Ef\|_{BL^p(O_i)}^p \lesssim d^{-3} \|Ef\|_{BL^p(B_R)}^P. 
\end{equation}

\textbf{The cellular case.} We are in the cellular case if $\|Ef\|_{BL^p(B_R)}^p \lesssim \sum_i \|Ef\|_{BL^p(O_i)}^p. $

Since the wave packets $Ef_{\theta,v}$ with $T_{\theta,v}\cap O_i=\emptyset$ have negligible contribution to $\|Ef\|_{BL^p(O_i)}$ in the sense that 
\[ \|\sum_{T_{\theta,v}\cap O_i=\emptyset} Ef_{\theta,v}\|_{BL^p (O_i)} = \mathit{RapDec}(R) \|f\|_{L^2}, \]
it is helpful to define 
\[ f_{O_i} : = \sum_{T_{\theta,v} \cap O_i \neq \emptyset } f_{\theta,v},\]
so  $\|Ef\|_{BL^p(O_i)} = \|Ef_{O_i}\|_{BL^p(O_i)} +\mathit{RapDec}(R)\|f\|_{L^2}.$ 

If a tube $T_{\theta,v}\cap O_i\neq \emptyset$, then the core line of $T_{\theta,v}$ must intersect $O_i'$. If a line $l \nsubseteq Z(Q)$, then $l$ intersects $Z(Q)$ $\lesssim d$ times. Hence, each $T_{\theta,v}$ intersects $\lesssim d+1$ cells $O_i$. 
So we have 
\[
\sum_i \|f_{O_i}\|_{L^2}^2 =\sum_{i} \sum_{T_{\theta,v}\cap O_i\neq \emptyset} \|f_{\theta,v}\|_{L^2}^2 \lesssim d \sum_{\theta,v} \|f_{\theta,v}\|_{L^2}^2 \lesssim d \|f\|_{L^2}^2.
\]
Since there are $\gtrsim d^3$ cells,  for $99 \%$ of the cells, 
\begin{equation}\label{eq: celll2}
\|f_{O_i}\|_{L^2} \lesssim d^{-1} \|f\|_{L^2}.
\end{equation}
In the cellular case, we have 
\[ \|Ef\|_{BL^p(B_R)}^p \lesssim \sum_i \|Ef\|_{BL^p(O_i)}^p. \]
Combining with inequality~\eqref{eq: cell}, we derive that for $\gtrsim d^3$ cells $O_i$, 
\begin{equation}\label{eq: popcell}
\|Ef\|_{BL^p(B_R)}^p \lesssim d^3 \|Ef_{O_i}\|_{BL^p(O_i)}^p + \mathit{RapDec}(R) \|f\|_{L^2}^p.
\end{equation}

Now we are ready to use Theorem~\ref{thm: guth} at scale $R/d$ as an inductive step for the functions $Ef_{O_i}$, which yields
\[ \|Ef_{O_i}\|_{BL^p(O_i)} \leq C_{\epsilon} (\frac{R}{d})^{\epsilon} \|f_{O_i}\|_{L^2}^{12/13} \underset{d(\tau)= (\frac{R}{d})^{-1/2}}{\max} \|f_{O_i, \tau} \|_{L^2_{avg}(\tau)}^{1/13}, \]
where $\tau$ denotes a cap of radius $(\frac{R}{d})^{-1/2}$. 

To estimate $\underset{ d(\tau)=(\frac{R}{d})^{-1/2}}{\max} \|f_{O_i, \tau}\|_{L^2_{avg}(\tau)}$, we shall apply the following two lemmas.

\begin{lemma}\label{lem: comparel2}
	Let $\tau$ be a cap of radius $r^{-1/2}$ with $r<R$, then 
	\[ \|f_{\tau}\|_{L^2_{avg}(\tau)} \lesssim \underset{\theta\subset \tau, d(\theta)= R^{-1/2}}{\max} \|f_{\theta}\|_{L^2_{avg}(\theta)}.\]
\end{lemma}
\begin{proof} Note that
	\begin{align*}
	\|f_{\tau}\|_{L^2_{avg}(\tau)}^2 &\lesssim \underset{\theta\subset \tau, d(\theta)=R^{-1/2}}{\mathit{Avg}} \|f_{\theta}\|_{L^2_{avg}(\theta)}^2\\
	&\lesssim \underset{\theta\subset \tau, d(\theta)=R^{-1/2}}{\max} \|f_{\theta}\|_{L^2_{avg}(\theta)}^2. 
	\end{align*}
\end{proof}
\begin{definition}\label{def: supT}
		Let $\mathbb{T}$ be a collection of tubes $T_{\theta,v}$, define 
		$$f^{\mathbb{T}} \coloneqq \sum_{T_{\theta,v}\in \mathbb{T}} f_{\theta,v}.$$
\end{definition}
\begin{lemma}\label{lem: fewerwp}
	$\|f^{\mathbb{T}}\|_{L^2} \lesssim \|f\|_{L^2}.$
\end{lemma}
\begin{proof}
	By inequality~\eqref{eq: l2orthogonality}, 
	\begin{align*}
	\|f^{\mathbb{T}}\|_{L^2}^2 \lesssim \sum_{T_{\theta,v}\in \mathbb{T}} \|f_{\theta,v}\|_{L^2}^2 
	\lesssim \sum_{T_{\theta,v}} \|f_{\theta,v}\|_{L^2}^2 \sim \|f\|_{L^2}^2. 
	\end{align*}
\end{proof}

We apply Lemma~\ref{lem: comparel2} to $f_{O_i, \tau}$ and then Lemma~\ref{lem: fewerwp} to the function $f_{O_i, \theta}$ to obtain:
\[ \underset{d(\tau) = (\frac{R}{d})^{-1/2} }{\max}\|f_{i, \tau}\|_{L^2_{avg}(\tau)} \lesssim  \underset{d(\theta)= R^{-1/2}}{\max} \|f_{O_i, \theta}\|_{L^2_{avg}(\theta)} \lesssim \underset{d(\theta)=R^{-1/2}}{\max} \|f_{\theta}\|_{L^2_{avg}(\theta)}. \]

We use a cell $O_i$ satisfying both inequality~\eqref{eq: celll2} and inequality~\eqref{eq: popcell} to close the induction,
\begin{align*}
\|Ef\|_{BL^p(B_R)}^p & \lesssim d^3 \|Ef_{O_i}\|_{BL^p(O_i)}^p\\
&\lesssim d^3 C_{\epsilon}^p (\frac{R}{d})^{\epsilon p} \|f_{O_i}\|_{L^2}^{\frac{12p}{13}} \underset{ d(\tau)=(\frac{R}{d})^{-1/2}}{\max} \|f_{O_i, \tau}\|_{L^2_{avg}(\tau)}^{\frac{p}{13}}\\
&\lesssim d^{3-\frac{12p}{13}} C_{\epsilon}^p ( \frac{R}{d})^{\epsilon p} \|f\|_{L^2}^{\frac{12p}{13}} \underset{d(\theta)=R^{-1/2}}{\max} \|f_{\theta}\|_{L^2_{avg}(\theta)}^{\frac{p}{13}}\\
&\leq C_{\epsilon}^p R^{\epsilon p} \|f\|_{L^2}^{\frac{12 p}{13}} \underset{d(\theta)=R^{-1/2}}{\max} \|f_{\theta}\|_{L^2_{avg}(\theta)}^{\frac{p}{13}},
\end{align*}
where the last inequality holds for any  $p\geq 13/4$ and any sufficiently large $d$ so that $d^{\epsilon p}\gtrsim 1$. 

\vspace{5pt}
\textbf{The algebraic case.} If we are not in the cellular case, then by \eqref{eq: alcell}, 
	$\|Ef\|_{BL^p(B_R)} \lesssim \|Ef\|_{BL^p(W)}.$ This is called the algebraic case because the $BL^p$--norm of $Ef$ is concentrated on the neighborhood of an algebraic surface. Similar to the previous argument concerning $f_{O_i}$, only wave packets $Ef_{\theta,v}$ whose essential support $T_{\theta,v}$ intersect $W$ contribute  to $\|Ef\|_{BL^p(W)}$. Depending on how they intersect, we identify a tangential part, which consists of wave packets tangential to $W$, and a transversal part, which consists of wave packets intersecting $W$ transversely. 

We give the precise definition of the tangential tubes and the transversal tubes as follows. We cover $W$ with finitely overlapping balls $B_k$ of radius $\rho \coloneqq R^{1-\delta}$. For each $B_k$, let 
\begin{equation}\label{eq: Sk}
S_k: = W\cap B_k.
\end{equation}
  We define $\mathbb{T}_{S_k}$ to be the collection of the tangential tubes for $S_k$ and $\mathbb{T}_{S_k, \mathit{trans}}$ to be  the collection of transversal tubes for $S_k$. Their definitions are the same as in Guth's paper \cite{Guth1}. 

\begin{definition}[{\cite[Definition 3.3]{Guth1}}]\label{def: tang} 
	$\mathbb{T}_{S_k}$ is the set of all tubes $T$ obeying the following two conditions: 
	\begin{itemize}
		\item $T\cap S_k \neq \emptyset$.
		\item If $z$ is any non-singular point of $Z(P)$ lying in $2B_k \cap 10 T$, then 
		\[ \mathit{Angle} (v(T), T_z Z(P)) \leq R^{-1/2+2\delta},\]
		where $v(T)$ is the unit vector in the direction of the tube $T$.
	\end{itemize}

\end{definition}

\begin{definition}[{\cite[Definition 3.4]{Guth1}}]\label{def: trans}
	$\mathbb{T}_{S_k, \mathit{trans}}$ is the set of all $T$ obeying the following two conditions: 
	\begin{itemize}
		\item $T\cap S_k \neq \emptyset$.
		\item There exists a non-singular point $z$ of $Z(P)$ lying in $2B_k \cap 10T$, such that 
		\[ \mathit{Angle}(v(T), T_z Z(P)) > R^{-1/2+2\delta}.\]
	\end{itemize}
\end{definition}
For each $S_k$, define 
\[ f_{S_k} \coloneqq \sum_{T_{\theta,v} \in \mathbb{T}_{S_k}}f_{\theta,v}, \text{ and } f_{S_k, \mathit{trans}} \coloneqq \sum_{T_{\theta,v} \in \mathbb{T}_{S_k, \mathit{trans}}} f_{\theta,v}. \]

We estimate the algebraic part by 
\begin{align}\label{eq: tangtrans}
\|Ef\|_{BL^p(W)}^p &\lesssim \sum_{S_k} \|Ef_{S_k} + Ef_{S_k, \mathit{trans}} \|_{BL^p (S_k)}^p + \mathit{RapDec}(R)\|f\|_{L^2}^p \\
&\lesssim \sum_{S_k} \|Ef_{S_k}\|_{BL^p(S_k)}^p + \sum_{S_k} \|Ef_{S_k, \mathit{trans}} \|_{BL^P(S_k)}^p + \mathit{RapDec}(R)\|f\|_{L^2}^p. 
\end{align}
The last inequality is due to  \eqref{eq: triangle}. 

Here is a remark regarding the parameter $A$ in the $BL^p$--norm.  Now we think of the argument as a recursion rather than an induction.  Inequality~\eqref{eq: triangle}  is only used in the algebraic cases.  This recursion only involves  $\leq N \coloneqq \lceil \frac{\log \delta}{\log (1-\delta)} \rceil$ steps in the algebraic cases because otherwise the radius would be $\leq R^{(1-\delta)^N}\leq R^{\delta}$ and we could apply Lemma~\ref{lem: tao}  using $R^{3\delta/p} <  R^{\epsilon^3}$. Hence, the parameter $A$ is reduced for  $\leq N$ times, which is a constant independent of $R$.  It suffices to set  $A\geq N^2$, and we can omit $A$ from our notation. 

The algebraic case can be further divided  into the tangential case and the transversal case. For the tangential case we are going to estimate $\|Ef_{S_k}\|_{BL^p(S_k)}$ directly,  and for the transversal case we will apply induction on scale as in the cellular case. 

\vspace{5pt}
\textbf{The tangential sub-case.}
We are in the tangential case if $\sum_{S_k} \|Ef_{S_k}\|_{BL^p(S_k)}^p \gtrsim \|Ef\|_{BL^p(W)}^p \gtrsim \|Ef\|_{BL^p(B_R)}^p.$ The function
$f_{S_k}$  is concentrated on wave packets from $\mathbb{T}_{S_k}$,  which heuristically implies that the support of $f_{S_k}$ is small: 

\begin{lemma}[{\cite[Lemma 4.9]{Guth1}}]\label{lem: geometric}
	The support $\supp f_{S_k}$ lies in a union of $\lesssim R^{1/2+O(\delta)}$ caps $\theta$ of radius $R^{-1/2}$. 
\end{lemma}

We shall apply this crucial geometric lemma to efficiently bound $\|f_{S_k}\|_{L^2}$.
 Indeed, by Lemma~\ref{lem: geometric} and Lemma~\ref{lem: fewerwp}, 
 \begin{equation}\label{eq: l2avg}
 \|f_{S_k}\|_{L^2}^2 \lesssim R^{-1/2+O(\delta)}  \underset{d(\theta)=R^{-1/2}}{\max}\|f_{S_k,\theta}\|_{L^2_{avg}(\theta)}^2 \lesssim R^{-1/2+O(\delta)} \underset{d(\theta)=R^{-1/2}}{\max} \|f_{\theta}\|_{L^2_{avg}(\theta)}^2. 
 \end{equation}
 
 By Lemma~\ref{lem: fewerwp}, we also have 
 \begin{equation}\label{eq: l2avg2}
 \|f_{S_k}\|_{L^2} \lesssim \|f\|_{L^2}.
 \end{equation}
  When $13/4 \leq p\leq 10/3$, we combine Lemma~\ref{lem: tao} and \eqref{eq: l2avg}  \eqref{eq: l2avg2}  to obtain
 \begin{align*}
 \|Ef_{S_k}\|_{BL^p(B_R)}^p &\lesssim \rho^{\frac{5}{2}-\frac{3p}{4}-\frac{p}{52} +\epsilon^{10}}R^{O(\delta)}  \|f\|_{L^2}^{\frac{12p}{13}} \underset{d(\theta)=R^{-1/2}}{\max} \|f_{\theta}\|_{L^2_{avg}(\theta)}^{p/13}\\&\lesssim R^{O(\delta)} \|f\|_{L^2}^{\frac{12p}{13}} \underset{d(\theta)=R^{-1/2}}{\max}  \|f_{\theta}\|_{L^2_{avg}(\theta)}^{p/13}. 
 \end{align*}
 
 Since we can cover $B_R$ by  $\lesssim R^{3\delta}$ balls $B_k$, we have 
 \begin{align*}
\|Ef\|_{BL^p(B_R)}^p &\lesssim \sum_{S_k} \|Ef_{S_k}\|_{BL^p(S_k)}^p\\
&\lesssim R^{O(\delta)} \|f\|_{L^2}^{12p/13} \underset{d(\theta)=R^{-1/2}}{\max} \|f_{\theta}\|_{L^2_{avg}(\theta)}^{p/13}\\
&\leq C_{\epsilon} R^{\epsilon} \|f\|_{L^2}^{12p/13} \underset{d(\theta)=R^{-1/2}}{\max} \|f_{\theta}\|_{L^2_{avg}(\theta)}^{p/13}. 
\end{align*}
Since $\delta \leq \epsilon^3$, when $\epsilon$ is small enough, $R^{O(\delta)} \leq R^{\epsilon}$. Then we choose $C_{\epsilon}$ large enough to dominate the implicit constants appeared in $\lesssim$.  This completes the estimate for the tangential case. 

 \vspace{5pt}
 \textbf{The transversal sub-case.}
 We are in the transversal case when $$\|Ef\|_{BL^p(B_R)}^p \lesssim \sum_{S_k} \|Ef_{S_k, \mathit{trans}}\|_{BL^p(S_k)}^p.$$ By inequalities~\eqref{eq: alcell} \eqref{eq: tangtrans}, if $\|Ef\|_{BL^p(B_R)}$ is not dominated by the cellular part or the tangential part, then it must be dominated by the transversal part. To deal with the transversal case, we shall apply induction using Theorem~\ref{thm: guth} on $Ef_{S_k, \mathit{trans}}$ at scale $\rho$, which is similar to the cellular case. 
 
  After induction in each $S_k$, we need to sum up $\|f_{S_{k}, \mathit{trans}}\|_{L^2}$, which requires the following Lemma~\ref{lem: trans} to obtain a similar version of inequality~\eqref{eq: celll2}.
 
 \begin{lemma}[{\cite[Lemma 3.5]{Guth1}}]\label{lem: trans}
 	Each tube $T_{\theta,v}$ belongs to $\lesssim \mathit{Poly}(d)$ different sets $\mathbb{T}_{S_k, \mathit{trans}}$. Here $\mathit{Poly}(d)$ means a quantity bounded by a  polynomial of $d$  whose  degree  and coefficients  are  constants independent of $d$. 
 \end{lemma}

By Lemma~\ref{lem: trans}, we have 
\begin{equation}\label{eq: trans}
\sum_{S_k} \|f_{S_k, \mathit{trans}}\|_{L^2}^2 \lesssim \mathit{Poly}(d) \|f\|_{L^2}^2.
\end{equation}

We apply induction on scale and then sum over $S_k$: 
\begin{align}\label{eq: esttrans}
\|Ef\|_{BL^p(B_R)}^p &\lesssim \sum_{S_k} \|Ef_{S_k, \mathit{trans}}\|_{BL^p(S_k)}^p \\
&\lesssim C_{\epsilon}^p \rho^{\epsilon p} \sum_{S_k} \|f_{S_k, \mathit{trans}}\|_{L^2}^{\frac{12p}{13}} \underset{d(\tau)=\rho^{-1/2}}{\max} \|f_{S_k, \mathit{trans}, \tau}\|_{L^2_{avg}(\tau)}^{\frac{p}{13}}.
\end{align}
To estimate $\|f_{S_k, \mathit{trans}, \tau}\|_{L^2_{avg}(\tau)}$, we apply Lemma~\ref{lem: comparel2} and then Lemma~\ref{lem: fewerwp} on $\|f_{S_k, \mathit{trans}, \theta}\|_{L^2}$:
\[ \|f_{S_k, \mathit{trans}, \tau}\|_{L^2_{avg}(\tau)}\lesssim \underset{ \theta\subset \tau, d(\theta)=R^{-1/2}}{\max} \|f_{S_k, \mathit{trans}, \theta}\|_{L^2_{avg}(\theta)}\lesssim \underset{d(\theta)=R^{-1/2}}{\max} \|f_{\theta}\|_{L^2_{avg}(\theta)}.\]

Then we sum up $\|f_{S_k, \mathit{trans}}\|_{L^2}$ over  all $S_k$ using \eqref{eq: esttrans} and  the fact that $\|\cdot \|_{l^q} \leq \|\cdot \|_{l^2}$ when $q\geq 2$: 
\begin{align*}
\|Ef\|_{BL^p(B_R)}^p &\lesssim C_{\epsilon}^p \rho^{\epsilon p} ( \sum_{S_k} \|f_{S_k, \mathit{trans}}\|_{L^2}^{\frac{12p}{13}}) \underset{d(\theta)=R^{-1/2}}{\max} \|f_{\theta}\|_{L^2_{avg}(\theta)}^{\frac{p}{13}}\\
& \lesssim C_{\epsilon}^p \rho^{\epsilon p} (\sum_{S_k} \|f_{S_k, \mathit{trans}}\|_{L^2}^2)^{\frac{6p}{13}} \underset{d(\theta)=R^{-1/2}}{\max} \|f_{\theta}\|_{L^2_{avg}(\theta)}^{\frac{p}{13}} \\
&\lesssim C_{\epsilon}^p R^{\epsilon p} R^{-\delta \epsilon p}  \mathit{Poly}(d) \|f\|_{L^2}^{\frac{12p}{13}} \underset{d(\theta)=R^{-1/2}}{\max} \|f_{\theta}\|_{L^2_{avg}(\theta)}^{\frac{p}{13}}\\
&\leq C_{\epsilon}^p R^{\epsilon p}  \|f\|_{L^2}^{\frac{12p}{13}} \underset{d(\theta)=R^{-1/2}}{\max} \|f_{\theta}\|_{L^2_{avg}(\theta)}^{\frac{p}{13}}.
\end{align*}
Since $d\leq \log R$, when $R$ is sufficiently large, $R^{\delta \epsilon p} $ is larger than $\mathit{Poly}(d)$ multiplying the implicit constant in $\lesssim$. 
This completes the induction step in the transversal case and the proof of Theorem~\ref{thm: guth}. 
\subsection{Two ends argument}\label{subsection: te}
In this subsection, we recall Wolff's two ends argument following Tao's paper \cite[Section 5]{Tao3}. We cover $B_R$ with balls $B_k$ of radius $\rho= R^{1-\epsilon_0}$ for some $\epsilon_0$ with $\delta/\epsilon \leq \epsilon_0 \leq \epsilon/100$.  The idea is to assign each tube $T_{\theta,v}$ a few  balls $B_k$ one can ``exclude'' via the inductive hypothesis (see Lemma~\ref{lem: telocal}). This assignment, which we denote by $B_k \sim T_{\theta,v}$, will be defined explicitly  in Section~\ref{section: sim}.   For now, it suffices to know that the relation $\sim$ satisfies the following property: 
\begin{itemize}
\item 	for each $T_{\theta,v}$, the number of $B_k$ with $B_k\sim T_{\theta,v}$ is $\lesssim_{\delta} 1$. 
\end{itemize}
We choose $R$ large enough such that $O_{\delta}(1) \leq R^{\delta^2}$. 

For each $B_k$, we define $\mathbb{T}_{k}^{\sim} =\{ (\theta,v): B_k\sim T_{\theta,v}\}$ and  $\mathbb{T}_{k}^{\nsim} =\{ (\theta,v): B_k\nsim T_{\theta,v}\}$. 
We define 
 the local part of $f$ as $$f_k^{\sim}  \coloneqq \sum_{(\theta,v)\in \mathbb{T}_{k}^{\sim}} f_{\theta,v},$$ 
 and the global part of $f$ as $$f^{\nsim}_k \coloneqq \sum_{(\theta,v)\in \mathbb{T}_{k}^{\nsim}} f_{\theta,v}.$$  
 For each $T_{\theta,v}$,  there are  $\lesssim_{\delta} 1$ balls $B_k$ such that $B_k\sim T_{\theta,v}$, so 
\begin{equation}\label{eq: telocal}
\sum_{B_k} \|f_k^{\sim}\|_{L^2}^2 \lesssim_{\delta} \|f\|_{L^2}^2.
\end{equation}
By the triangle inequality, 
\begin{equation}\label{eq: snsim}
\|Ef\|_{BL^p(B_R)}^p \lesssim \sum_{B_k} \|Ef^{\sim}_k \|_{BL^p(B_k)}^p + \sum_{B_k} \|Ef^{\nsim}_{k}\|_{BL^p(B_k)}^p +\mathit{RapDec}(R)\|f\|_{L^2}^p.
\end{equation}

The local part is relatively easy and can be controlled using  induction on scale: 
\begin{lemma}\label{lem: telocal}
	If Theorem~\ref{thm: induction2} holds at scale $\rho=R^{1-\epsilon_0}$ for each  $f_k^{\nsim}$ and $$\sum_{B_k}\|Ef_k^{\sim}\|_{BL^p(B_k)}^p \gtrsim R^{-\delta} \|Ef\|_{BL^p(B_R)}^p,$$ then Theorem~\ref{thm: induction2} holds at scale $R$ for $f$. 
\end{lemma}
\begin{proof}
	
	Since  Theorem~\ref{thm: induction2} holds at scale $\rho=R^{1-\epsilon_0}$ for the functions $f^{\sim}_k$,  we have 
\begin{align*}
\|Ef\|_{BL^p(B_R)}^p & \lesssim R^{\delta} \sum_{B_k} \|Ef_k^{\sim}\|_{BL^p(B_k)}^p \\
&\lesssim R^{\delta} C_{\epsilon}^p \rho^{\epsilon p }  \sum_{B_k} \|f_k\|_{L^2}^2 \underset{d(\theta')=\rho^{-1/2}}{\max} \|f_{k,\theta'}\|_{L^2_{avg}(\theta')}^{p-2}\\
&\lesssim C_{\epsilon}^p R^{\epsilon p} R^{-p \epsilon \epsilon_0 +O(\delta)} \|f\|_{L^2}^2  \underset{d(\theta)=R^{-1/2}}{\max} \|f_{\theta}\|_{L^2_{avg}(\theta)}^{p-2}.
\end{align*}
The last inequality follows from  \eqref{eq: telocal}, Lemma~\ref{lem: comparel2} and Lemma~\ref{lem: fewerwp}. Since $\delta/\epsilon \lesssim \epsilon_0$,  when $R$ is sufficiently large, we have  $C R^{-\epsilon \epsilon_0 +O(\delta)} <1$, where $C$ is the implicit constant  in $\lesssim$. 
\end{proof}

The estimate of the global part is more involved, which  will occupy the rest of the paper.  The   tools that are used in this estimate  will be developed in the following  sections.

In fact,  we shall combine the two ends argument with the polynomial partitioning techniques to obtain Theorem~\ref{thm: induction2};  see Section~\ref{section: proofmain} for more details. 
\section{Polynomial structure lemmas}\label{section: poly}

In this section, we explain how to decompose a function concentrated on wave packets from a set of tubes $T_{\theta,v}$ using polynomial partitioning iteratively.  This is carried out in two steps. For the first step, we use Lemma~\ref{lem: recursion} to  decompose $\|Ef\|_{BL^p(B_R)}^p$  by iteratively applying Guth's polynomial partitioning theorem. For the second step, we show in Lemma~\ref{lem: properties} how to use the partitioning polynomials given by the first step to decompose a function concentrated on wave packets from a set $\mathbb{T}$ of tubes $T_{\theta,v}$. 

To state Lemma~\ref{lem: recursion}, let us introduce the notions of $r$--fat surfaces and tree structures. 
\begin{definition}\label{def: fatsurf}
	Fix a large integer $d\sim \log R$ and some  $r$ with $R^{\delta}\leq r\leq R$. Let $B_{r^{1-\delta}}$ be a ball of radius $r^{1-\delta}$ and $Z_S$ be a union of smooth algebraic surfaces with  $\deg Z_S \leq d$. 
	We define a  fat $r$--surface $S$  to be  $$S\coloneqq N_{r^{1/2+\delta} }Z_S\cap B_{r^{1-\delta}}.$$
\end{definition}
\begin{remark}
	$S_k$ defined in \eqref{eq: Sk} is an example of a fat $R$--surface. 
\end{remark}

 In Lemma~\ref{lem: recursion}, we use a tree structure to describe the hierarchy of partitioning surfaces.  A tree structure $\mathcal{O}$ consists of tree elements, which we call nodes. Below are some terminologies used in the tree structure   in Lemma~\ref{lem: recursion}. 
\renewcommand\labelitemi{\tiny$\bullet$}
\begin{itemize}
	\item Root: the top node in a tree, the prime ancestor. 
	\item Child: a node directly connected to another node when moving away from the root, an immediate descendant. 
	\item Parent: the converse notion of a child, an immediate ancestor.
	\item Descendant: a node reachable by repeated proceeding from parent to child. 
	\item Ancestor: a node reachable by repeated proceeding from child to parent.
	\item Leaf: a node with no children.
	\item Edge: the connection between one node and another.
	\item Path: a sequence of nodes and edges connecting a node with a descendant.
	\item Depth: the number of edges along the shortest path between a node and the root.\
	\item Height: the number of edges on the longest path between a node and a descendant leaf.  The height of tree is the height of the root node in the tree. 
\end{itemize}


\begin{lemma}\label{lem: recursion}
	If $f$ is supported in the unit disk, then there exists a tree structure $\mathcal{O}$ of height $J\lesssim \frac{\log R}{\log \log R}$ satisfying the following properties. 
	\begin{enumerate}
		\item \label{en: str1} The root of $\mathcal{O}$ is $O_0 \coloneqq  B_{R_0}$ with $R_0\coloneqq R$.
		\item  For each $0\leq j\leq J-1$,  the children of a node $O_j$  of depth $j$ are  some subsets $O_{j+1}$ of $O_j$, and each $O_{j+1}$ lies in some ball $B_{R_{j+1}}$ of radius $R_{j+1}\leq R_j/d$. Here the radii $R_{j+1}$ are the same for all nodes  $O_{j+1}$ of depth $j+1$.  Moreover,  $R_J\leq R^{\delta}, R_{J-1}> R^{\delta}.$
		\item There exist a number $n\leq \delta^{-2}$ and  indices $0=j_0 < j_1 < \dots < j_n\leq J$  such that each node $O_{j_t}$ of depth $j_t$ with  $1\leq t\leq n$ is a fat $R_{j_t-1}$--surface $S_t$,  and   $R_{j_t} = R_{j_t-1}^{1-\delta}$. Let $\mathcal{S}_t$ denote the set of all  fat $R_{j_t-1}$--surfaces in $\mathcal{O}$. In other words, $\mathcal{S}_t$ is  the set of  all $j_t$--th nodes of $\mathcal{O}$. 
	\item Define the function $T(j): = \max \{ t: j_t\leq j\}$.   For each $1\leq j \leq J$, the following estimates hold:
	\begin{equation} \label{eq: indd}
	\|Ef\|_{BL^p(O_j)}^p \lesssim d^{-3(j-T(j))}  \|Ef\|_{BL^p(B_R)}^p, 
	\end{equation}
	and 
	\begin{equation}\label{eq: indc}
	 \|Ef\|_{BL^p(B_R)}^p \leq 2^{j} (\log R)^{T(j)} \sum_{O_j\in \mathcal{O}} \|Ef\|_{BL^p(O_j)}^p,
	\end{equation}
	where the sum $\sum_{O_j\in \mathcal{O}}$ is  over all nodes $O_j$ of depth $j$ in $\mathcal{O}$.
	\item If $j\neq j_t$ for any $1\leq t\leq n$, then for each node $O_{j-1}$ of depth $j-1$,  a tube $T_{\tau_j, w_j}$ of length $R_{j-1}$ and radius $R_{j-1}^{1/2+\delta}$ intersects  $\lesssim d$  children $O_j$ of $O_{j-1}$. 
\end{enumerate}

\end{lemma}
\begin{proof}
	Assume that we have defined $\mathcal{O}$ up to depth $j-1$ with the radius $R_{j-1}>R^{\delta}$ and we have verified inequalities \eqref{eq: indd} and \eqref{eq: indc} for all $0\leq j' \leq j-1$.  The base case is  when  $j-1=0$.   In  the statement of  Lemma~\ref{lem: recursion},  we have defined $O_0= B_R$, $R_0=R$ and $j_0=0$.  So  \eqref{eq: indd} and \eqref{eq: indc} hold trivially for $j-1=0$ and $T(0)=0$. 
	
For each node $O_{j-1}$ of depth $j-1$, we now  define its children. 
We apply polynomial partitioning on $\|Ef\|_{BL^p(O_{j-1})}^p$ as in Subsection~\ref{subsection: poly}. More precisely,  the parameter $R$ in Subsection~\ref{subsection: poly} becomes $R_{j-1}$ and the degree $d\sim \log R$.  Let $Z_{j}$ be the zero set of the degree $d$ partitioning polynomial and set  $W_{j}  \coloneqq N_{R_{j-1}^{1/2+\delta}} Z_j\cap B_{R_{j-1}}$. Then  we have
	\begin{equation}\label{eq: decca}
	\|Ef\|_{BL^p(O_{j-1})}^p \leq \|Ef\|_{BL^p(W_j)}^p + \|Ef\|_{BL^p(O_{j-1}\setminus W_j)}^p.
	\end{equation}

	\textbf{The Cellular case.}
	 We are in the cellular case if
	\begin{equation}\label{eq: inddj}
	\sum_{O_{j-1}\in \mathcal{O}} \|Ef\|_{BL^p(O_{j-1})}^p \leq 2\sum_{O_{j-1}\in \mathcal{O}}\|Ef\|_{BL^p(O_{j-1}\setminus W_j)}^p. 
	\end{equation}
	
	For each node $O_{j-1}\in \mathcal{O}$ of depth $j-1$, the set $O_{j-1}\setminus W_j$ is decomposed into $\lesssim d^3$ cells  $O_j$: 
     $$O_{j-1}\setminus W_j =\bigsqcup O_j.$$
     The children of $O_{j-1}$ then consist of all cells $O_j$ that appear in this decomposition. 
     Each tube $T_{\tau_j, w_j}$ of length $R_{j-1}$ and  radius $R_{j-1}^{1/2+\delta}$ intersects $\lesssim d$ cells $O_j$ in $O_{j-1}$. So we have verified (5) in the statement of Lemma~\ref{lem: recursion} for $j$. 
	For each cell $O_j$, \eqref{eq: cell}  implies that
\begin{equation}\label{eq: cellpoly}
\|Ef\|_{BL^p(O_j)}^p \lesssim d^{-3} \|Ef\|_{BL^p(O_{j-1})}^p.
\end{equation}
	 Define $R_j \coloneqq R_{j-1}/d$. 
Then \eqref{eq: indd} follows from  \eqref{eq: cellpoly} and  \eqref{eq: indc} follows from \eqref{eq: inddj}. If $R_j\leq R^{\delta}$, then we stop and define $J=j$; otherwise we proceed to the next step for each node $O_j$ of depth $j$. 
	
	\textbf{The Algebraic case.}
	If \eqref{eq: inddj} does not hold, then by \eqref{eq: decca}, 
	\begin{equation}\label{eq: indaj}
	\sum_{O_{j-1}\in \mathcal{O}}\|Ef\|_{BL^p(O_{j-1})}^p \leq 2 \sum_{O_{j-1}\in \mathcal{O}}\|Ef\|_{BL^p(W_j)}^p.
	\end{equation}
	Define  $j_{ T(j-1)+1} : =j$, so $T(j)=T(j-1)+1$. 
	 We cover $W_j$ with finitely overlapping balls $B_k$ of radius $R_j \coloneqq R_{j-1}^{1-\delta}$.  
	 Since $R_{j-1}\geq R^{\delta}$,  $R_{j} = R_{j-1}/ R_{j-1}^{\delta} \leq  R_{j-1}/d$.  We define  the fat $R_{j-1}$--surface $S_{T(j)} \coloneqq W_j \cap B_{R_j}$. 
	We sort the fat $R_{j-1}$--surfaces $S_{T(j)}$ according to $\|Ef\|_{BL^p(S_{T(j)})}$.
	Since for each $O_{j-1}$ the number of $S_{T(j)} \subset O_{j-1}$  is $\lesssim R_{j-1}^{3\delta}$,  the sum of $\|Ef\|_{BL^p(S_{T(j)})}^p$ over those fat $R_{j}$--surfaces $S_{T(j)}$ with $$\|Ef\|_{BL^p(S_{T(j)})}\leq R_{j-1}^{-1} \|Ef\|_{BL^p(O_{j-1})}$$ is $\lesssim R_{j-1}^{-p+3\delta} \|Ef\|_{BL^p(O_{j-1})}^p$. Therefore we can discard these $S_{T(j)}$  from consideration. 
	There exists a dyadic number $\lambda$ with  $R_{j-1}^{-1} \leq \lambda \leq 1$  such that the corresponding set $$\mathcal{S}_{T(j)} \coloneqq \{ S_{T(j)}: \|Ef\|_{BL^p(S_{T(j)})} \sim \lambda \|Ef\|_{BL^p(O_{j-1})}, \text{ where }O_{j-1} \supset S_{T(j)} \}$$ satisfies
	\begin{equation}\label{eq: pg1}
	\sum_{O_{j-1}\in \mathcal{O}}\|Ef\|_{BL^p(O_{j-1})}^p\leq 2\log R_{j-1} \sum_{S_{T(j)}\in \mathcal{S}_{T(j)}} \|Ef\|_{BL^p(S_{T(j)})}^p.
	\end{equation}
	We define the children $O_j$ of $O_{j-1}$ as the fat $R_{j-1}$--surfaces $S_{T(j)}\in \mathcal{S}_{T(j)}$ that are contained in $O_{j-1}$.  If $R_j\leq R^{\delta}$, we stop and define $J=j$, $n=T(j)$; otherwise, we proceed to the next step for each $O_j$.  The estimate  \eqref{eq: indc} now follows from  \eqref{eq: indaj} and \eqref{eq: pg1} with $R_{j-1}\leq R$. Since $T(j)=T(j-1)+1$, inequality~\eqref{eq: indd} holds by the induction hypothesis  that \eqref{eq: indd} holds for $j-1$.
	
	\vspace{5pt}
	
	Finally, 
	since $R_j\leq R_{j-1}/d$ and $d\sim \log R$, we have $J\lesssim  \frac{\log R}{\log \log R}$.  Moreover,  for each number $1\leq t\leq n$,  we have $R_{j_t}= R_{j_t-1}^{1-\delta}$ and $R_{j_n-1}\geq R^{\delta}$; so   $n\leq \delta^{-2}$. 

\end{proof}

By the construction in Lemma~\ref{lem: recursion}, each leaf  in $\mathcal{O}$ has depth $J$. 
In the next corollary, we will refine the tree structure $\mathcal{O}$ by removing some leaves. For the rest of the paper, we consider only the refined tree structure defined in Corollary~\ref{cor: tree}.  To ease notation, write $O=O_J$. 

\begin{cor}\label{cor: tree}
	For the tree $\mathcal{O}$ defined in Lemma~\ref{lem: recursion}, we can remove some leaves of $\mathcal{O}$, such that 
	\begin{enumerate}
		\item for each  remaining  leaf $O \in \mathcal{O}$, 
		\begin{equation} \label{eq: leaf}
		\|Ef\|_{BL^p(O)}\sim \lambda_0\|Ef\|_{BL^p(B_R)};
		\end{equation}
		\item for the  remaining leaves in $ \mathcal{O}$, \begin{equation}
		\|Ef\|_{BL^p(B_R)} \lesssim R^{\delta} \sum_{O\in \mathcal{O}} \|Ef\|_{BL^p(O)}^p.
		\end{equation}
	\end{enumerate}
\end{cor}
\begin{proof}
			We sort the leaves $O\in \mathcal{O}$ according to $\|Ef\|_{BL^p(O)}$. Since the number of leaves is $\lesssim R^3$, it suffices to consider only the leaves $O$ with $\|Ef\|_{BL^p(O)} \geq R^{-4} \|Ef\|_{BL^p(B_R)}$. 
			 By pigeonholing and inequality~\eqref{eq: indc}, there exists  a dyadic number $\lambda_0$ with  $R^{-4 }  \leq \lambda_0 \leq 1$ such that if we keep only the leaves satisfying \eqref{eq: leaf},  then
	\begin{equation}
	\|Ef\|_{BL^p(B_R)}^p \leq 2^{J+4} (\log R)^{n +1}  \sum_{O\in \mathcal{O}} \|Ef\|_{BL^p(O)}^p.
	\end{equation}
	Finally,   $2^{J+4} (\log R)^{n+1}\lesssim R^{\delta}$ for sufficiently  large  $R$ because $J\lesssim \frac{\log R}{\log \log R}$ and $n\leq \delta^{-2}$. 
\end{proof}

\begin{definition}\label{def: tangsurf}
	Let $S$ be a fat $r$--surface.  Let $\mathbb{T}_S$ be the  set of all tubes $T$ of length $r$ and  radius $r^{1/2+\delta}$ obeying  the following two conditions: 
	\begin{itemize}
		\item $T\cap S\neq \emptyset$,
		\item If $z$ is any non-singular point of $Z_S$ lying in $2B_{r^{1-\delta}} \cap 10 T$ for the surface  $Z_S$ and the ball  $B_{r^{1-\delta}}$ as  in Definition~\ref{def: fatsurf}, then 
		\[ \mathit{Angle}(v(T), T_z Z_S)\leq r^{-1/2+2\delta}.\]
		Here $v(T)$ is the unit vector in the direction of the tube $T$.
	\end{itemize}
\end{definition}

\begin{definition}\label{def: transsurf}
		Let $S$ be a fat $r$--surface. Let $\mathbb{T}_{S,\mathit{trans}}$ be the  set of all tubes $T$ of length $r$ and radius $r^{1/2+\delta}$ obeying  the following two conditions: 
	\begin{itemize}
		\item $T\cap S\neq \emptyset$.
		\item If $z$ is any non-singular point of $Z_S$ lying in $2B_{r^{1-\delta}} \cap 10 T$ for the surface $Z_S$ and the ball $B_{r^{1-\delta}}$ as in Definition~\ref{def: fatsurf}, then 
		\[ \mathit{Angle}(v(T), T_z Z_S)\geq r^{-1/2+2\delta}.\]
	\end{itemize}
	\end{definition}
Note that for the  fat $R$--surfaces $S_k$,  Definition~\ref{def: tangsurf} and Definition~\ref{def: transsurf}  are the same as Definition~\ref{def: tang} and Definition~\ref{def: trans}.
Recall that for any set $\mathbb{T}$ of tubes $T_{\theta,v}$, we defined  $g^{\mathbb{T}}\coloneqq \sum_{T_{\theta,v}\in \mathbb{T}} g_{\theta,v}$. In particular, $f^{\nsim}_k$ fits into this definition with $g=f$ and  $\mathbb{T}=\mathbb{T}_k^{\nsim}$.
\begin{lemma}\label{lem: properties}
	Let $\mathcal{O}$ be the tree defined in Corollary~\ref{cor: tree}. For any function $g$ supported in the unit disk in $\mathbb{R}^2$, we perform the wave packet  decomposition  $g=\sum_{\theta,v} g_{\theta,v} +\mathit{RapDec}(R)$.  Given any set  $\mathbb{T}$ of tubes $T_{\theta,v}$,  we can define $g^{\mathbb{T}}_{O}$, $g^{\mathbb{T}}_{S_t}$ for each leaf $O \in \mathcal{O}$ and  each  $S_t\in \mathcal{S}_t$ with  $1\leq  t\leq n$ satisfying the following properties. 
	\begin{enumerate}[label=(\roman*)]
		\item Inside each leaf $O\in \mathcal{O}$, consider the ancestors of $O$: $S_1\supset\cdots \supset S_n \supset O$. Then
		\begin{equation}\label{eq: deccell}
Eg^{\mathbb{T}} = Eg^{\mathbb{T}}_O +\sum_{t=1}^n Eg^{\mathbb{T}}_{S_t} + \mathit{RapDec}(R)\|g\|_{L^2}.
\end{equation}
\item For each fat $R_{j_t-1}$--surfaces $S_t\in \mathcal{S}_t$, we define 
\begin{equation}\label{eq: defpi}
g^{\mathbb{T}}_{\Pi_{S_t}} \coloneqq \sum_{T_{\tau_{j_t}, w_{j_t}} \in \mathbb{T}_{S_t}} (\sum_{T_{\theta,v}\in \mathbb{T}, T_{\theta,v}\cap S_t\neq\emptyset} g_{\theta,v})_{\tau_{j_t}, w_{j_t}}.
\end{equation}
Consider the ancestors of $S_t$:  $S_1 \supset\cdots \supset S_t$. Then we have 
 \begin{equation}\label{eq: relpi}
 g^{\mathbb{T}}_{\Pi_{S_t}} = g^{\mathbb{T}}_{S_t} + \sum_{l=1}^{t-1} g^{\mathbb{T}}_{S_l, S_t} +\mathit{RapDec}(R)\|g\|_{L^2},
 \end{equation}
 where $g^{\mathbb{T}}_{S_l, S_t}  \coloneqq \sum_{T_{\tau_{j_t}, w_{j_t}}\in \mathbb{T}_{S_t}} (g^{\mathbb{T}}_{S_l})_{\tau_{j_t}, w_{j_t}}$. 
 \item The following $L^2$--estimates hold: 
 \begin{equation}\label{eq: strcell}
 \sum_{S_t\in \mathcal{S}_t} \|g^{\mathbb{T}}_{S_t}\|_{L^2}^2 \lesssim d^{j_t-t} \mathit{Poly}(d)^{t-1} R^{3\delta}\|g\|_{L^2}^2, 
 \end{equation}
 and 
 \begin{equation}\label{eq: strsurf}
 \sum_{O \in \mathcal{O}} \|g^{\mathbb{T}}_O\|_{L^2}^2 \lesssim d^{J-n} \mathit{Poly}(d)^n\||g\|_{L^2}^2.
 \end{equation}
 
	\end{enumerate}
\end{lemma}
Note that  $S_t$ is some node $O_{j_t}$ of depth $j_t$, but $g^{\mathbb{T}}_{S_t}$ and $g^{\mathbb{T}}_{O_{j_t}}$ represent different functions. Morally speaking, we use the subscript $S_t$ to indicate the tangential part and $O_{j_t}$ to indicate the transversal part.  Their precise definitions are given by  \eqref{eq: surfta} and \eqref{eq: surftr} below. 
\begin{proof}
	Define $g^{\mathbb{T}}_{O_0}: =g^{\mathbb{T}}$.  We are going to define $g^{\mathbb{T}}_{O_j}$  inductively  and show that  
	\begin{equation}\label{eq: recursion}
	Eg^{\mathbb{T}} =  Eg^{\mathbb{T}}_{O_j} + Eg^{\mathbb{T}}_{S_{T(j)}} +\cdots +Eg^{\mathbb{T}}_{S_1} +\mathit{RapDec}(R)\|g\|_{L^2} \text{ on each } O_j, 
	\end{equation}
	and 
	\begin{equation}\label{eq: recl2}
	\sum_{O_j\in \mathcal{O}}\|g^{\mathbb{T}}_{O_j}\|_{L^2}^2  \lesssim d^{j-T(j)} \mathit{Poly}(d)^{T(j)} \|g\|_{L^2}^2, 
	\end{equation}
	where $S_t$,  $1\leq t\leq T(j)$, are ancestors of $O_j$. 
 Assume that for some $j$ with $1\leq j\leq J$, 
	\begin{itemize}
		\item 	we have defined $g^{\mathbb{T}}_{O_{j-1}}$ and $g^{\mathbb{T}}_{S_t}$ for $t\leq T(j-1)$, 
		\item  inequalities~\eqref{eq: recursion} and  \eqref{eq: recl2}  hold for $j-1$, 
		\item equation~\eqref{eq: relpi} holds for $t\leq T(j-1)$. 
	\end{itemize}

	The base case is when  $j-1=0$, which holds automatically because $O_0$ has no ancestor.
	
	\vspace{5pt}
	
 \textbf{ Case A.} If $T(j)=T(j-1)$, then \eqref{eq: relpi} holds. We perform the wave packet decomposition of $g^{\mathbb{T}}_{O_{j-1}}$ on the ball $B_{R_{j-1}}$ containing $O_{j-1}$:
 \[
 g^{\mathbb{T}}_{O_{j-1}} = \sum_{(\tau_j, w_j): d(\tau_j)=R_{j-1}^{-1/2}} (g^{\mathbb{T}}_{O_{j-1}})_{\tau_j, w_j} +\mathit{RapDec}(R_{j-1}) \|g^{\mathbb{T}}_{O_{j-1}}\|_{L^2}. 
 \]
Since $R_{j-1}> R^{\delta}$, $\mathit{RapDec}(R_{j-1}) = \mathit{RapDec}(R)$. By Lemma~\ref{lem: fewerwp},  
$\|g^{\mathbb{T}}_{O_{j-1}}\|_{L^2}\lesssim  \|g\|_{L^2}$.  As a result, 
\[
 g^{\mathbb{T}}_{O_{j-1}} = \sum_{(\tau_j, w_j): d(\tau_j)=R_{j-1}^{-1/2}} (g^{\mathbb{T}}_{O_{j-1}})_{\tau_j, w_j} +\mathit{RapDec}(R) \|g\|_{L^2}. 
\]
 Then we
 define 
 \begin{equation}\label{eq: gcell}
 g^{\mathbb{T}}_{O_j} \coloneqq \sum_{T_{\tau_j, w_j}\cap O_j\neq \emptyset}(g^{\mathbb{T}}_{O_{j-1}})_{\tau_j, w_j}, 
\end{equation}
where $O_{j-1}$ is the parent of $O_j$. 
By Lemma~\ref{lem: wpsupp}, the function $Eg^{\mathbb{T}}_{O_j}$ restricted to $O_j$ is: 
\begin{align*}
Eg^{\mathbb{T}}_{O_j}  = Eg^{\mathbb{T}}_{O_{j-1}} + \mathit{RapDec}(R)\|g\|_{L^2}.
\end{align*}

Therefore \eqref{eq: recursion} holds for $j$ by the assumption that \eqref{eq: recursion} holds for $(j-1)$.
Since $T(j) = T(j-1)$, we know that  $j\neq j_t$ for all $1\leq t\leq n$. By (5) in Lemma~\ref{lem: recursion},  a tube  $T_{\tau_j, w_j}$ of length $R_{j-1}$ and radius $R_{j-1}^{1/2+\delta}$ intersects $\lesssim d$ children $O_j$ of each $O_{j-1}$. Hence, 
\begin{equation}\label{eq: l2jj}
\sum_{O_j\subset O_{j-1}}\|g^{\mathbb{T}}_{O_j}\|_{L^2}^2\lesssim d\|g^{\mathbb{T}}_{O_{j-1}}\|_{L^2}^2.
\end{equation}

Inequality~\eqref{eq: recl2} for $j$ follows from \eqref{eq: l2jj}  and the assumption that \eqref{eq: recl2} holds for $(j-1)$. 

\vspace{5pt}

\textbf{Case B.}	If $T(j)= T(j-1)+1$, then for $S_{T(j)}= O_j$ and the parent $O_{j-1}$ of $O_{j}$, we define 
	\begin{equation}\label{eq: surfta}
	g^{\mathbb{T}}_{S_{T(j)}} \coloneqq \sum_{T_{\tau_j, w_j}\in \mathbb{T}_{S_{T(j)}}} (g^{\mathbb{T}}_{O_{j-1}})_{\tau_j, w_j}, 
	\end{equation}
	and 
	\begin{equation}\label{eq: surftr}
	g^{\mathbb{T}}_{O_j} \coloneqq \sum_{T_{\tau_j, w_j}\in \mathbb{T}_{S_{T(j)}, \mathit{trans}}} (g^{\mathbb{T}}_{O_{j-1}})_{\tau_j, w_j}. 
	\end{equation}

	By Lemma~\ref{lem: wpsupp} and \eqref{eq: surfta} \eqref{eq: surftr}, 
	\begin{equation}\label{eq: decj}
	Eg^{\mathbb{T}}_{O_{j-1}} = Eg^{\mathbb{T}}_{O_j} + Eg^{\mathbb{T}}_{S_{T(j)}} + \mathit{RapDec}(R)\|f\|_{L^2} \text{ on each } O_j.
	\end{equation}
	Therefore \eqref{eq: recursion} for $j$ follows from ~\eqref{eq: decj}  and the assumption that \eqref{eq: recursion} holds for $(j-1)$. 
	By Lemma~\ref{lem: fewerwp}, 
	\[ \|g^{\mathbb{T}}_{O_j}\|_{L^2} \lesssim \| g^{\mathbb{T}}_{O_{j-1}}\|_{L^2}, \]
	and so \eqref{eq: recl2} for $j$ holds by  the assumption that \eqref{eq: recl2} holds for $(j-1)$. 
	
Now we are ready to verify \eqref{eq: deccell}, \eqref{eq: strcell} and \eqref{eq: strsurf} using the fact that  \eqref{eq: recursion} and \eqref{eq: recl2} hold for all $j$ with  $1\leq j\leq J$.  Equation~\eqref{eq: recursion} for $j=J$ implies  \eqref{eq: deccell}.  Inequality~\eqref{eq: recl2} for $j=J$ implies  \eqref{eq: strsurf}. Since there are $O(R_{j_t-1}^{3\delta})\leq R^{3\delta}$  fat $R_{j_t-1}$--surfaces  $S_{t}$ in each $O_{j_t-1}$,  inequality~\eqref{eq: recl2} for $j=j_t-1$ implies \eqref{eq: strcell}:
	\begin{align*}
	\sum_{S_t\in \mathcal{S}_t}\|g^{\mathbb{T}}_{S_t}\|_{L^2}^2 & = \sum_{S_t\subset O_{j_t-1}} \sum_{O_{j_t-1}\in \mathcal{O}}\|g^{\mathbb{T}}_{S_t}\|_{L^2}^2 \\
	&\leq R^{3\delta}\sum_{O_{j_t-1}\in \mathcal{O}}\|g^{\mathbb{T}}_{O_{j_t}-1}\|_{L^2}^2 \\
	&\leq R^{3\delta} d^{j_t-t} \mathit{Poly}(d)^{t-1} \|g\|_{L^2}^2.
	\end{align*}
	
 Finally,  it remains to verify  \eqref{eq: relpi}.  If $n(j)=1$, then consider the ancestors of $S_1$:  $O_1\supset \cdots \supset O_{j_1-1}\supset  S_1$. 
We decompose the function $g^{\mathbb{T}}_{\Pi_{S_1}}$ as
	\begin{align}\label{eq: dectail}
	g^{\mathbb{T}}_{\Pi_{S_1}}=& \sum_{T_{\tau_{j_1}, w_{j_1}}\in \mathbb{T}_{S_1}} (\sum_{T_{\theta,v}\in \mathbb{T}, T_{\theta,v}\cap O_1\neq \emptyset}  g_{\theta,v})_{\tau_{j_1}, w_{j_1}}  \\ &-  \sum_{T_{\tau_{j_1}, w_{j_1}}\in \mathbb{T}_{S_1}} (\sum_{T_{\theta,v}\in \mathbb{T}, T_{\theta,v}\cap S_1= \emptyset, T_{\theta,v} \cap O_1\neq \emptyset}  g_{\theta,v})_{\tau_{j_1}, w_{j_1}}. \nonumber
	\end{align}

		 The second term on the right-hand side of \eqref{eq: dectail} is bounded above by $\mathit{RapDec}(R)\|g\|_{L^2}$ by Claim~\ref{claim: tail}, whose proof will be presented later. 
	Recall that $O_0 = B_R, R_0=R$ and   a tube $T_{\tau_1, w_1}$ is of length $R_0$ and radius $R_0^{1/2+\delta}$, so   $T_{\theta,v}$ can also be viewed as a tube of the form $T_{\tau_1, w_1}$.  By \eqref{eq: gcell}, we rewrite the first term on the right-hand side of \eqref{eq: dectail} as follows 
	\begin{align*}
	&\sum_{T_{\tau_{j_1}, w_{j_1}}\in \mathbb{T}_{S_1}} (\sum_{T_{\theta,v}\in \mathbb{T}, T_{\theta,v}\cap O_1\neq \emptyset}  g_{\theta,v})_{\tau_{j_1}, w_{j_1}}\\
	=& \sum_{T_{\tau_{j_1}, w_{j_1}}\in \mathbb{T}_{S_1}} ( \sum_{T_{\tau_1,w_1} \cap O_1\neq \emptyset} g^{\mathbb{T}}_{O_0, \tau_1, w_1})_{\tau_{j_1}, w_{j_1}} + \mathit{RapDec}(R)\|g\|_{L^2}\\
	=& \sum_{T_{\tau_{j_1}, w_{j_1}}\in \mathbb{T}_{S_1}} (g^{\mathbb{T}}_{O_1})_{ \tau_{j_1}, w_{j_1}} +\mathit{RapDec}(R)\|g\|_{L^2}\\
	=& \sum_{T_{\tau_{j_1}, w_{j_1}}\in \mathbb{T}_{S_1}} (\sum_{T_{\tau_2, w_2}\cap O_2\neq \emptyset} g^{\mathbb{T}}_{O_1, \tau_2, w_2})_{ \tau_{j_1}, w_{j_1}} + \sum_{T_{\tau_{j_1}, w_{j_1}}\in \mathbb{T}_{S_1}} (\sum_{T_{\tau_2, w_2}\cap O_2= \emptyset} g^{\mathbb{T}}_{O_1, \tau_2, w_2})_{ \tau_{j_1}, w_{j_1}} +\mathit{RapDec}(R)\|g\|_{L^2}.
	\end{align*}
	For the second term on the right-hand side,  we apply Lemma~\ref{lem: slwp} 
			 \footnote{Strictly speaking, to apply Lemma~\ref{lem: slwp}, $T_{\tau_2,w_2}\cap O_2 = \emptyset$ under the summation  should be replaced by  $C T_{\tau_2,w_2}\cap  O_2 =  \emptyset$ for a constant $C$ due to ``$\lesssim$'' in (1) of Lemma~\ref{lem: slwp} . It is convenient to  omit the constant $C$ from our notation as it does not affect the main arguments.   }
	repeatedly to bound  it by $\mathit{RapDec}(R)\|g\|_{L^2}$ using the fact that $S_1\subset O_{2}$.  To see this, if $R_{j_1-1}\geq R_1^{1/2}$, then it follows from  a direction application of Lemma~\ref{lem: slwp}.  Otherwise, we take $2<j < j_1$ such that $R_{j_1-1}^{1/2} \geq R_{j-1}^{1/2}$, by Lemma~\ref{lem: slwp}, it suffices to show that 
	\[
	\sum_{T_{\tau_j, w_j} \cap O_j \neq \emptyset} (\sum_{T_{\tau_2, w_2}\cap O_2=\emptyset} g^{\mathbb{T}}_{O_1, \tau_2, w_2})_{\tau_j, w_j} \leq \mathit{RapDec}(R)\|g\|_{L^2}. 
	\]
	This again follows by applying Lemma~\ref{lem: slwp}  repeatedly as above. 
	
	Now we apply the arguments  above  recursively  to show that 
	\begin{align*}
		&\sum_{T_{\tau_{j_1}, w_{j_1}}\in \mathbb{T}_{S_1}} (\sum_{T_{\theta,v}\in \mathbb{T}, T_{\theta,v}\cap O_1\neq \emptyset}  g_{\theta,v})_{\tau_{j_1}, w_{j_1}}\\
	=& \sum_{T_{\tau_{j_1}, w_{j_1}}\in \mathbb{T}_{S_1}} (\sum_{T_{\tau_2, w_2}\cap O_2\neq \emptyset} g^{\mathbb{T}}_{O_1, \tau_2, w_2})_{\tau_{j_1}, w_{j_1}} +\mathit{RapDec}(R)\|g\|_{L^2} \\
	=& \cdots = \sum_{T_{\tau_{j_1}, w_{j_1}}\in \mathbb{T}_{S_1}} (g^{\mathbb{T}}_{O_{j_1-1}})_{ \tau_{j_1}, w_{j_1}} +\mathit{RapDec}(R)\|g\|_{L^2} \\
	= & g^{\mathbb{T}}_{S_1}  +\mathit{RapDec}(R)\|g\|_{L^2}. 
	\end{align*}

\vspace{5pt}

If $n(j)=t>1$, then consider the ancestors of $S_t$: $O_1\supset \cdots \supset O_{j_1-1} \supset O_{j_1}=S_1\supset \cdots \supset S_t$. 
 By the argument above  for $n(j)=1$ and Claim~\ref{claim: tail},  we have 
\begin{align*}
g^{\mathbb{T}}_{\Pi_{S_t}} &=\sum_{T_{\tau_{j_t}, w_{j_t}}\in \mathbb{T}_{S_t}} ( \sum_{T_{\theta,v}\in \mathbb{T}, T_{\theta,v}\cap S_t\neq \emptyset} g_{\theta,v})_{\tau_{j_t}, w_{j_t}} \\
& = \sum_{T_{\tau_{j_t}, w_{j_t}}\in \mathbb{T}_{S_t}} ( g^{\mathbb{T}}_{O_{j_1-1}})_{\tau_{j_t}, w_{j_t}}  +\mathit{RapDec}(R)\|g\|_{L^2}
	\\
&= \sum_{T_{\tau_{j_t}, w_{j_t}}\in \mathbb{T}_{S_t}} ( g^{\mathbb{T}}_{O_{j_1}} + g^{\mathbb{T}}_{S_{1}})_{\tau_{j_t}, w_{j_t}}    + \sum_{T_{\tau_{j_t}, w_{j_t}} \in \mathbb{T}_{S_t}} ( \sum_{T_{\tau_{j_1}, w_{j_1}}\cap S_1=\emptyset } g^{\mathbb{T}}_{O_{j_1-1}, \tau_{j_1}, w_{j_1}})_{\tau_{j_t}, w_{j_t}}+\mathit{RapDec}(R)\|g\|_{L^2} \\
&= \sum_{T_{\tau_{j_t}, w_{j_t}}\in \mathbb{T}_{S_t}} ( g^{\mathbb{T}}_{O_{j_1}} + g^{\mathbb{T}}_{S_{1}})_{\tau_{j_t}, w_{j_t}}  + \mathit{RapDec}(R)\|g\|_{L^2}.
\end{align*}
The second to last equality is due to 	\eqref{eq: decj}, while  the last inequality is due to Claim~\ref{claim: tail}. 
By the definition of $g^{\mathbb{T}}_{S_1, S_t}$, we conclude that 
\begin{align*}
g^{\mathbb{T}}_{\Pi_{S_t}} & =  \sum_{T_{\tau_{j_t}, w_{j_t}}\in \mathbb{T}_{S_t}} ( g^{\mathbb{T}}_{O_{j_1}} )_{\tau_{j_t}, w_{j_t}}   + g^{\mathbb{T}}_{S_1, S_t}+\mathit{RapDec}(R)\|g\|_{L^2}\\
&=\cdots = \sum_{T_{\tau_{j_t}, w_{j_t}}\in \mathbb{T}_{S_t}} ( g^{\mathbb{T}}_{O_{j_{t-1}}} )_{\tau_{j_t}, w_{j_t}}   +g^{\mathbb{T}}_{S_{t-1}, S_t}+\cdots +  g^{\mathbb{T}}_{S_1, S_t}+\mathit{RapDec}(R)\|g\|_{L^2}\\
&= g^{\mathbb{T}}_{S_t} +g^{\mathbb{T}}_{S_{t-1}, S_t}+\cdots +  g^{\mathbb{T}}_{S_1, S_t}+\mathit{RapDec}(R)\|g\|_{L^2}.
\end{align*}
To complete the proof of Lemma~\ref{lem: properties}, it remains to verify Claim~\ref{claim: tail}.

\end{proof}

\begin{claim}\label{claim: tail}
	For $1\leq t\leq n$  and any $O_{l}\supseteq O_{j}\supset S_{j_t}$, we have 
	\begin{equation}\label{eq: indtail}
	\sum_{T_{\tau_{j_t}, w_{j_t}}\in \mathbb{T}_{S_t}} (\sum_{T_{\tau_j, w_j} \cap S_{t}=\emptyset , T_{\tau_j, w_j}\cap O_j\neq \emptyset } g^{\mathbb{T}}_{O_l,\tau_j, w_j} )_{\tau_{j_t}, w_{j_t}}  = \mathit{RapDec}(R)\|g\|_{L^2}.
	\end{equation}
\end{claim}
	\begin{proof}
		We are going to induct on the ratio $\frac{\log R_{j-1}}{\log R_{j_t-1}}$.  Since $j_t-1<J$ and $R_{j_t-1}\geq R^{\delta}$,  the ratio  $\frac{\log R_{j-1}}{\log R_{j_t-1}}$ is between $1$ and $\delta^{-1}$.

		 When $\frac{\log R_{j-1}}{\log R_{j_t-1}} \leq 2$, 
	 inequality~\eqref{eq: indtail} follows from a direct application of   Lemma~\ref{lem: slwp}.

		 When  $\frac{\log R_{j-1}}{\log R_{j_t-1}}\geq 2$, 
		 we pick some  $j<j'<j_t$ such that  $0<\frac{\log R_{j-1}}{\log R_{j'-1}} \leq 2$. 
		 	 By the induction hypothesis,  \eqref{eq: indtail} holds when $j$ is  replaced by $j'$. 
		 For $S_t\subset O_{j'} \subset O_j$, we have
		 	\begin{align}\label{eq: 3.25}
		 &\sum_{T_{\tau_{j_t}, w_{j_t}}\in \mathbb{T}_{S_t}} ( \sum_{T_{\tau_j,w_j} \cap S_t= \emptyset, T_{\tau_j,w_j} \cap O_j\neq \emptyset} g^{\mathbb{T}}_{O_l, \tau_j,w_j})_{\tau_{j_t}, w_{j_t}}\\
		 = & \sum_{T_{\tau_{j_t}, w_{j_t}}\in \mathbb{T}_{S_t}} ( \sum_{T_{\tau_j,w_j} \cap S_t= \emptyset, T_{\tau_j,w_j} \cap O_{j'}\neq \emptyset} g^{\mathbb{T}}_{O_l, \tau_j,w_j})_{\tau_{j_t}, w_{j_t}} +  
		  \mathit{RapDec}(R)\|g\|_{L^2} \nonumber \\
		 &+ \sum_{T_{\tau_{j_t}, w_{j_t}}\in \mathbb{T}_{S_t}} ( \sum_{T_{\tau_j,w_j} \cap O_{j'}= \emptyset, T_{\tau_j,w_j} \cap O_{j}\neq \emptyset} g^{\mathbb{T}}_{O_l, \tau_j,w_j})_{\tau_{j_t}, w_{j_t}}.  \nonumber
		 \end{align}
		 Now we apply Lemma~\ref{lem: slwp} repeatedly to bound the last term on the right-hand side of \eqref{eq: 3.25}  by $\mathit{RapDec}(R)\|g\|_{L^2}$ using the fact that $S_t\subset O_{j'}$.   
		 
		 To deal with the first term on the right-hand side of \eqref{eq: 3.25}, we apply  Lemma~\ref{lem: slwp} to show that 
		 \begin{align*}
		 &\sum_{T_{\tau_{j_t}, w_{j_t}}\in \mathbb{T}_{S_t}} ( \sum_{T_{\tau_j,w_j} \cap S_t= \emptyset, T_{\tau_j,w_j} \cap O_{j'}\neq \emptyset} g^{\mathbb{T}}_{O_l, \tau_j,w_j})_{\tau_{j_t}, w_{j_t}} \\
		 =&  \sum_{T_{\tau_{j_t}, w_{j_t}}\in \mathbb{T}_{S_t}} ( \sum_{T_{\tau_{j'},w_{j'}} \cap S_t= \emptyset, T_{\tau_{j'},w_{j'}} \cap O_{j'}\neq \emptyset} g^{\mathbb{T}}_{O_l, \tau_{j'},w_{j'}})_{\tau_{j_t}, w_{j_t}} + \mathit{RapDec}(R)\|g\|_{L^2}\\
		 =& \mathit{RapDec}(R)\|g\|_{L^2}.
		 \end{align*}
		 The last inequality follows from the  induction hypothesis since   $\frac{\log R_{j'-1}}{\log R_{j_t-1}}< \frac{\log R_{j-1}}{\log R_{j_t-1}}$.
	\end{proof}
The proof of Lemma~\ref{lem: properties} is now completed. 

\vspace{5pt}

Recall that in \eqref{eq: snsim} of  Subsection~\ref{subsection: te},  we have decomposed the function into a local part and a global part. In addition,  we have estimated  the local part $\sum_{B_k} \|Ef^{\sim}_k\|_{BL^p(B_k)}^p$  by induction on scale in Lemma~\ref{lem: telocal}. The following lemma is  a further decomposition of the global part using   Lemma~\ref{lem: properties}. 

%
\begin{lemma}\label{lem: findscales}
For the function $Ef^{\nsim}_k$ in the two ends reduction in Subsection~\ref{subsection: te}, either of the following happens:
\begin{enumerate}[label=(\roman*)]
	\item \label{case1}
 there exist 
a $\geq 2R^{-\delta}$--fraction of the leaves $O\in \mathcal{O}$ \footnote{ In fact, one should read inequality~\eqref{eq: celltrans} as $(\sum_{B_k} \|Ef^{\nsim}_k\|_{BL^p(O)}^p)^{1/p} \lesssim R^{\delta }(\sum_{B_k} \|Ef^{\nsim}_{k, O}\|_{BL^p(O)}^p)^{1/p}$. But each $O$ lies in  $\lesssim 1$ balls $B_k$ defined in Subsection~\ref{subsection: te}, so we drop the summation and write $(\sum_{B_k} \|Ef^{\nsim}_k\|_{BL^p(O)}^p)^{1/p}$ as $\|Ef^{\nsim}\|_{BL^p(O)}$ to ease the notation. }, 
\begin{equation}\label{eq: celltrans}
\|Ef^{\nsim}\|_{BL^p(O)} \leq  R^{\delta}\|Ef^{\nsim}_O\|_{BL^p(O)};
\end{equation}
\item \label{case2}
we  can choose the smallest integer $t$ with $1\leq t\leq n$ such that there exists a subset $\mathcal{O}(t)$ of the leaves in $\mathcal{O}$ with
\begin{equation}\label{eq: Otsize}
|\mathcal{O}(t)|\geq 2 R_{j_{t-1}}^{-\delta} | \{O: O\in \mathcal{O}\}|
\end{equation}
and for each $O\in \mathcal{O}(t)$, we have 
\begin{equation}\label{eq: 55}
\|Ef^{\nsim}\|_{BL^p(O)} \leq R_{j_{t-1}}^{\delta} \|Ef^{\nsim}_{S_t}\|_{BL^p(O)}, 
\end{equation}
and for all $1\leq l<t$, 
\begin{equation}\label{eq: tang}
\|Ef^{\nsim}\|_{BL^p(O)} \geq R_{j_{l-1}}^{\delta} \|Ef^{\nsim}_{S_l}\|_{BL^p(O)}.
\end{equation}

\end{enumerate}
\end{lemma}
\begin{proof}
Assume that \eqref{eq: celltrans}  holds for  only a $\leq 2R^{-\delta}$   fraction of the leaves $O\in \mathcal{O}$, then we show that the  integer $t$  in \ref{case2} exists by pigeonholing. Indeed, we first look at $S_1$. If there exists a subset $\mathcal{O}(1)$ of the leaves $O$ with $$|\mathcal{O}(1)|\geq 3R^{-\delta} |\{O: O\in \mathcal{O}\}|$$ and 
\begin{equation}\label{eq: con}
\|Ef^{\nsim}\|_{BL^p(O)} \leq R^{\delta} \|Ef^{\nsim}_{S_1}\|_{BL^p(O)},
\end{equation}
then we choose $t=1$. Otherwise, there exist  a $\geq (1-3R^{-\delta})$--fraction of the leaves $O$ that do not satisfy inequality~\eqref{eq: con}.  In this case, we focus on the aforementioned  $(1-3R^{-\delta})$--fraction of leaves $O\in \mathcal{O}$.  We choose $t=2$ if  there exists a subset $\mathcal{O}(2)$  such that 
\eqref{eq: 55} holds for $t=2$ for each $O\in \mathcal{O}(2)$ and 
 $$|\mathcal{O}(2)|\geq 3R_{j_1}^{-\delta} (1-3R^{-\delta}) |\{O: O\in \mathcal{O}\}|.$$
 Otherwise, we focus  on the leaves not satisfying \eqref{eq: con} or \eqref{eq: 55} for $t=2$, 
 which form a $$\geq (1-3R_{j_1}^{-\delta} ) (1-3R^{-\delta}) \geq 1-3R^{-\delta}  - 3R_{j_1}^{-\delta}$$ fraction of the original leaves. As we iterate this process, at each step,  either we find an integer $t$ satisfying  \eqref{eq: Otsize}, \eqref{eq: 55} and \eqref{eq: tang}, or the estimates
\[ \|Ef^{\nsim}\|_{BL^p(O)}\geq R_{j_{l-1}}^{\delta} \|Ef^{\nsim}_{S_l}\|_{BL^p(O)} , \quad 1 \leq l\leq  n, \]
hold  for   a $\geq 1-3R^{-\delta} - 3R_{j_{1}}^{-\delta} - \cdots - 3R_{j_{n-1}}^{-\delta} \geq 1- 3n R^{-\delta} \geq 2/3$ fraction of the original leaves.  However, by the triangle inequality and  \eqref{eq: deccell}, this contradicts  the assumption that \eqref{eq: celltrans}  holds for only a $\leq 2R^{-\delta}$--fraction of the leaves. 
\end{proof}
\begin{lemma}\label{lem: find scale}
	Let $\mathcal{O}$ be the tree defined in Corollary~\ref{cor: tree}. If there is a leaf  $O\in \mathcal{O}$ satisfying both \eqref{eq: 55} and \eqref{eq: tang}, 
then
\begin{equation}\label{eq: twin}
			 \|Ef^{\nsim}_{S_t}\|_{BL^p(O)}\sim \|Ef^{\nsim}_{\Pi_{S_t}}\|_{BL^p(O)}.
			 \end{equation}
\end{lemma}

The proof of Lemma~\ref{lem: find scale} is independent of the rest of the argument, and  is deferred to Section~\ref{section: 38}.

\section{Easy cases}\label{section: easy}

Let us continue the discussion of the global part following Lemma~\ref{lem: findscales} and assuming that the global part dominates: 
\begin{equation}\label{eq: Stdominate}
\text{ for a }\geq (1-R^{-\delta})\text{--fraction of the leaves } O \in \mathcal{O},   \, \,  \|Ef\|_{BL^p(O)} \lesssim \|Ef^{\nsim}\|_{BL^p(O)}.
\end{equation}

 In this section, we deal with some easy cases in Lemma~\ref{lem: findscales}. We start with \ref{case1} of  Lemma~\ref{lem: findscales}. 

\begin{lemma}\label{lem: small} Let $\mathcal{O}$ be the tree defined in Corollary~\ref{cor: tree}. If \eqref{eq: Stdominate} holds and case	\ref{case1} of  Lemma~\ref{lem: findscales} happens, 	then 
	\begin{equation}\label{eq: smallgoal}
	\|Ef\|_{BL^p(O)} \lesssim R^{O(\delta)} \|f\|_{L^2} \text{ for any } p\geq 3.
	\end{equation}
	\end{lemma}
\begin{proof}
	By \eqref{eq: Stdominate} and \ref{case1} of  Lemma~\ref{lem: findscales}, there exists a subset $\mathcal{O}(0)$ of the leaves in $\mathcal{O}$ such that
	\begin{itemize}
		\item $|\mathcal{O}(0)|\geq R^{-\delta} |\{ O: O\in \mathcal{O}\}|$;
		\item $\|Ef\|_{BL^p(O)} \lesssim \|Ef^{\nsim}\|_{BL^p(O)} \lesssim R^{\delta} \|Ef^{\nsim}_{O}\|_{BL^p(O)}.$
	\end{itemize}
By Corollary~\ref{cor: tree}, for each leaf $O \in \mathcal{O}(0)$, 
\begin{equation}\label{eq: small1}
\|Ef\|_{BL^p(B_R)}^p  \lesssim R^{\delta}  |\{O: O\in \mathcal{O}\}|  \cdot \|Ef \|_{BL^p(O)}^p \lesssim R^{2\delta}  |\{O: O\in \mathcal{O}\}|  \cdot  \|Ef^{\nsim}_O\|_{BL^p(O)}^p.
\end{equation}

By Lemma~\ref{lem: tao}, since $O\subset B_{R_J}$ and $R_J\leq R^{\delta}$, $\epsilon^7 \leq \delta \leq \epsilon^3$, 
\begin{equation}\label{eq: small2}
\|Ef^{\nsim}_O\|_{BL^p(O)}^p \lesssim R_J^{\frac{5}{2}-\frac{3p}{4}}  R^{ \epsilon^{10}} \|f^{\nsim}_O\|_{L^2}^p \lesssim R^{O(\delta)} \|f^{\nsim}_O\|_{L^2}^p.
\end{equation}

Using inequality~\eqref{eq: strsurf} in Lemma~\ref{lem: properties}, for a typical leaf $O\in \mathcal{O}(0)$, 
\begin{equation}\label{eq: small3}
\|f^{\nsim}_O\|_{L^2}^2 \lesssim \frac{d^{J-n} \mathit{Poly}(d)^n R^{\delta}}{|\{O: O\in \mathcal{O}\}|}\|f\|_{L^2}^2 \lesssim \frac{d^J}{ |\{O: O\in \mathcal{O}\}|} R^{2\delta} \|f\|_{L^2}^2.
\end{equation}
The last inequality follows from the fact that  $d\lesssim \log R$ and $n\leq \delta^{-2}$. 

Combining  \eqref{eq: small1}, \eqref{eq: small2} and \eqref{eq: small3}, we have 
\begin{equation}\label{eq: small4}
\|Ef\|_{BL^p(B_R)}^p \lesssim R^{O(\delta)}  |\{ O: O\in \mathcal{O}\}|^{1-p/2} d^{Jp/2}\|f\|_{L^2}^p.
\end{equation}
Using \eqref{eq: indd} in Lemma~\ref{lem: recursion} and inequality~\eqref{eq: small1}, we derive that 
\begin{equation}\label{eq: small5}
|\{O: O\in \mathcal{O}\} |\gtrsim d^{3(J-n) } R^{-2\delta}\gtrsim d^{3J} R^{-3\delta}
\end{equation}
because $d\lesssim \log R$ and $n\leq \delta^{-2}$. 

Now  \eqref{eq: smallgoal} follows from  combining \eqref{eq: small4} and \eqref{eq: small5}.

\end{proof}

Then we deal with  \ref{case2} of Lemma~\ref{lem: findscales} when the radius $R_{j_t-1}$ is large. 
\begin{lemma}\label{lem: large}
	Let $\mathcal{O}$ be the tree defined in  Corollary~\ref{cor: tree} and assume \eqref{eq: Stdominate}.  Then for any $p\geq 3$,  if  case \ref{case2} of Lemma~\ref{lem: findscales} holds for  the integer  $t$ and   $R_{j_t-1}\geq R^{1/(p-2)}$, 
	then the estimate \eqref{eq: induction2} in Theorem~\ref{thm: induction2} holds for $f$. 
\end{lemma}
\begin{proof}
	By Corollary~\ref{cor: tree}, \eqref{eq: Stdominate} and \ref{case2} of Lemma~\ref{lem: findscales}, 
	\begin{align*}
	\|Ef\|_{BL^p(B_R)}^p  \lesssim R^{\delta} \sum_{O\in \mathcal{O}} \|Ef\|_{BL^p(O)}^p \lesssim R^{O(\delta)}  \sum_{O\in \mathcal{O}} \|Ef^{\nsim}_{S_t}\|_{BL^p(O)}^p
	\lesssim R^{O(\delta)}  \sum_{S_t\in \mathcal{S}_t} \|Ef^{\nsim}_{S_t}\|_{BL^p(S_t)}^p.
	\end{align*} 
	We apply  Lemma~\ref{lem: tao} to bound  each  $\|Ef^{\nsim}_{S_t}\|_{BL^p(S_t)}$ that appears in the last summation. Since $S_t$ is a fat $R_{j_t-1}$--surface and $R^{\epsilon^{10}} \leq R^{\delta}$, we have\footnote{Strictly speaking, when $R_{j_t}\geq R^{1-\epsilon_0}$, there should be a summation over $B_k$:
		\[ \|Ef\|_{BL^p(B_R)}^p \lesssim R^{2\delta} \sum_{B_k} \sum_{S_t\in \mathcal{S}_t} \|Ef^{\nsim}_{k, S_t}\|_{BL^p(S_t\cap B_k)}^p. \]
		 This causes a loss of $\lesssim R^{3\epsilon_0}$. But the $R^{3\epsilon_0}$--loss is acceptable since $3\epsilon_0< \epsilon/10$. }
	\begin{equation}\label{eq: 65}
	\|Ef\|_{BL^p(B_R)}^p\lesssim R^{O(\delta)}  \sum_{S_t\in \mathcal{S}_t} \|Ef^{\nsim}_{S_t}\|_{BL^p(S_t)}^p\lesssim R^{O(\delta)} R_{j_t}^{\frac{5}{2}-\frac{3p}{4}}\sum_{S_t\in \mathcal{S}_t} \|f^{\nsim}_{S_t}\|_{L^2}^p.
	\end{equation}
Now we apply Lemma~\ref{lem: geometric}, Lemma~\ref{lem: comparel2} and  Lemma~\ref{lem: fewerwp} to obtain
	\begin{equation}\label{eq: 66}
	\|f^{\nsim}_{S_t}\|_{L^2}^2 \lesssim R_{j_t}^{-1/2+O(\delta)} \underset{d(\tau_t)=r_t^{-1/2}}{\max} \|f^{\nsim}_{S_t, \tau_t}\|_{L^2_{avg}(\tau_t)}^2 \lesssim      R_{j_t}^{-1/2+O(\delta)}     \underset{d(\theta)=R^{-1/2}}{\max} \|f_{\theta}\|_{L^2_{avg}(\theta)}^2. 
	\end{equation}
	
	By inequality~\eqref{eq: strcell} in Lemma~\ref{lem: properties} and the fact that $d\lesssim \log R, t \leq \delta^{-2}$, 
	\begin{equation}\label{eq: 67}
       \sum_{S_t\in \mathcal{S}_t}\|f^{\nsim}_{S_t}\|_{L^2}^2 \lesssim d^{j_t-t} \mathit{Poly}(d)^{t-1} R^{3\delta} \|f\|_{L^2}^2 \lesssim d^{j_t} R^{4\delta} \|f\|_{L^2}^2. 
	\end{equation}
	
 Using  (2) of Lemma~\ref{lem: recursion}, we have $R_{j_t} \leq R/d^{j_t}$.  By our assumption,   $R_{j_t} \geq R^{1/(p-2)}$, so  $d^{j_t} \leq R_{j_t}^{p-3}$.  Now we combine \eqref{eq: 65}, \eqref{eq: 66} and \eqref{eq: 67} to finish the proof: 
	\begin{align*}
	\|Ef\|_{BL^p(B_R)}^p &\lesssim R^{O(\delta)} R_{j_t}^{\frac{5}{2}-\frac{3p}{4}} d^{j_t} R_{j_t}^{-\frac{p-2}{4}} \|f\|_{L^2}^2 \underset{d(\theta)=R^{-1/2}}{\max} \|f_{\theta}\|_{L^2_{avg}(\theta)}^{p-2}\\
	&\lesssim R^{O(\delta)} \|f\|_{L^2}^2 \underset{d(\theta)=R^{-1/2}}{\max} \|f_{\theta}\|_{L^2_{avg}(\theta)}^{p-2}.
	\end{align*}
	Since  $\delta \leq \epsilon^3$, the constant term  $ R^{O(\delta)}$ is much smaller than $ C_{\epsilon} R^{\epsilon}$ when $R$ is sufficiently large. 
\end{proof}

\section{Brooms and bushes}\label{section: broom}

In this section, we introduce the key geometric objects: brooms and bushes. 

\begin{definition}\label{def: ess}
	Let $S$ be a fat $r$--surface, we define $\mathbb{T}_{S, ess}$ to be the set of tubes $T_{\tau,w}\in \mathbb{T}_{S}$ (see Definition~\ref{def: tangsurf}) such that there exists another $T_{\tau',w'}\in \mathbb{T}_S$ with
	\begin{itemize}
		\item $T_{\tau,w}\cap T_{\tau',w'}\cap S\neq \emptyset$,
		\item $\mathit{Angle}(G(\omega_{\tau}), G(\omega_{\tau'})) \geq K^{-1}$.
	\end{itemize}
Recall that $T_{\tau,w}$ has length $r$ and radius $r^{1/2+\delta}$,  and $\omega_{\tau}$ is the center of $\tau$. The function $G$ is defined by \eqref{eq: Gtheta}  and   $1\leq K^{\epsilon^2} \leq R^{\epsilon^{12}/100}$.
\end{definition}

Let $f_S$ be a function concentrated on wave packets from $\mathbb{T}_S$. We define $f_{S, ess}  \coloneqq \sum_{T_{\tau, w}\in \mathbb{T}_{S, ess}}  f_{S, \tau, w}$. By the definition of $BL^p$--norm, $Ef_{S, ess}$ essentially represents $Ef_{S}$ in the sense that
\begin{equation}\label{eq: ess}
\|Ef_S\|_{BL^p(S)} \leq \|Ef_{S, ess}\|_{BL^p(S)} +\mathit{RapDec}(r)\|f_S\|_{L^2}.
\end{equation}

 \begin{definition}\label{def: fatplane}
 	Let $Z$ be a plane, a fat $r$--plane $\Sigma$ is defined as $$ \Sigma \coloneqq N_{r^{1/2+\delta}} Z
  \cap B_{r}.$$ Given a fat $r$--plane $\Sigma$, we use $Z_{\Sigma}$ to denote the plane $Z$. 
  Let $\mathbb{T}_{\Sigma}$ denote the set of tubes $T_{\tau,w} \subset \Sigma$ of length $r$ and radius $r^{1/2+\delta}$.  
 \end{definition}
Here a fat $r$--plane is slightly different from  a  fat $r$--surface with $Z_S=Z$ because we intersect with a ball  $B_r$ instead of $B_{r^{1-\delta}}$.  
 \begin{prop}\label{lem: plane}
 	Let $S$ be a fat $r$--surface.  For a fixed cap $\tau$ of radius $r^{-1/2}$, let $\mathbb{T}_{S,\tau} \subset \mathbb{T}_{S, \mathit{ess}}$ be the set of tubes $T_{\tau, w}\in \mathbb{T}_S$. Then there exists a set $\Omega(S, \tau)$ of  fat $r$--planes $\Sigma$, such that
 	\begin{enumerate}
 		\item $|\Omega(S, \tau)| \lesssim_d r^{O(\delta)}$,
 		\item 
 	 $\mathbb{T}_{S, \tau} \subset \bigcup_{\Sigma \in \Omega(S, \tau)} \mathbb{T}_{\Sigma},$
 	 where $\mathbb{T}_{\Sigma}$ is defined in Definition~\ref{def: fatplane}. 
 	 \end{enumerate}
 \end{prop}

It is worth noting that Proposition~\ref{lem: plane} fails without restricting to $\mathbb{T}_{S, \mathit{ess}}$. For example, if $Z_S$ is a cylinder of radius $< r^{1-2\delta}$ and $G(\omega_{\tau})$ is parallel to the axis, then Proposition~\ref{lem: plane} does not hold for the collection of tubes  $\mathbb{T}_S$. But in this case, the set $\mathbb{T}_{S, \mathit{ess}}$ is empty, so Proposition~\ref{lem: plane} holds trivially. 
We postpone the proof of Proposition~\ref{lem: plane} to Section~\ref{section: plane}.

From now on, we write $f_S=f_{S, ess}$ and $\mathbb{T}_S = \mathbb{T}_{S, ess}$. 

\subsection{Definition of a broom}
In this subsection, we assume that $R^{1/2+\delta}\leq r\leq R$ and fix a cap $\tau$ of radius $r^{-1/2}$. Recall that $\omega_{\tau} = (\omega_{\tau,1}, \omega_{\tau,2})$ is  the center of $\tau$.   Let $f_{\tau}$ be a function supported on $\tau$ and $\tau_0$ be another cap of the same radius and centered at the origin.  Then we can find a  Galilean transformation  
\begin{equation}\label{eq: A}
\mathcal{A}: (\xi_1, \xi_2, \xi_3)\mapsto (\xi_1 + \omega_{\tau, 1}, \,\, \xi _2 +\omega_{\tau, 2}, \,\, \xi_3+ 2\xi_1\omega_{\tau, 1} + 2\xi_2\omega_{\tau, 2} + |\omega_{\tau}|^2 ), 
\end{equation}
such that  $Ef_{\tau} ( (\mathcal{A}^{-1})^{T}x) = Eg_{\tau_0}(x)$ for some function  $g_{\tau_0}$  supported on the cap $\tau_0$.  Once we have defined the geometric objects associated to $\tau_0$, we can define the ones associated to a general cap $\tau$ of radius $r^{-1/2}$ using the  Galilean transformation \eqref{eq: A}. 
 For  the rest of this subsection, we assume that $\tau=\tau_0$ is centered at the origin. Note that $G(\omega_{\tau})=(0, 0 ,1)$. 
 

Let $\Sigma$ be a fat $r$--plane such that $ \mathit{Angle}( (0, 0, 1), Z_{\Sigma}) \leq r^{-1/2}$.  Then there is a plane of the form
\begin{equation}\label{eq: Z}
 Z=\{ (x_1, x_2, x_3): k_1 x_1 + k_2 x_2 = a\}
\end{equation}
 and a ball $B_{2r}$ such that  $\Sigma \subset N_{2r^{1/2+\delta}} Z \cap B_{2r}$, which is morally a fat $r$--plane defined by $Z$.  For the rest of this subsection, we  think of   $\Sigma$ as  defined by $Z_{\Sigma}=Z$. 

 Consider a large tube $T_{\theta,v}$ with $\theta\subset \tau$. If $T_{\theta,v}\cap \Sigma \neq \emptyset$, then $2T_{\theta,v}\cap \Sigma$ is equal to a plank $T$ of dimensions approximately $r \times R^{1/2+\delta}\times r^{1/2+\delta}$.  If two large tubes $T_{\theta,v}\cap T_{\theta', v'}\cap \Sigma \neq \emptyset$ and $\theta, \theta'\subset \tau$, then  $2T_{\theta,v} \cap \Sigma \approx 2T_{\theta', v'} \cap \Sigma \approx T$ in the sense that $T \subset 10 (2T_{\theta,v} \cap \Sigma), 10 (2T_{\theta', v'} \cap \Sigma) \subset 100 T$. 
 
 We decompose $\Sigma$ into a disjoint union of parallel planks $T$ whose long direction is  along $G(\omega_{\tau})$:
 \begin{equation}\label{eq: SigmaT}
 \Sigma =\bigsqcup T.
 \end{equation}
     For  $\Sigma = N_{r^{1/2+\delta}} Z \cap B_{r}$ with $Z$ defined in \eqref{eq: Z}, we decompose $\tau$ into a union of finitely overlapping strips $s$, where $s: = N_{R^{-1/2}}l\cap \tau$  for some line  $l$ of the form $$l=\{ (\xi_1, \xi_2): k_2 \xi_1 -k_1 \xi_2 =b\}.$$ Note that $l$ is orthogonal to $Z$.
 
\noindent
\begin{center}  
	\begin{minipage}{0.4\textwidth}
		\begin{tikzpicture}
		\draw(0,0)--(1,2)--(3.5, 2)--(2.5,0)--(0,0);
		\draw(0.2,-0.5)--(2.7, -0.5)--(3.7, 1.5);
		\draw[dashed] (0.2,-0.5)--(1.2, 1.5)--(3.7,1.5);
		\draw[dashed] (1.2, 1.5)--(1,2); 
		\draw (3.7, 1.5)--(3.5, 2) (0,0)--(0.2, -0.5) (2.7, -0.5)--(2.5, 0);
		\draw (3.5,0) node[below left]{$\Sigma$};
		\draw[|<-](0.2, -1)--(1, -1) node[right]{$r$};
		\draw[->|](1.5, -1)--(2.7, -1);
		\draw[|<->|](-0.2, -0.5)--(-0.2, 0) node[left]{$r^{1/2+\delta}$};
		\draw[blue, very thick] (1.5, 2)--(2.5,2)--(1.5,0)--(0.5,0)--(1.5, 2);
		\draw[blue, very thick] (1.5, 0)--(1.7, -0.5)--(0.7, -0.5)--(0.5,0);
		\draw[blue, very thick, dashed] (0.7, -0.5)--(1.7, 1.5)--(1.5, 2) (1.7, -0.5)--(2.7, 1.5)--(2.5, 2) (1.7, 1.5)--(2.7, 1.5);
		\draw[|<->|, blue] (1.5, 2.2)--(2.5, 2.2) node[blue, above left]{$R^{1/2+\delta}$};
		\draw (3, 2.2) node[blue]{$T$};
		\draw[blue, thick, ->] (3.5,0.4)--(4.1,1.6) node[blue, right]{$G(\omega_{\tau})$};
		\draw[|<-] (-0.2,0)--(0.2, 0.8) node[above]{$r$};
		\draw[->|] (0.4, 1.2)--(0.8, 2);
		\end{tikzpicture}
	\end{minipage}\qquad
	\begin{minipage}{0.4\textwidth}
		\begin{tikzpicture}[scale=0.5]
		\draw circle[radius=3];
		\draw(-3, 3) node{$\tau$};
		\draw (0.5,0) circle[radius=0.5] node{$\theta$};
		\draw (0.5,-1) circle[radius=0.5];
		\draw(0.5, 2) node{$s$};
		\foreach \x in {-2,...,0}
		\draw[blue, very thick] (\x, -0.42*\x-3)--(\x, 3+0.42*\x);
		\foreach \x in {1,2}
		\draw[blue, very thick] (\x, 0.42*\x-3)--(\x, 3-0.42*\x);
		\draw[|<->|](3.7, -3)--(3.7,3) node[below right]{$r^{-1/2}$};
		\draw[|<->|](0,3.2)--(1,3.2) node[above]{$R^{-1/2}$};
		\end{tikzpicture}
	\end{minipage}
\end{center}

\begin{definition}\label{def: broom}
	For each plank $T$ and strip $s$ defined as above in this subsection, we define the broom $\mathcal{B}$ associated to the pair $(s, T)$ as 
	\begin{equation}
	\mathcal{B}  \coloneqq  \{ T_{\theta,v}: T\subset 2T_{\theta,v} \text{ and } \theta\subset s\}.
	\end{equation}
	We call $T$ the root of the broom $\mathcal{B}$ and we denote it by $T_{\mathcal{B}}$. We say that a broom $\mathcal{B}$ is rooted at a fat $r$--plane $\Sigma$ if $T_{\mathcal{B}}\subset \Sigma$. 
\end{definition}
\begin{lemma}\label{lem: brsize}
There are $\lesssim (\frac{R}{r})^{1/2}$ tubes in a broom. 
\end{lemma}
\begin{proof}
	A broom $\mathcal{B}$ is associated to a unique pair $( s, T_{\mathcal{B}} )$. 
	Each $R^{-1/2}$--cap $\theta \subset s$ corresponds to $\lesssim 1$ tubes $T_{\theta,v}$ intersecting $T_{\mathcal{B}}$ because the long direction of  $T_{\mathcal{B}}$ points on  $G(\omega_{\tau})$ for the $r^{-1/2}$--cap $\tau$ containing $s$. In addition,  there are $\lesssim (\frac{R}{r})^{1/2}$ caps $\theta$ in $s$. 
\end{proof}
\begin{lemma}\label{remark: brplane}
	Let $\mathcal{B}$ be a broom defined in Definition~\ref{def: broom} and let $C_{T_{\mathcal{B}}}$ denote the center of $T_{\mathcal{B}}$. 	Then the tubes $T_{\theta,v}$ in $\mathcal{B}$ lie in the $R^{1/2+\delta}$--neighborhood of  the plane
	$$Z_{\mathcal{B}}\coloneqq \{ k_2(x_1- C_{T_{\mathcal{B}},1}) - k_1(x_2- C_{T_{\mathcal{B}},2}) +2b(x_3-C_{T_{\mathcal{B}},3})=0\}.  $$
	 Moreover, $Z_{\mathcal{B}}$ is orthogonal to the plane $Z$, which is parallel to  $Z_{\Sigma}$ up to an angle difference of $O(r^{-1/2})$ (see Figure~\ref{sbroom} on page 3). 
\end{lemma}
\begin{proof}
	By a translation in $\mathbb{R}^3$, we can assume that $C_{T_{\mathcal{B}}}$ is the origin. If a tube $T_{\theta,v}$ lies in the broom $\mathcal{B}$ associated to the pair $(s, T_{\mathcal{B}})$, then there is a point $\omega \in \theta$, such that $k_2\omega_1-k_1\omega_2=b$.   Then the direction of the core line of $T_{\theta,v}$ is parallel to $G(\omega)$  up to an angle difference of $\leq R^{-1/2}$.  We finish the proof of the lemma by observing  that  $G(\omega)$ is parallel to the plane $Z_{\mathcal{B}}$ and $T_{\theta,v}$ passes through the origin. 
\end{proof}
For example, suppose  that $Z$ is the plane $\{ x_1 =0\}$,   $T$ is the  plank defined by 
\[ T =\{  (x_1, x_2, x_3): |x_1|\leq r^{1/2+\delta}, |x_2|\leq R^{1/2+\delta}, |x_3|\leq r \},\]
and $s$ is the strip defined by $s= N_{R^{-1/2} }l \cap \tau$ for 
\[ l= \{ (\xi_1, \xi_2): \xi_2 = b \} \text{ with }  0< b< r^{-1/2}.\]
Then the broom $\mathcal{B}$ associated to $(s, T)$ consists of  tubes  $T_{\theta,v}$  of the form
\[ T_{\theta,v}= \{(x', x_3) \in B(0, R),  \, |x' + 2x_3(\omega_{\theta, 1}, b)|\leq R^{1/2+\delta}\},\]
where $\omega_{\theta, 1}$ is the first coordinate of $\omega_{\theta}$. 
And the tubes $T_{\theta,v}\in \mathcal{B}$ lie in  the $R^{1/2+\delta}$--neighborhood of 
\[ Z_{\mathcal{B}} = \{ (x_1, x_2, x_3):  x_2+ 2x_3b =0 \}. \]

\subsection{ Properties of brooms.}
\begin{lemma}\label{lem: br}
	Let $\mathcal{B}$ be a broom defined in Definition~\ref{def: broom} and $g$ be a function concentrated on wave packets from $\mathbb{T}\subset \mathcal{B}$. Write $b= | \mathbb{T}|$ as the number of tubes in $\mathbb{T}$.  Then
	\begin{equation}
	\|Eg\|_{L^2(2T_{\mathcal{B}})}^2 \lesssim b (\frac{R}{r})^{-1/2} R^{\delta}\|Eg\|_{L^2(10 B_r)}^2,
	\end{equation}
	where $B_r$ is a ball containing $T_{\mathcal{B}}$. 
\end{lemma}
\begin{proof}
	For any $T_{\theta,v}\in \mathcal{B}$, by Lemma~\ref{lem: lc1} and Lemma~\ref{lem: lc2},
	\[ 
	\int_{2T_{\mathcal{B}}} |Eg_{\theta,v}|^2 \lesssim  \frac{|2T_{\mathcal{B}}|}{|T_{\theta,v}|} \int |Eg_{\theta,v}|^2 w_{T_{\theta,v}}, 
	\]
	where $w_{T_{\theta,v}}$ is a rapidly decaying function with $w_{T_{\theta,v}}(x) \leq C_N ( 1+ n(x, T_{\theta,v}))^{-N} $  for any $N\in \mathbb{N}$. 
	By Lemma~\ref{lem: wpsupp}, the function $Eg_{\theta,v}$ is essentially supported on $T_{\theta,v}$. Hence,
	\[ \int |Eg_{\theta,v}|^2 w_{T_{\theta,v}} \lesssim \int_{10B_R} |Eg_{\theta,v}|^2 \sim R \|g_{\theta, v}\|_{L^2}^2, \]
	where the last inequality is due to Lemma~\ref{lem: l2}. Therefore, 
	\[
	\int_{2T_{\mathcal{B}}} |Eg_{\theta,v}|^2 \lesssim  \frac{|2T_{\mathcal{B}}|}{|T_{\theta,v}|} R \|g_{\theta, v}\|_{L^2}^2. 
	\]
	
	We apply the Cauchy-Schwarz inequality, \eqref{eq: l2orthogonality} and Lemma~\ref{lem: l2} using the fact that $r\geq R^{1/2+\delta}$ to obtain
	\begin{align*}
	\int_{2T_{\mathcal{B}}} |Eg|^2 & =\int_{2T_{\mathcal{B}}} |\sum_{T_{\theta,v}\in \mathbb{T}}Eg_{\theta,v}|^2\\
	&\leq b \sum_{T_{\theta,v}\in \mathbb{T}} \int_{2T_{\mathcal{B}}} |Eg_{\theta,v}|^2 \\
	&\lesssim b R \frac{|2T_{\mathcal{B}}|}{|T_{\theta,v}|} \sum_{T_{\theta,v}\in \mathbb{T}} \|g_{\theta,v}\|_{L^2}^2 \\
	& \sim b  \frac{|2T_{\mathcal{B}}|}{|T_{\theta,v}|}  \frac{R}{r} \|Eg\|_{L^2(10B_r)}^2\\
	&\lesssim b(\frac{R}{r})^{-1/2} R^{\delta} \|Eg\|_{L^2 (10 B_r)}^2. 
	\end{align*} 
\end{proof}

Lemma~\ref{lem: br} says that if a function $g$ is concentrated on wave packets from  $\mathbb{T}$, which is a subset of a single broom $\mathcal{B}$,  then we have some quantitative bound on $\|Eg\|_{L^2(T_{\mathcal{B}})}$ depending on $|\mathbb{T}|$.

 The following lemma is the motivation for defining brooms.

\begin{lemma}\label{lem: l2rescale} 
Suppose that $T$ is a plank. 	Let $\mathbb{T}_{\tau}$ be a set of tubes $T_{\theta,v}$ such that $T_{\theta,v}\cap T\neq \emptyset$ and  $\theta\subset \tau$. If $f$ is a function concentrated on wave packets from $\mathbb{T}_{\tau}$, and we set $f_s \coloneqq\underset{(\theta,v): \theta \subset s}{\sum} f_{\theta,v}$, then 
	\begin{equation}\label{eq: l2rescale}
	\int_T |Ef|^2 \lesssim  \sum_{s\subset \tau} \int |Ef_s|^2  w_{T}.
	\end{equation}
	Here 	$w_{T}$ is a positive function equal to $1$ on $T$ and rapidly decaying outside of $T$: $w_{T}(x) \leq C_N (1+ n(x, T))^{-N}$ for any $N\in \mathbb{N}$. 
\end{lemma}
\begin{proof}
Lemma~\ref{lem: l2rescale} follows from Lemma~\ref{lem: l2loc} by an affine change of variables. More precisely, let $q_s$ be a rectangular box containing $\{ (\omega, |\omega|^2): \omega\in s\}$ with  dimensions approximately $r^{-1}\times r^{-1/2}\times R^{-1/2}$.  
Let $A$ be the affine transform that maps $q_s$ to a ball of radius $r^{-1}$.  Then the affine transformation $A$ satisfies 
\begin{itemize}
	\item $A$ is independent of $s\subset \tau$ and maps each $s$ to a ball of radius $r^{-1}$, 
	\item $(A^{-1})^t$ maps $T$ to a box of dimensions $r\times rR^{\delta}\times r^{1+\delta}$, which is approximately a union of finitely overlapping  balls of radius $r$. 
\end{itemize} 
In the definition of the strip  $s$, the orientation of $s$ is chosen according to the plank  $T$ so that we can find an affine transformation $A$ satisfying  the two properties above. 
Now we apply Lemma~\ref{lem: l2loc} with $g_q((A^{-1})^t x) = Ef_s(x)$ to conclude the proof. 
\end{proof}

Define $f^{\mathcal{B}}  \coloneqq  \sum_{T_{\theta,v}\in \mathcal{B}} f_{\theta,v}$ (which is a special case of $f^{\mathbb{T}}$ with $\mathbb{T}=\mathcal{B}$), 
the right-hand side of \eqref{eq: l2rescale} can be rewritten as 
\begin{equation}\label{eq: inteprete}
\sum_{s\subset \tau} \int |Ef_s|^2  w_{T} = \sum_{\mathcal{B}: T_{\mathcal{B}}=T} \int |Ef^{\mathcal{B}}|^2 w_{T_{\mathcal{B}}}. 
\end{equation}
\begin{remark}\label{rem: decbroom}
	Let $\Sigma$ be a fat $r$--plane, and let  $\mathbb{T}_{\tau}$ be a set of tubes $T_{\theta,v}$ such that $\theta\subset \tau$ and $T_{\theta,v}\cap \Sigma \neq \emptyset$. 
If $f$ is a function concentrated on wave packets from $\mathbb{T}_{\tau}$, then by \eqref{eq: SigmaT},  Lemma~\ref{lem: l2rescale} and \eqref{eq: inteprete}, we have
\begin{equation}\label{eq: decbroom}
\int_{\Sigma} |Ef|^2 \lesssim  \sum_{\mathcal{B}: T_{\mathcal{B}}\subset  \Sigma} \int |Ef^{\mathcal{B}}|^2 w_{T_{\mathcal{B}}}. 
\end{equation}
Here the notation $\sum_{\mathcal{B}: T_{\mathcal{B}}\subset \Sigma}  \coloneqq  \sum_{T\subset \Sigma} \sum_{\mathcal{B}: T_{\mathcal{B}}=T}$ means summing over the planks $T$ in \eqref{eq: SigmaT}, then summing over the brooms $\mathcal{B}$ such that $T_{\mathcal{B}}=T$ for each plank $T$. 
\end{remark}

\begin{lemma}\label{lem: brestimate}
	Let $\tau$ be  a cap of radius $r^{1/2}$, and let $\Sigma$ be a fat $r$--plane such that $\mathit{Angle} (G(\omega_{\tau}), Z_{\Sigma})\leq r^{-1/2}$. Let $\mathbb{T}_{\tau}$ be a set of tubes $T_{\theta,v}$ with $\theta\subset \tau$ such that 
		for any broom $\mathcal{B}$ rooted at $\Sigma$, $$ |(\mathbb{T}_{\tau}\cap \mathcal{B})|\leq b.$$
	If $f$ is concentrated on wave packets from $\mathbb{T}_{\tau}$, then 
	\begin{equation}
\int_{\Sigma} |Ef|^2 \lesssim b (\frac{R}{r})^{-1/2}R^{\delta}  \int_{10B_r} |Ef|^2 +\mathit{RapDec}(R)\|h\|_{L^2}^2. 
\end{equation}
\end{lemma}
\begin{proof}
	By Remark~\ref{rem: decbroom} and  the fact that  $w_{T_{\mathcal{B}}}$ is rapidly decaying outside $T_{\mathcal{B}}$, 
	\begin{align*}
	\int_{\Sigma} |Ef|^2 &\lesssim \sum_{\mathcal{B}: T_{\mathcal{B}}\subset \Sigma} \int |Ef^{\mathcal{B}}|^2 w_{T_{\mathcal{B}}} \\
	&\lesssim \sum_{\mathcal{B}: T_{\mathcal{B}}\subset \Sigma} \int_{2T_{\mathcal{B}}}  |Ef^{\mathcal{B}}|^2 + \mathit{RapDec}(R)\|f\|_{L^2}^2. 
	\end{align*}
	By Lemma~\ref{lem: br}, 
	\begin{equation}
	\|Ef^{\mathcal{B}}\|_{L^2(2T_{\mathcal{B}})}^2 \lesssim b(\frac{R}{r})^{-1/2} R^{\delta}  \|Ef^{\mathcal{B}}\|_{L^2(10 B_r)}^2. 
	\end{equation}
	Using Lemma~\ref{lem: l2} and inequality~\eqref{eq: l2orthogonality}, we have 
	\begin{align*}
	\int_{\Sigma} |Ef|^2 &\lesssim b (\frac{R}{r})^{-1/2} R^{\delta} \sum_{\mathcal{B}: T_{\mathcal{B}}\subset \Sigma} r\|f^{\mathcal{B}}\|_{L^2}^2 +\mathit{RapDec}(R)\|f\|_{L^2}^2 \\
		 &\lesssim  b (\frac{R}{r})^{-1/2} R^{\delta} \int_{10B_r} |Ef|^2 +\mathit{RapDec}(R)\|f\|_{L^2}^2. 
	\end{align*}
\end{proof}
\subsection{Bushes}
In this section, assume that $R^{\delta}\leq r\leq R^{1/2+\delta}$.  So a fat $r$--plane $\Sigma$ can be contained in a large tube $T_{\theta,v}$. 

\begin{definition}\label{def: bush}
	For a fat $r$--plane $\Sigma$ and a cap $\tau$ of radius $r^{-1/2}$ such that $\mathit{Angle}(G(\omega_{\tau}), Z_{\Sigma})\leq r^{-1/2}$, we define the bush $\mathcal{U}$ associated to $(\tau, \Sigma)$ as 
	\begin{equation}
	\mathcal{U}: =\{ T_{\theta,v}: \Sigma \subset 2T_{\theta,v}  \text{ and }  \theta\subset \tau \}.
	\end{equation}
	We say that the bush $\mathcal{U}$ is rooted at $\Sigma$. 
\end{definition}

\begin{lemma}\label{lem: bush}
	Let $\mathcal{U}$ be a bush associated to a pair $(\tau, \Sigma)$ defined in Definition~\ref{def: bush}, and let $g$ be a function concentrated on wave packets from $\mathbb{T}\subset \mathcal{U}$. Write $u=|\mathbb{T}|$, the number of tubes in $\mathbb{T}$, then
	\begin{equation}
	\|Eg\|_{L^2 (\Sigma)}^2 \lesssim r^{-1/2+\delta} u \|Eg\|_{L^2(10B_r)}^2. 
	\end{equation}
\end{lemma}
\begin{proof}
	The bush $\mathcal{U}$ is associated to the pair $(\tau, \Sigma)$ with $\Sigma = N_{r^{1/2+\delta}}Z_{\Sigma} \cap B_{r}$. 
		We do wave packet decomposition of $g$ on a ball $B_{r^{2 -4\delta }}$ containing $\Sigma$:
	\[ g = \sum_{(\theta', v'): d(\theta')= r^{-1 +2\delta}} g_{\theta', v'} +\mathit{RapDec}(r)\|g\|_{L^2}.\]
	
	Let $\mathbb{T}'$ be the set of of tubes $T_{\theta',v'}$ of radius $r^{1- 4\delta^2}$ and length $r^{2-4\delta}$ such that $T_{\theta',v'}\cap \Sigma\neq \emptyset$, and there exists $T_{\theta,v}\in \mathbb{T}$ with $  \theta\subset 2\theta'$. 
Then
	 $|\mathbb{T}'| \lesssim |\mathbb{T}|  = u $.  
	We apply  Lemma~\ref{lem: wpsupp},  Lemma~\ref{lem: slwp} and the Cauchy--Schwarz inequality,
	\begin{align*}
	\int_{\Sigma} |Eg|^2 
	&\leq  \int_{\Sigma} | \sum_{(\theta',v') \in \mathbb{T}'} Eg_{\theta',v'}|^2 +\mathit{RapDec}(r)\|g\|_{L^2}^2 \\ 
	&\lesssim u  \int_{\Sigma} \sum_{(\theta',v')} |Eg_{\theta',v'}|^2 +\mathit{RapDec}(r)\|g\|_{L^2}^2.
	\end{align*}

	Then we apply Lemma~\ref{lem: lc1} and Lemma~\ref{lem: lc2} for the locally constant property of $Eg_{\theta', v'}$. Finally, we apply Lemma~\ref{lem: wpsupp}, \eqref{eq: l2orthogonality} and Lemma~\ref{lem: l2} as in the proof of  Lemma~\ref{lem: br},
	\begin{align*}
	\int_{\Sigma} \sum_{(\theta', v') \in \mathbb{T}'} |Eg_{\theta',v'}|^2 &\lesssim \frac{|\Sigma|}{| T_{\theta',v'}|} r^{2-4\delta}  \sum_{ (\theta',v') \in \mathbb{T}'} \|g_{\theta',v'}\|_{L^2}^2 \\
	&\sim r^{-1/2 +\delta}  \int_{10B_r} |Eg|^2.
	\end{align*}
	Since $\mathit{RapDec}(r)\|g\|_{L^2}^2$  is small compared to $r^{-1/2+\delta} u\|g\|_{L^2}^2$, we can remove it from the previous inequality. 
\end{proof}

\section{The relation $B_k \sim T_{\theta,v}$}\label{section: sim}
In this section, we define the relations $B_k\sim T_{\theta,v}$ and $B_k\nsim T_{\theta,v}$, which were used in the two ends argument in  Subsection~\ref{subsection: te}. 

Recall that we take a finitely overlapping cover of $B_{R}$ using the balls $B_k$  of radius $R^{1-\epsilon_0}$. 
Let $1\leq t\leq n$ such that $R_{j_t-1} \leq R^{1-\epsilon_0}$.  To ease the notation, we write $r=R_{j_t-1}$, \,  $S = S_t$ and $\mathcal{S}=\mathcal{S}_t$.  By Proposition~\ref{lem: plane}, for each $r^{-1/2}$--cap $\tau$, 
\begin{equation}
\|Ef^{\nsim}_{\Pi_{S}, \tau}\|_{L^2(S)}^2 \lesssim \sum_{ \Sigma\in \Omega(S, \tau)} \|Ef^{\nsim}_{\Pi_{S}, \tau}\|_{L^2(\Sigma)}^2 +\mathit{RapDec}(R)\|f\|_{L^2}, 
\end{equation}
where $\Omega(S, \tau)$  is the set of fat $r$--planes given by Proposition~\ref{lem: plane}. 
Define 
$$\Omega(\tau) \coloneqq \bigcup_{S\in \mathcal{S}} \Omega(S, \tau). $$ 
We decompose  the unit sphere in $\mathbb{R}^3$ into a union of caps $\alpha$ of radius $1/100$. This decomposition will be useful later for a geometric argument in  Lemma~\ref{lem: geoob}.
Let $\Gamma_b$ denote the set of numbers $\{ 1, (\frac{R}{r})^{\delta}, \dots, (\frac{R}{r})^{M\delta}\}$ for an integer  $M\sim \lceil \frac{1}{\delta}\rceil$. 
By Lemma~\ref{lem: brsize}, the number of tubes in a broom is $\lesssim (\frac{R}{r})^{1/2}$. The number of tubes in a bush is $\lesssim (\frac{R}{r})$ for similar reasons. 
Let $\Gamma_{\gamma}$ denote the set of numbers $\{ 1, (\frac{R}{r})^{\delta}, \dots, (\frac{R}{r})^{N\delta}\}$ for $N = \lceil 3\delta^{-1}\epsilon_0^{-1} \rceil \lesssim \delta^{-2} $. Since $r\leq R^{1-\epsilon_0}$, we have $(\frac{R}{r})^{ N\delta} \geq R^3$.  

Recall that $f$ is concentrated on wave packets from $\mathbb{T}_0$, where each wave packet $f_{\theta,v}$  from $\mathbb{T}_0$ has about the same $L^2$--norm (see Subsection~\ref{subsection: reduction}). To describe the relation $B_k\sim T_{\theta,v}$, we are going to define a family of  functions by induction  on the space $$\mathbb{T}_0\times \Omega(\tau).$$
To explain the initial step,  for each cap $\alpha$, each  $b_1\in \Gamma_b$ and $\Sigma \in \Omega(\tau)$, we define $\chi_{\alpha, t, b_1}(T_{\theta,v}, \Sigma) =1$ if $T_{\theta,v}$ belongs to a broom $\mathcal{B}$ when $ r\geq R^{1/2+\delta}$ (or a bush $\mathcal{U}$ when $r \leq R^{1/2+\delta}$) rooted at $\Sigma$ such that 
\begin{enumerate}
	\item the normal direction of $Z_{\Sigma}$ lies in $\alpha$;
	\item  $b_1 \leq |\mathbb{T}_0\cap \mathcal{B}|< (\frac{R}{r})^{\delta} b_1$ when $r\geq R^{1/2+\delta}$ (or  $b_1 \leq |\mathbb{T}_0\cap \mathcal{U}|< (\frac{R}{r})^{\delta} b_1$ when $r\leq R^{1/2 +\delta}$).
\end{enumerate}
Otherwise, we define $\chi_{\alpha, t, b_1}(T_{\theta,v}, \Sigma)=0$. 
Here we use $t$ in $\chi_{\alpha, t, b_1}$ to record that $r=R_{j_t-1}, \, S=S_t$ and $\mathcal{S}=\mathcal{S}_t$.  Note that if $\theta\subset \tau$ and $T_{\theta,v} \cap \Sigma\neq \emptyset$, we can always find a cap $\alpha$ and $b_1\in \Gamma_b$ such that $\chi_{\alpha, t, b_1}(T_{\theta,v}, \Sigma)=1$. 

Now we proceed to the second step. For each $\gamma_1\in \Gamma_{\gamma}$, we define $\chi_{\alpha, t, b_1, \gamma_1}(T_{\theta,v}, \Sigma)=1$ if 
\begin{enumerate}
	\item $\chi_{\alpha, t, b_1}(T_{\theta,v}, \Sigma)=1$;
	\item  the following estimate holds $$\gamma_1 \leq \sum_{\Sigma'\in\Omega(\tau)} \chi_{\alpha, t, b_1}(T_{\theta,v}, \Sigma') < (\frac{R}{r})^{\delta} \gamma_1.$$
\end{enumerate}
Otherwise, we define $\chi_{\alpha, t, b_1, \gamma_1}(T_{\theta,v}, \Sigma)=0$. Since  $\Omega(\tau)  \lesssim R^{O(\delta)}  |\mathcal{S}|\leq R^3$, the summation  holds \begin{equation}\label{eq: exgamma}
\underset{\Sigma\in \Omega(\tau)}{\sum} \chi_{\alpha, t, b_1} (T_{\theta,v}, \Sigma) \leq R^3.
\end{equation}
It follows that  if $\chi_{\alpha, t, b_1}(T_{\theta,v}, \Sigma)=1$, then we can always find $\gamma_1\in \Gamma_{\gamma}$  such that $\chi_{\alpha, t, b_1, \gamma_1}(T_{\theta,v}, \Sigma)=1$. 

As for the induction step, let $\kappa=(\alpha, t, b_1, \gamma_1, \dots, b_l, \gamma_l)$ and $\kappa' = (\alpha, t, b_1, \gamma_1, \dots, b_l)$ for some integer $2\leq l \leq 100 \delta^{-2}\epsilon_0^{-2}$. We define functions $\chi_{\kappa}$ and $\chi_{\kappa'}$ inductively as above until the condition that 
\begin{equation}\label{eq: stop}
b_l \geq (\frac{R}{r})^{-\delta} b_{l-1} \text{ and } \gamma_{l} \geq (\frac{R}{r})^{-\delta} \gamma_{l-1}
\end{equation}
is fulfilled. 

More precisely, we define $\chi_{\kappa}(T_{\theta,v}, \Sigma)=1$ if 
\begin{enumerate}
	\item $\chi_{\kappa'}(T_{\theta,v}, \Sigma) =1$;
	\item the following estimate holds
	\[ \gamma_l \leq \sum_{ \Sigma'\in \Omega(\tau)}\chi_{\kappa'}(T_{\theta,v}, \Sigma') < \gamma_l (\frac{R}{r})^{\delta}.\] 
\end{enumerate}
Otherwise, we define $\chi_{\kappa}(T_{\theta,v}, \Sigma) =0$. 
For each $b_{l+1}$, we define $\chi_{\alpha, t, b_1, \gamma_1, \dots, b_l, \gamma_l, b_{l+1}}(T_{\theta,v}, \Sigma)=1$ if 
\begin{enumerate}
	\item $\chi_{\kappa}(T_{\theta,v}, \Sigma)=1$;
	\item for the broom $\mathcal{B}$  (or the bush $\mathcal{U}$ when $r\leq R^{1/2+\delta}$) containing $T_{\theta,v}$ rooted at $\Sigma$, 
	$$b_{l+1}\leq \sum_{T_{\theta',v'}\in \mathcal{B}} \chi_{\kappa}(T_{\theta',v'}, \Sigma) < b_{l+1} (\frac{R}{r})^{\delta}. $$
\end{enumerate}
From our definition of $\chi_{\kappa}$ and $\chi_{\kappa'}$, only when $b_1\geq b_2\geq \cdots \geq b_l$ and $\gamma_1\geq \cdots \geq \gamma_l$ the functions $\chi_{\kappa}$ and $\chi_{\kappa'}$ can be nontrivial (not identically zero).  There are  $\lesssim_{\delta} 1$ vectors $\kappa$ and $\kappa'$ satisfying \eqref{eq: stop}.  For any $T_{\theta,v} \cap \Sigma \neq \emptyset$ and $\theta\subset \tau$, there exists a unique  $\kappa$  satisfying \eqref{eq: stop} such that $\chi_{\kappa}(T_{\theta,v}, \Sigma)=1$.   To see this, if it does not exist  $\kappa$ satisfying \eqref{eq: stop} with $\chi_{\kappa}(T_{\theta,v}, \Sigma)=1$, we would have a sequence of pairs $(b_l, \gamma_l)$ such that $l \geq 100 \delta^{-2} \epsilon_0^{-2} $, and either
\begin{itemize}
	\item $ (\frac{R}{r})^{\delta} \leq  b_{l-1}/b_l $  and $\gamma_l \leq \gamma_{l-1}$,  or
	\item $b_l \leq b_{l-1}$ and $ (\frac{R}{r})^{\delta} \leq \gamma_{l-1}/\gamma_l$.
\end{itemize}
This is impossible by pigeonholing for  $(b_l, \gamma_l)\in \Gamma_b\times \Gamma_{\gamma}$. 

\begin{definition}\label{def: simkappa}
	For each $\kappa$ and each tube $T_{\theta,v}$, let $B_k^{*}$ be a ball of radius $R^{1-\epsilon_0}$ that maximizes the quantity
	\[ \sum_{\Sigma\in \Omega(\tau), \Sigma\subset B_k^{*}} \chi_{\kappa} (T_{\theta,v}, \Sigma).\]
	If there are multiple maximizers, then we choose only one. We define $B_k\sim_{\kappa} T_{\theta,v}$ if $ B_k\subset 10 B_{k}^{*}$. Otherwise, we define  $B_k \nsim_{\kappa} T_{\theta,v} $.  We define $B_{k}\sim_{\kappa'} T_{\theta,v}$ and $B_k\nsim_{\kappa'} T_{\theta,v}$ according to the same rule with $\kappa'$ in place of $\kappa$. 
\end{definition}

\begin{definition}\label{def: sim}
	We define $B_k\sim T_{\theta,v}$ if $B_{k} \sim_{\kappa} T_{\theta,v}$ or $B_k \sim_{\kappa'}$ for any $\kappa $ or  $\kappa'$  satisfying \eqref{eq: stop}. We define $B_k\nsim T_{\theta,v}$ if   $B_{k}\nsim_{\kappa} T_{\theta,v} $ and $B_{k} \nsim_{\kappa'} T_{\theta,v}$ for all  $\kappa$ and $\kappa'$ satisfying~\eqref{eq: stop}. 
\end{definition}

\begin{lemma}
	For each tube $T_{\theta,v}$, the number of balls $B_k$ such that $B_k\sim T_{\theta,v}$ is $\lesssim_{\delta} 1$. 
\end{lemma}
\begin{proof}
	By the definition of $\sim_{\kappa}$, for each $T_{\theta,v}$, the number of $B_k$ such that $B_{k}\sim_{\kappa} T_{\theta,v}$ is $\lesssim 1$. Similarly, the number of $B_k$ such that $B_k\sim_{\kappa'} T_{\theta,v}$ is $\lesssim 1$. Since there are $\lesssim_{\delta} 1$ choices for $\kappa$ and $\kappa'$, the number of $B_k$ such that $B_k\sim T_{\theta,v}$ is $\lesssim_{\delta}1$. 
\end{proof}
\section{Estimating the global part $Ef^{\nsim}$}\label{section: l2}

Let $\mathcal{O}$ be the tree defined  in  Corollary~\ref{cor: tree}. In this section, assume that 
\begin{itemize}
\item $\|Ef\|_{BL^p(O)} \lesssim \|Ef^{\nsim}\|_{BL^p(O)}$ holds for a $\geq (1-R^{-\delta})$--fraction of the leaves $O\in \mathcal{O}$;
\item 
case \ref{case2} of  Lemma~\ref{lem: findscales} happens   for some $t$ with $1\leq t\leq n$ and $R_{j_t-1}\leq R^{1-\epsilon_0}$, which implies that \eqref{eq: twin} in Lemma~\ref{lem: find scale} holds. 
\end{itemize}

Our goal is to estimate the global part $Ef^{\nsim}$. The main results are Lemma~\ref{lem: keylarge2} and Lemma~\ref{lem: keysmall2}.

To ease the notation, we write $r= R_{j_t-1}$, \, $S=S_t$ and $\mathcal{S}=\mathcal{S}_t$. 
We can assume that each small tube $T_{\tau, w} \in \mathbb{T}_{S}$ belongs to only one $\mathbb{T}_{\Sigma}$ with  $\Sigma \in \Omega(S, \tau)$. Otherwise, we assign $T_{\tau,w}$ to only one $\mathbb{T}_{\Sigma}$.

For each fat $r$--surface $S\in \mathcal{S}$,  let $\tau(S)$ denote a  $r^{-1/2}$--cap that attains the maximum
\begin{equation}\label{eq: tauS}
\max_{d(\tau)=r^{-1/2}} \|f^{\nsim}_{\Pi_S, \tau} \|_{L^2_{avg}(\tau)}. 
\end{equation}
If there are multiple choices for $\tau(S)$, we choose only one. 
 For each  $\kappa$ satisfying the condition \eqref{eq: stop}, each $r$--fat plane $\Sigma \in \Omega(S, \tau(S)) $ and the ball $B_k$ containing $\Sigma$, we define  the set of tubes\footnote{Recall that $\mathbb{T}_0$ is defined as  in Subsection~\ref{subsection: reduction}.}
 \begin{equation}\label{eq: defsigmatu}
 \mathbb{T}^{\nsim}_{\kappa}(\Sigma) \coloneqq  \{ T_{\theta,v}\in \mathbb{T}_0: \chi_{\kappa}(T_{\theta,v}, \Sigma)=1, B_k\nsim T_{\theta,v}\}, 
 \end{equation}
  and define 
 \begin{equation}\label{eq: defsigma}
 f^{\nsim}_{\kappa, \Sigma} : =  \sum_{T_{\theta,v} \in \mathbb{T}^{\nsim}_{\kappa}(\Sigma)} f_{\theta,v},
 \end{equation}
 and
 \begin{equation}\label{eq: defsigmat}
  f^{\nsim}_{\kappa, \Sigma, \mathit{tang}} : = \sum_{w: T_{\tau,w}\in \mathbb{T}_{\Sigma}} ( f^{\nsim}_{\kappa, \Sigma} )_{\tau, w}.
 \end{equation}
 For the ball $B_k$ containing $S$, we apply  Lemma~\ref{lem: properties} with $g=f$ and  $\mathbb{T}= \mathbb{T}_k^{\nsim}= \{ T_{\theta,v}: B_k\nsim T_{\theta,v}\}$ to define
 \begin{equation}
 f^{\nsim}_{\Pi_{S}} = \sum_{T_{\tau,w} \in \mathbb{T}_S}  (\sum_{(\theta,v): T_{\theta,v}\cap S\neq \emptyset, B_k\nsim T_{\theta,v}} f_{\theta,v})_{\tau, w}.
\end{equation}
Hence, for $\tau=\tau(S)$, 
\begin{equation}\label{eq: deck}
f^{\nsim}_{\Pi_S, \tau}  =  \sum_{w: T_{\tau,w} \in \mathbb{T}_S}  (\sum_{(\theta,v): T_{\theta,v}\cap S\neq \emptyset, B_k\nsim T_{\theta,v}} f_{\theta,v})_{\tau, w}= \sum_{\Sigma\in \Omega(S, \tau)} \sum_{\kappa} f^{\nsim}_{\kappa, \Sigma, \mathit{tang}} + \mathit{RapDec}(R)\|f\|_{L^2},
\end{equation}
where  $\sum_{\kappa}$ means  summing over all $\kappa$ satisfying \eqref{eq: stop}. 
Recall that  the number of such  $\kappa$ is $\lesssim_{\delta} 1$ and $|\Omega(S, \tau)|\lesssim r^{O(\delta)}$.  For each $S\in \mathcal{S}$ and $\tau=\tau(S)$, by the triangle inequality and \eqref{eq: deck}, there exist some $\Sigma\in \Omega(S, \tau(S))$ and  $\kappa$ such that 
	\begin{equation}\label{eq: choosekappa}
\|f^{\nsim}_{\Pi_S, \tau(S)}\|_{L^2(\tau(S))} \lesssim_{\delta} r^{O(\delta)} \|f^{\nsim}_{\kappa, \Sigma, \mathit{tang}} \|_{L^2(\tau(S))}.
\end{equation}
We decompose $\mathcal{S}=\cup \mathcal{S}(\kappa)$, where $\mathcal{S}(\kappa)$ denotes the set of $S\in\mathcal{S}$ such that \eqref{eq: choosekappa} holds for $\kappa$ with some $\Sigma\in \Omega(S, \tau(S))$.
Hence, 
\begin{equation}
\sum_{S\in \mathcal{S}} \sum_{O\in \mathcal{O}, O\subset S} \|Ef^{\nsim}_{\Pi_S}\|_{BL^p(O)}^p \leq \sum_{\kappa} \sum_{S\in \mathcal{S}(\kappa)} \sum_{O\in \mathcal{O}, O\subset S} \|Ef^{\nsim}_{\Pi_S}\|_{BL^p(O)}^p. 
\end{equation}
By pigeonholing, there exists  $\kappa$ with  $\mathcal{S}(\kappa) \subset \mathcal{S}$ such that 
\begin{itemize}
	\item for each $S\in \mathcal{S}(\kappa)$, there exists a fat $r$--plane $\Sigma\in \Omega(S, \tau(S))$ such that 
\eqref{eq: choosekappa} holds, 
	\item \begin{equation}\label{eq: chooseS}
\sum_{S\in \mathcal{S}}  \sum_{O\in\mathcal{O}, O\subset S}	\|Ef^{\nsim}_{\Pi_S}\|_{BL^p(O)}^p \lesssim_{\delta}  \sum_{S\in \mathcal{S}(\kappa)}  \sum_{O\in \mathcal{O}, O\subset S}	\|Ef^{\nsim}_{\Pi_S}\|_{BL^p(O)}^p. 
	\end{equation}
\end{itemize}

We separate the estimate of the global part $Ef^{\nsim}$ in three cases according to the size of $r$. The case when $r\geq R^{1-\epsilon_0}$ has been treated in Lemma~\ref{lem: large}. In this section, we focus on the remaining two cases,  which are more difficult:  the $R^{1/2+\delta}\leq r \leq R^{1-\epsilon_0}$ case and the $R^{\delta}\leq r\leq R^{1/2+\delta}$ case. 

\subsection{The  $R^{1/2+\delta}\leq r\leq R^{1-\epsilon_0}$ case.}

We start with Lemma~\ref{lem: geoob} for the geometric property of a broom. 

\begin{lemma}\label{lem: geoob}
	Let $\Sigma_1$ and $\Sigma_2$ be fat $r$--planes such that
	\begin{enumerate}
		\item 	$\mathit{dist}(\Sigma_1, \Sigma_2)\geq R^{1-\epsilon_0}$, 
		\item  the normal directions of $Z_{\Sigma_1}$ and $Z_{\Sigma_2}$ form an angle of $\leq 1/100$. 
	\end{enumerate}
	If $\mathcal{B}$ is a broom rooted at $\Sigma_2$, then there are $\lesssim R^{O(\epsilon_0)}$ tubes in $\mathcal{B}$ intersecting $\Sigma_1$. 
\end{lemma}
\begin{proof}
	By Lemma~\ref{remark: brplane},  for any broom $\mathcal{B}$ rooted at $\Sigma_2$, the tubes $T_{\theta,v}\in \mathcal{B}$ lie in the $R^{1/2+\delta}$--neighborhood of a plane  $Z_{\mathcal{B}}$, and  $Z_{\mathcal{B}}$ is orthogonal to $Z_{\Sigma_2}$ up to an angle difference of $O(r^{-1/2})$.  Since the normal directions of $Z_{\Sigma_1}$ and $Z_{\Sigma_2}$ form an angle of $\leq 1/100$, the normal directions of $Z_{\Sigma_1}$ and $Z_{\mathcal{B}}$ form an angle of  $\gtrsim 1$, which implies that the volume  
	\begin{equation}\label{eq: vol}
	\mathit{Vol} (N_{2R^{1/2+\delta}} Z_{\mathcal{B}} \cap \Sigma_1) \lesssim r^{3/2+\delta} R^{1/2+\delta}.
	\end{equation} 
	Recall that $\Sigma_1$ has dimensions about $r\times r \times r^{1/2+\delta}$.
	Furthermore, for  each tube  $T_{\theta,v}\in \mathcal{B}$, $2T_{\theta,v}$ contains the plank  $T_{\mathcal{B}}$,  and $T_{\mathcal{B}}$  is $\geq R^{1-\epsilon_0}$ away from $\Sigma_1$.  We claim that at any point $x$ with $\mathit{dist} (x, \Sigma_2) \geq R^{1-\epsilon_0}$, there are $\lesssim R^{\epsilon_0}$ tubes $T_{\theta,v}\in \mathcal{B}$  such that $x\in 2T_{\theta,v}$. The reason is that if $2T_{\theta,v}$ contains both  $x$ and $T_{\mathcal{B}}$ (recall that $\mathit{dist}(x, T_{\mathcal{B}})\geq R^{1-\epsilon_0}$), then $\theta$ is contained in a small cap $\beta$ of radius $\lesssim R^{-1/2+\epsilon_0 +\delta}$ determined by $x$ and $T_{\mathcal{B}}$.   Since $\mathcal{B}$ is associated to a unique pair $(s, T_{\mathcal{B}})$, there are $\lesssim R^{O(\epsilon_0)}$ caps $\theta$ such that $\theta\subset \beta\cap s$.  Moreover,  for each $\theta$, there are $\lesssim 1$ tubes $T_{\theta,v}\in \mathcal{B}$. 
	For each tube $T_{\theta,v}\in\mathcal{B}$, if $T_{\theta,v}\cap \Sigma_1\neq \emptyset$, then \begin{equation}\label{eq: vol2}
	\mathit{Vol} (2T_{\theta,v} \cap \Sigma_1)\gtrsim  r^{3/2+\delta} R^{1/2+\delta}.
	\end{equation}
	Combining \eqref{eq: vol},  \eqref{eq: vol2} and our discussion above,   there are $\lesssim R^{O(\epsilon_0)}$ tubes from $\mathcal{B}$ intersecting $\Sigma_1$. 
	
\end{proof}

Then we apply Lemma~\ref{lem: geoob} to prove Lemma~\ref{lem: keylarge}.

\begin{lemma}\label{lem: keylarge}
	If $R^{1/2+\delta}\leq r\leq R^{1-\epsilon_0}$, $\tau$ is a cap of radius $r^{-1/2}$ and  $\Sigma\in \Omega(\tau)$, then for each $\kappa = (\alpha, t, b_1,\gamma_1, \dots, b_l, \gamma_l)$ satisfying~\eqref{eq: stop}, 
	\begin{equation}
	\|f^{\nsim}_{\kappa, \Sigma, \mathit{tang}} \|_{L^2(\tau)}^2 \lesssim (\frac{R}{r})^{-1/2} R^{O(\epsilon_0)} \|f_{\tau}\|_{L^2}^2. 
	\end{equation}
\end{lemma}
\begin{proof}
	Let $B_k$ be the ball of radius $R^{1-\epsilon_0}$ containing $\Sigma$ and write $\Sigma_1=\Sigma$. The idea is to double count the number of wave packets shared by $\Sigma_1$ and fat $r$--planes $\Sigma_2\nsubseteq  5B_k$: specifically, the quantity
	\begin{equation}\label{eq: dc}
	\sum_{\Sigma_2\nsubseteq 5B_k} \sum_{(\theta,v): \theta\subset \tau} \chi_{\kappa} (T_{\theta,v}, \Sigma_1) \chi_{\kappa'}(T_{\theta,v}, \Sigma_2),
	\end{equation}
	where  $\sum_{\Sigma_2\nsubseteq 5B_k}$ means summing over all $\Sigma_2\in \Omega(\tau)$ such that $\Sigma_2\nsubseteq 5B_k$ and $\sum_{(\theta,v): \theta\subset \tau}$ means summing over all  tubes $T_{\theta,v}\in \mathbb{T}_0$ with  $\theta\subset \tau$.  Recall that $\Omega(\tau) = \cup_{S\in \mathcal{S}} \Omega(S, \tau)$. 
	
	The definition of $B_k\nsim T_{\theta,v}$ implies that $B_k\nsim_{\kappa'} T_{\theta,v}$ with $\kappa'=(\alpha, t, b_1,\gamma_1, \dots, b_l)$. For each $T_{\theta,v}$, we have 
	\begin{equation}\label{eq: dctw}
	\sum_{\Sigma_2\nsubseteq 5B_k} \chi_{\kappa'} (T_{\theta,v}, \Sigma_2) \gtrsim \sum_{\Sigma'} \chi_{\kappa'}(T_{\theta,v}, \Sigma'),
	\end{equation}
	otherwise the ball $B_k^*$ that maximizes $\underset{\Sigma'\in \Omega(\tau), \Sigma \subset B_k^*}{\sum} \chi_{\kappa'}(T_{\theta,v}, \Sigma')$ should belong to $5B_k$, which yields $B_k \sim_{\kappa'} T_{\theta,v}$. This violates the assumption that $B_k\nsim T_{\theta,v}$.  Note that \eqref{eq: dctw} is the only place where  we need to use the information that $B_k\nsim T_{\theta,v}$.  For each tube $T_{\theta,v}$ satisfying $\chi_{\kappa}(T_{\theta,v}, \Sigma)=1$, we have 
	\begin{equation}\label{eq: lower}
	\sum_{\Sigma'} \chi_{\kappa'} (T_{\theta,v}, \Sigma') \geq \gamma_l. 
	\end{equation}
	Define $N$ to be the quantity $$N\coloneqq \sum_{(\theta, v): \theta\subset \tau} \chi_{\kappa}(T_{\theta,v}, \Sigma_1) .$$ We have the following lower bound for \eqref{eq: dc} by combining \eqref{eq: dctw} and \eqref{eq: lower}:
	\begin{equation}\label{eq: lowerbd}
	\sum_{\Sigma_2\nsubseteq 5B_k} \sum_{(\theta,v): \theta\subset \tau} \chi_{\kappa} (T_{\theta,v}, \Sigma_1) \chi_{\kappa'}(T_{\theta,v}, \Sigma_2) \geq \gamma_l N. 
	\end{equation}
	
	Next we are going to give an upper bound for \eqref{eq: dc}. 
Recall that for a fixed $\Sigma_2$, each tube $T_{\theta,v}$ belongs to at most one broom rooted at $\Sigma_2$. Let $\tilde{\kappa}=(\alpha, t, b_1, \dots, b_{l-1}, \gamma_{l-1})$ and $\tilde{\kappa}'=(\alpha, t, b_1, \dots, b_{l-1})$. 
If $\chi_{\kappa'}(T_{\theta,v}, \Sigma_2)=1$, then $\chi_{\tilde{\kappa}}(T_{\theta,v}, \Sigma_2)=1$ and  for the broom $\mathcal{B}$ containing $T_{\theta,v}$ rooted at $\Sigma_2$, 
$$\sum_{T_{\theta',v'}\in \mathcal{B}} \chi_{\tilde{\kappa}}(T_{\theta',v'}, \Sigma_2)\geq b_l.$$
  If $\chi_{\kappa}(T_{\theta,v}, \Sigma_1) \chi_{\kappa'} (T_{\theta,v}, \Sigma_2) =1$, 
then the normal vectors of $Z_{\Sigma_1}$ and $Z_{\Sigma_2}$ are in the same cap $\alpha$ of radius $1/100$. So we can apply Lemma~\ref{lem: geoob} for each $\Sigma_2\nsubseteq B_k$ to show that 
\begin{align}\label{eq: upper} 
\sum_{(\theta,v): \theta\subset \tau} \chi_{\kappa}(T_{\theta,v}, \Sigma_1)\chi_{\kappa'}(T_{\theta,v}, \Sigma_2) &\leq \sum_{(\theta,v): \theta\subset \tau, T_{\theta,v} \cap \Sigma_1\neq \emptyset} \chi_{\kappa'} (T_{\theta,v}, \Sigma_2)\\
& \lesssim R^{O(\epsilon_0)} b_l^{-1} \sum_{(\theta,v): \theta\subset \tau} \chi_{\tilde{\kappa}}(T_{\theta,v}, \Sigma_2). \nonumber
\end{align}

By the definition of $\chi_{\tilde{\kappa}}$, we have $$\sum_{\Sigma_2} \chi_{\tilde{\kappa}}(T_{\theta,v}, \Sigma_2) \leq  \sum_{\Sigma_2} \chi_{\tilde{\kappa}'}(T_{\theta,v}, \Sigma_2) \leq  \gamma_{l-1} (\frac{R}{r})^{\delta}. $$
Define $N_0$ to be
\[ N_0\coloneqq | \{  T_{\theta,v}\in \mathbb{T}_0: \theta\subset \tau \}|. \]
 Hence,
\begin{equation}\label{eq: upper2}
\sum_{\Sigma_2} \sum_{(\theta,v): \theta\subset\tau} \chi_{\tilde{\kappa}}(T_{\theta,v}, \Sigma_2) < \gamma_{l-1} N_0 (\frac{R}{r})^{\delta}. 
\end{equation}
Summing both sides of inequality~\eqref{eq: upper} over $\Sigma_2\nsubseteq 5B_k$ and then applying inequality~\eqref{eq: upper2}, we have the following upper bound for \eqref{eq: dc}:
\begin{equation}\label{eq: upperbd}
\sum_{\Sigma_2\nsubseteq 5B_k} \sum_{(\theta,v): \theta\subset \tau} \chi_{\kappa}(T_{\theta,v}, \Sigma_1) \chi_{\kappa'}(T_{\theta,v}, \Sigma_2) \lesssim R^{O(\epsilon_0)}  \gamma_{l-1}N_0 b_l^{-1}.
\end{equation}
Since $\kappa$ satisfies \eqref{eq: stop}, we have $\gamma_l \geq  (\frac{R}{r})^{-\delta} \gamma_{l-1}$.

We compare the lower bound \eqref{eq: lowerbd} and the upper bound \eqref{eq: upperbd} for \eqref{eq: dc} to obtain
\begin{equation}\label{eq: qq}
 N b_l \lesssim R^{O(\epsilon_0)} N_0. 
\end{equation}

To finish the proof, we apply Lemma~\ref{lem: brestimate}  for the function  $f^{\nsim}_{\kappa, \Sigma}$ defined in \eqref{eq: defsigma}, which is concentrated on wave packets from $\mathbb{T}^{\nsim}_{\kappa}(\Sigma)$.  By the definition of $\mathbb{T}_{\kappa}^{\nsim}$ in \eqref{eq: defsigmatu},  $|\mathbb{T}^{\nsim}_{\kappa}(\Sigma)\cap \mathcal{B}| \leq b_l (\frac{R}{r})^{\delta}$ for any  broom $\mathcal{B}$ rooted at $\Sigma$.  So we have 
\begin{align*}
\|Ef^{\nsim}_{\kappa, \Sigma, \mathit{tang}}\|_{L^2(10 B_r)}^2 &\lesssim \|Ef^{\nsim}_{\kappa, \Sigma, \mathit{tang}} \|_{L^2(\Sigma )}^2  +\mathit{RapDec}(R)\|f\|_{L^2}^2 \\ &\lesssim R^{2\delta} (\frac{R}{r})^{-1/2} b_l \|Ef^{\nsim}_{\kappa,\Sigma}\|_{L^2(10 B_r)}^2  +\mathit{RapDec}(R)\|f\|_{L^2}^2.
\end{align*}

By Lemma~\ref{lem: l2}, we have $\|Ef^{\nsim}_{\kappa, \Sigma, \mathit{tang}}\|_{L^2(10B_r)}^2\sim r\|f^{\nsim}_{\kappa, \Sigma, \mathit{tang}}\|_{L^2}^2$ and $\|Ef^{\nsim}_{\kappa, \Sigma}\|_{L^2(10B_r)}^2 \sim r\|f^{\nsim}_{\kappa, \Sigma}\|_{L^2}^2$, which yield 
\begin{equation}\label{eq: scte}
\|f^{\nsim}_{\kappa, \Sigma, \mathit{tang}}\|_{L^2}^2 \lesssim R^{2\delta} (\frac{R}{r})^{-1/2} b_l \|f^{\nsim}_{\kappa, \Sigma}\|_{L^2}^2 +\mathit{RapDec}(R)\|f\|_{L^2}^2.
\end{equation}
Recall that all wave packets from $\mathbb{T}_0$  have approximately the same $\|f_{\theta,v}\|_{L^2}$. Since there are $N$ out of $N_0$ tubes $T_{\theta,v}\in \mathbb{T}_0$ with $\chi_{\kappa}(T_{\theta,v}, \Sigma)=1$, we have 
\begin{equation}\label{eq: sscte}
\|f^{\nsim}_{\kappa, \Sigma}\|_{L^2}^2 \lesssim N/N_0 \|f_{\tau}\|_{L^2}^2 +\mathit{RapDec}(R)\|f\|_{L^2}^2. 
\end{equation}
Combining \eqref{eq: scte}, \eqref{eq: sscte} and \eqref{eq: qq}, we have 
\begin{equation}
\|f^{\nsim}_{\kappa, \Sigma, \mathit{tang}}\|_{L^2}^2 \lesssim R^{O(\epsilon_0) } (\frac{R}{r})^{-1/2} \|f_{\tau}\|_{L^2}^2. 
\end{equation}
We can drop the $\mathit{RapDec}(R) \|f\|_{L^2}^2$--term  because it is much smaller than $R^{O(\epsilon_0)} (\frac{R}{r})^{-1/2} \|f_{\tau}\|_{L^2}^2.$
Since $\delta <\epsilon_0$, we obtain the desired bound. 
\end{proof}

Using Lemma~\ref{lem: keylarge} we prove Lemma~\ref{lem: keylarge2}.

\begin{lemma}\label{lem: keylarge2}
	Let $\mathcal{O}$ be the tree defined in Corollary~\ref{cor: tree}.  Assume that 
	\begin{equation}\label{eq: nsim}
	\|Ef\|_{BL^p(O)}\lesssim \|Ef^{\nsim}\|_{BL^p(O)}
	\end{equation}
	 holds for a $\geq (1-R^{-\delta})$--fraction of the leaves $O\in \mathcal{O}$ and  \eqref{eq: 55}  in Lemma~\ref{lem: findscales} and \eqref{eq: twin} in Lemma~\ref{lem: find scale} hold   for some $t$ with $R^{1/2+\delta} \leq R_{j_t-1}\leq R^{1-\epsilon_0}$.  Then for $p\geq 3+3/13$, inequality~\eqref{eq: induction2} in Theorem~\ref{thm: induction2} holds for the function $f$.
\end{lemma}
\begin{proof}
	To ease the notation, write  $r=R_{j_t-1}$,\,  $S=S_t$,  \,$\mathcal{S}=\mathcal{S}_t$ and $D=d^{j_t}$. 
	By Corollary~\ref{cor: tree} and inequality~\eqref{eq: nsim},
	\begin{align*}
	\|Ef\|_{BL^p(B_R)}^p &\lesssim R^{\delta} \sum_{O \in \mathcal{O}} \|Ef\|_{BL^p(O)}^p\\
		&\lesssim R^{\delta} \sum_{O\in \mathcal{O}} \|Ef^{\nsim}\|_{BL^p(O)}^p.
		\end{align*}
		
We are going to estimate $\|Ef\|_{BL^p(B_R)}^p$ using two different approaches. 

\vspace{5pt}

\textbf{(I).} For the first approach, by inequalities \eqref{eq: 55}, \eqref{eq: twin}  and \eqref{eq: chooseS},  we have
		\begin{align*}
	\|Ef\|_{BL^p(B_R)}^p&\lesssim R^{O(\delta)} \sum_{O\in \mathcal{O}} \|Ef^{\nsim}_{\Pi_S}\|_{BL^p(O)}^{p}\\
	&\lesssim R^{O(\delta)} \sum_{S\in \mathcal{S}(\kappa)} \sum_{O\in \mathcal{O}, \, O\subset S} \|Ef^{\nsim}_{\Pi_S}\|_{BL^p(O)}^p\\
	 &\lesssim R^{O(\delta)} \sum_{S\in \mathcal{S}(\kappa)} \sum_{O\in \mathcal{O},  \, O\subset S } \|Ef^{\nsim}_S\|_{BL^p(O)}^2 \|Ef^{\nsim}_{\Pi_S}\|_{BL^p(O)}^{p-2}\\
	 &\lesssim R^{O(\delta)} \sum_{S\in \mathcal{S}(\kappa)} \|Ef^{\nsim}_S\|_{BL^p(S)}^2 \|Ef^{\nsim}_{\Pi_S}\|_{BL^p(S)}^{p-2}.
	\end{align*}
	The last inequality is due to H\"{o}lder's inequality, and the set  $\mathcal{S}(\kappa)\subset \mathcal{S}$ was chosen such that  \eqref{eq: choosekappa} holds for each $S \in \mathcal{S}(\kappa)$. 
	We apply Lemma~\ref{lem: tao} to both $\|Ef^{\nsim}_S\|_{BL^p(S)}$ and $\|Ef^{\nsim}_{\Pi_S}\|_{BL^p(S)}$ to derive
	\begin{align}\label{eq: est}
	\|Ef\|_{BL^p(B_R)}^p &\lesssim R^{O(\delta)} \sum_{S\in \mathcal{S}(\kappa)} \|Ef^{\nsim}_S\|_{BL^p(S)}^2 \|Ef^{\nsim}_{\Pi_S}\|_{BL^p(S)}^{p-2}\\& \lesssim R^{O(\delta)} r^{\frac{5}{2}-\frac{3p}{4}} \sum_{S\in \mathcal{S}(\kappa)} \|f^{\nsim}_S\|_{L^2}^2 \|f^{\nsim}_{\Pi_S}\|_{L^2}^{p-2}. \nonumber
	\end{align}
By Lemma~\ref{lem: geometric} and the definition of $\tau(S)$ (see \eqref{eq: tauS}), 
\begin{equation}\label{eq: geo}
\|f^{\nsim}_{\Pi_S}\|_{L^2}^2 \lesssim r^{-1/2+O(\delta)} \max_{d(\tau)=r^{-1/2}} \|f^{\nsim}_{\Pi_S, \tau}\|_{L^2_{avg}(\tau)}^2 = r^{-1/2+O(\delta)} \|f^{\nsim}_{\Pi_S, \tau(S)}\|_{L^2_{avg}(\tau(S))}^2. 
\end{equation}

Plugging \eqref{eq: geo}  in \eqref{eq: est}, we have 
\begin{equation}\label{eq: est01}
\|Ef\|_{BL^p(B_R)}^p \lesssim R^{O(\delta)} r^{\frac{5}{2}-\frac{3p}{4}}r^{-\frac{p-2}{4}}(\sum_{S\in \mathcal{S}(\kappa)} \|f^{\nsim}_S\|_{L^2}^2 )\cdot \max_{S\in \mathcal{S}(\kappa)} \|f^{\nsim}_{\Pi_{S}, \tau(S)}\|_{L^2_{avg}(\tau(S))}^{p-2}.
\end{equation}

We apply inequality~\eqref{eq: strcell} in Lemma~\ref{lem: properties} for $\mathcal{S}=\mathcal{S}_t$, 
\begin{equation}\label{eq: 728}
\sum_{S\in \mathcal{S}(\kappa)} \|f^{\nsim}_S\|_{L^2}^2\leq \sum_{S\in \mathcal{S}} \|f^{\nsim}_S\|_{L^2}^2 \lesssim d^{j_t-t} \mathit{Poly}(d)^{t-1} R^{3\delta} \|f\|_{L^2}^2.
\end{equation}
Since $1\leq t\leq n\leq \delta^{-2}$,  $d\lesssim \log R$ and $D=d^{j_t}$, we have 
\begin{equation}\label{eq: keyl2}
\sum_{S\in \mathcal{S}(\kappa)} \|f^{\nsim}_S\|_{L^2}^2\leq D R^{4\delta} \|f\|_{L^2}^2. 
\end{equation}
Inserting \eqref{eq: keyl2} into \eqref{eq: est01},  we obtain the first estimate on $\|Ef\|_{BL^p(B_R)}^p$:
\begin{equation}\label{eq: est11}
\|Ef\|_{BL^p(B_R)}^p \lesssim R^{O(\epsilon_0)} D r^{\frac{5}{2}-\frac{3p}{4}}r^{-\frac{p-2}{4}} \|f\|_{L^2}^2 \cdot \max_{S\in \mathcal{S}(\kappa)} \|f^{\nsim}_{\Pi_S, \tau(S)}\|_{L^2_{avg}(\tau(S))}^{p-2}. 
\end{equation}
Combining \eqref{eq: choosekappa} with Lemma~\ref{lem: keylarge} and taking average on $\tau(S)$,  we have
\begin{equation}\label{eq: keyest}
\max_{S\in \mathcal{S}(\kappa)} \|f^{\nsim}_{\Pi_S, \tau(S)}\|_{L^2_{avg}(\tau(S))}^2 \lesssim (\frac{R}{r})^{-1/2} R^{O(\epsilon_0)} \max_{d(\tau)=r^{-1/2}}\|f_{\tau}\|_{L^2_{avg}(\tau)}^2. 
\end{equation}
Inserting \eqref{eq: keyest} in \eqref{eq: est11}, we conclude that 
\begin{equation}\label{eq: est1}
\|Ef\|_{BL^p(B_R)}^p \lesssim R^{O(\epsilon_0)} D r^{\frac{5}{2}-\frac{3p}{4}}R^{-\frac{p-2}{4}} \|f\|_{L^2}^2 \cdot \max_{d(\tau)=r^{-1/2}} \|f_{\tau}\|_{L^2_{avg}(\tau)}^{p-2}. 
\end{equation}

\vspace{5pt}

\textbf{(II).} For the second approach, we sort the fat $r$--surfaces $S\in \mathcal{S}$ according to  the size 
$$\sum_{O\in \mathcal{O}(t), \, O\subset S} \|Ef\|_{ BL^p(O)}^p,$$  
where the set $\mathcal{O}(t)$ is as defined in Lemma~\ref{lem: findscales}.  By pigeonholing and  \eqref{eq: indc}  in Lemma~\ref{lem: recursion} with $j=j_t$ and $O_j=S_t$, there exist a dyadic number $R^{-10} \leq  \lambda \leq 1$ and  a subset $\mathcal{S}'\subset \mathcal{S}$ such that for each $r$--fat surface  
$S\in \mathcal{S}'$
$$\sum_{O\in \mathcal{O}(t), \, O\subset S} \|Ef\|_{ BL^p(O)}^p \sim \lambda \|Ef\|_{BL^p(B_R)}^p, $$ 
 and
\begin{align*}\|Ef\|_{BL^p(B_R)}^p &\lesssim  2^{j_t} (\log R)^t \sum_{S\in \mathcal{S}} \|Ef\|_{BL^p(S)}^p \\
	&\lesssim 2^{j_t} (\log R)^t R^{\delta} \sum_{S\in \mathcal{S}}  \,\, \sum_{O\in \mathcal{O}(t), \, O\subset S} \|Ef\|_{ BL^p(O)}^p\\&\lesssim 2^{j_t} (\log R)^{t+1} R^{\delta}  \sum_{S\in \mathcal{S}'}  \,\,\sum_{O\in \mathcal{O}(t), \, O\subset S} \|Ef\|_{ BL^p(O)}^p.
	\end{align*}

Since $j_t\leq \frac{\log R}{\log \log R}$ and $t\leq \delta^{-2}$, for $R$ sufficiently large, $ 2^{j_t} (\log R)^t \log R\lesssim R^{\delta}$. 
Using \eqref{eq: indd} in Lemma~\ref{lem: recursion} with  $j=j_t$ and the fact that  $t\leq \delta^{-2}$ and $d\lesssim \log R$,  we derive that
\begin{equation}\label{eq: card}
|\mathcal{S}'|\gtrsim R^{-2\delta} d^{-3t}  D^3\gtrsim R^{-3\delta} D^3.
\end{equation}
So for each $S\in \mathcal{S}'$, 
\begin{equation}
\|Ef\|_{BL^p(B_R)}^p \lesssim R^{2\delta} |\mathcal{S}'| \sum_{O\in \mathcal{O}(t), \, O\subset S} \|Ef\|_{BL^p(O)}^p.
\end{equation}
By \eqref{eq: nsim} and  \eqref{eq: 55}, 
\begin{equation}
\|Ef\|_{BL^p(B_R)}^p \lesssim R^{2\delta} |\mathcal{S}'| \sum_{O\in \mathcal{O}(t), \, O\subset S} \|Ef^{\nsim}_S\|_{BL^p(O)}^p\lesssim R^{2\delta} |\mathcal{S}'|\cdot \|Ef^{\nsim}_S\|_{BL^p(S)}^p.
\end{equation}
We apply Lemma~\ref{lem: tao} on the function $f^{\nsim}_S$ to show that
\begin{equation}\label{eq: est02}
\|Ef\|_{BL^p(B_R)}^p \lesssim R^{3\delta} |\mathcal{S}'| r^{\frac{5}{2}-\frac{3p}{4}} \|f^{\nsim}_S\|_{L^2}^p. 
\end{equation}
By (the second half of) inequality~\eqref{eq: 728}, there exists $S\in\mathcal{S}'$ such that 
\begin{equation}\label{eq: l2find}
\|f^{\nsim}_S\|_{L^2}^2 \lesssim \frac{D}{|\mathcal{S}' |} R^{4\delta}|f\|_{L^2}^2. 
\end{equation}

Combining \eqref{eq: est02} with \eqref{eq: l2find} and  using  \eqref{eq: card},  we obtain the second estimate for $\|Ef\|_{BL^p(B_R)}^p$:
\begin{equation}\label{eq: est2}
\|Ef\|_{BL^p(B_R)}^p \lesssim R^{O(\delta)} D^{p/2}|\mathcal{S}'|^{1-p/2} r^{\frac{5}{2}-\frac{3p}{4}}\|f\|_{L^2}^p\lesssim R^{O(\delta)}  D^{3-p} r^{\frac{5}{2}-\frac{3p}{4}} \|f\|_{L^2}^p.
\end{equation}

Having estimated $\|Ef\|_{BL^p(B_R)}^p$ using two different approaches, 
we analyze the resulting estimates  \eqref{eq: est1} and  \eqref{eq: est2}. 
By Lemma~\ref{lem: recursion}, $R_j\leq R_{j-1}/\log R$. Recall that we set at the beginning of the proof of Lemma~\ref{lem: keylarge2} that $r=R_{j_t-1}$ and  $D=d^{j_t}$, so  we have  $r\leq (R\log R)/D $.  Now we combine \eqref{eq: est1} and \eqref{eq: est2} into a single estimate:
\begin{equation}
\|Ef\|_{BL^p(B_R)}^p \lesssim R^{O(\epsilon_0)}  \min\{  D R^{-\frac{p-2}{4}}, D^{3-p} \} (\frac{R}{D})^{\frac{5}{2}-\frac{3p}{4}} \|f\|_{L^2}^p. 
\end{equation}
The maximum of the quantity $\min \{  D R^{-\frac{p-2}{4}}, D^{3-p} \}$ is attained when $D=R^{1/4}$. In this case,  when $p \geq 3+3/13$, the constant term is bounded by
\[ R^{O(\epsilon_0)} D^{3-p} D^{3(\frac{5}{2}-\frac{3p}{4})} < C_{\epsilon} R^{\epsilon},\]
which completes the proof of Lemma~\ref{lem: keylarge2}.
\end{proof}

\subsection{The $R^{\delta}\leq r\leq R^{1/2+\delta}$ case.}
The proof in this subsection is similar to the case addressed in the previous subsection and yet is simpler.

\begin{lemma}\label{lem: keysmall}
	If $R^{\delta}\leq r\leq R^{1/2+\delta}$, $\tau$ is a cap of radius $r^{-1/2}$ and  $\Sigma\in \Omega(\tau)$, then for each $\kappa=(\alpha, t, b_1, \gamma_1, \dots, b_l, \gamma_l)$ satisfying~\eqref{eq: stop}, 
	\begin{equation}
	\|f^{\nsim}_{\kappa, \Sigma, \mathit{tang}} \|_{L^2(\tau)}^2 \lesssim r^{-1/2}R^{O(\epsilon_0)} \|f_{\tau}\|_{L^2}^2. 
	\end{equation}
\end{lemma}
\begin{proof}
	Let $B_k$ be the ball of radius $R^{1-\epsilon_0}$ containing $\Sigma$ and write $\Sigma_1=\Sigma$. We double count the quantity
	\begin{equation}\label{eq: dcbush}
	\sum_{\Sigma_2\nsubseteq 5B_k} \sum_{(\theta,v): \theta\subset \tau} \chi_{\kappa}(T_{\theta,v}, \Sigma_1) \chi_{\kappa'}(T_{\theta,v}, \Sigma_2),
	\end{equation}
	where $\sum_{\Sigma_2\nsubseteq 5B_k}$ means summing over all $\Sigma_2\in \Omega(\tau)$ such that $\Sigma_2\nsubseteq 5B_k$ and $\sum_{(\theta,v): \theta\subset \tau}$ means summing over all tubes $T_{\theta,v}\in \mathbb{T}_0$ with $\theta\subset \tau$. 
	Define $N$ to be the number  $$N\coloneqq \sum_{(\theta,v): \theta\subset \tau} \chi_{\kappa} (T_{\theta,v}, \Sigma_1).$$ Using the same argument as in the proof of Lemma~\ref{lem: keylarge} when we derive \eqref{eq: lowerbd}, we have 
	\begin{equation}\label{eq: lowerbdbush}
	\sum_{\Sigma_2\nsubseteq 5B_k} \sum_{(\theta,v): \theta\subset \tau} \chi_{\kappa} (T_{\theta,v}, \Sigma_1) \chi_{\kappa'}(T_{\theta,v}, \Sigma_2) \geq \gamma_l N. 
	\end{equation}
	
	To show  an upper bound for \eqref{eq: dcbush}, we shall use a simple geometric observation. Since $\mathit{dist}(\Sigma_1, \Sigma_2) \geq R^{1-\epsilon_0}$ and each $\Sigma_j$ lies in a ball of radius  $R^{1/2+\delta}$,  the number of tubes $T_{\theta,v}$ intersecting both $\Sigma_1$ and $\Sigma_2$ is $\lesssim R^{O(\epsilon_0)}$. Specifically,
	\begin{equation}\label{eq: geobush}
	\sum_{(\theta,v): \theta\subset \tau} \chi_{\kappa}(T_{\theta,v}, \Sigma_1) \chi_{\kappa'}(T_{\theta,v}, \Sigma_2) \lesssim R^{O(\epsilon_0)}.
	\end{equation} 
	Note that unlike in Lemma~\ref{lem: geoob},  the cap $\alpha$ plays no role in the above geometric argument. 
	
	 Let $\tilde{\kappa}=(\alpha, t, b_1, \dots, b_{l-1}, \gamma_{l-1})$ and $\tilde{\kappa}'=(\alpha, t, b_1, \dots, b_{l-1})$. 
	If $\chi_{\kappa'}(T_{\theta,v}, \Sigma_2)=1$, then  $\chi_{\tilde{\kappa}}(T_{\theta,v}, \Sigma_2)=1$ and  for the  bush $\mathcal{U}$ containing $T_{\theta,v}$ rooted at $\Sigma_2$, 
	$$\sum_{T_{\theta',v'}\in \mathcal{U}} \chi_{\tilde{\kappa}}(T_{\theta',v'}, \Sigma_2)\geq b_l.$$
%
	 So we can rewrite \eqref{eq: geobush} as 
	\begin{equation}\label{eq: geobush2}
	\sum_{(\theta,v): \theta\subset \tau} \chi_{\kappa}(T_{\theta,v}, \Sigma_1) \chi_{\kappa'}(T_{\theta,v}, \Sigma_2) \lesssim R^{O(\epsilon_0)} b_l^{-1} \sum_{(\theta,v): \theta\subset\tau} \chi_{\tilde{\kappa}}(T_{\theta,v}, \Sigma_2). 
	\end{equation}
	Define $N_0$ to be the number
	$$N_0 \coloneqq | \{ T_{\theta,v}\in \mathbb{T}_0, \theta\subset \tau\}|.$$
	
If $\chi_{\tilde{\kappa}}(T_{\theta,v}, \Sigma_2)=1$, then 
	\begin{equation}\label{eq: gammabush}
	\sum_{\Sigma'} \chi_{\tilde{\kappa}}(T_{\theta,v}, \Sigma') \leq \sum_{\Sigma'} \chi_{\tilde{\kappa}'}(T_{\theta,v}, \Sigma') \lesssim \gamma_{l-1} (\frac{R}{r})^{\delta}. 
	\end{equation}
	Summing \eqref{eq: gammabush} over all $T_{\theta,v}\in \mathbb{T}_0$ with $\theta\subset \tau$, we have
	\begin{equation}\label{eq: totalbush}
	\sum_{\Sigma'} \sum_{(\theta,v): \theta\subset v} \chi_{\tilde{\kappa}}(T_{\theta,v}, \Sigma')\lesssim \gamma_{l-1} (\frac{R}{r})^{\delta} N_0. 
	\end{equation}
	Summing \eqref{eq: geobush2} over all $\Sigma_2\nsubseteq 5B_k$ and applying inequality~\eqref{eq: totalbush},  we obtain the following upper bound for \eqref{eq: dcbush},
	\begin{equation}\label{eq: upperbdbush}
	\sum_{\Sigma_2\nsubseteq 5B_k} \sum_{(\theta,v): \theta\subset \tau} \chi_{\kappa}(T_{\theta,v}, \Sigma_1) \chi_{\kappa'}(T_{\theta,v}, \Sigma_2) \lesssim R^{O(\epsilon_0)} b_l^{-1} \gamma_{l-1} N_0. 
	\end{equation}
Since $\kappa$ satisfies \eqref{eq: stop}, we have $\gamma_l\geq \gamma_{l-1} (\frac{R}{r})^{-\delta}$. Comparing \eqref{eq: lowerbdbush} and \eqref{eq: upperbdbush}, we have 
\begin{equation}\label{eq: estbush}
N/N_0 b_l \lesssim R^{O(\epsilon_0)}. 
\end{equation}

Now we are ready to estimate $\|f^{\nsim}_{\kappa, \Sigma, \mathit{tang}}\|_{L^2(\tau)}$. We apply Lemma~\ref{lem: bush} for $f^{\nsim}_{\kappa, \Sigma}$ defined in \eqref{eq: defsigma}, which is concentrated on wave packets from $\mathbb{T}^{\nsim}_{\kappa}(\Sigma)$.  Since $|\mathbb{T}^{\nsim}_{\kappa}(\Sigma)|\leq b_l (\frac{R}{r})^{\delta}$ and $\mathbb{T}^{\nsim}_{\kappa}(\Sigma) \subset \mathcal{U}$ for the unique bush $\mathcal{U}$ rooted at $\Sigma$ for a fixed $\tau$, we have 
\begin{align*}
\|Ef^{\nsim}_{\kappa, \Sigma, \mathit{tang}}\|_{L^2(10B_r)}^2 &\lesssim \|Ef^{\nsim}_{\kappa, \Sigma, \mathit{tang}} \|_{L^2(\Sigma)}^2 +\mathit{RapDec}(R)\|f\|_{L^2}^2 \\&\lesssim r^{-1/2+2\delta} R^{\delta} b_l \|Ef^{\nsim}_{\kappa, \Sigma}\|_{L^2(10B_r)}^2 +\mathit{RapDec}(R)\|f\|_{L^2}^2. 
\end{align*}
We can drop the $\mathit{RapDec}(R) \|f\|_{L^2}^2$--term  because it is much smaller than $  r^{-1/2+2\delta} R^{\delta}b_l \|Ef^{\nsim}_{\kappa,\Sigma}\|_{L^2(10 B_r)}^2$. 

Since by Lemma~\ref{lem: l2} $\|Ef^{\nsim}_{\kappa, \Sigma, \mathit{tang}}\|_{L^2(10B_r)}^2\sim r\|f^{\nsim}_{\kappa, \Sigma, \mathit{tang}}\|_{L^2}^2$  and $\|Ef^{\nsim}_{\kappa, \Sigma}\|_{L^2(10B_r)}^2 \sim r\|f^{\nsim}_{\kappa, \Sigma}\|_{L^2}^2$, we have 
\begin{equation}\label{eq: sctebush}
\|f^{\nsim}_{\kappa, \Sigma, \mathit{tang}}\|_{L^2}^2 \lesssim r^{-1/2+2\delta} R^{\delta} b_l \|f^{\nsim}_{\kappa, \Sigma}\|_{L^2}^2. 
\end{equation}
Recall that all wave packets from $\mathbb{T}_0$  have approximately the same $L^2$--norm. And there are $N$ out of $N_0$ tubes $T_{\theta,v}\in \mathbb{T}_0$ with $\chi_{\kappa}(T_{\theta,v}, \Sigma)=1$, we have 
\begin{equation}\label{eq: ssctebush}
\|f^{\nsim}_{\kappa, \Sigma}\|_{L^2}^2 \lesssim N/N_0 \|f_{\tau}\|_{L^2}^2. 
\end{equation}
Combining \eqref{eq: sctebush}, \eqref{eq: ssctebush},  \eqref{eq: estbush} and the fact that  $\delta <\epsilon_0$, we have 
\begin{equation}
\|f^{\nsim}_{\kappa, \Sigma, \mathit{tang}}\|_{L^2(\tau)}^2 \lesssim r^{-1/2} R^{O(\epsilon_0)} \|f_{\tau}\|_{L^2}^2,
\end{equation}
which completes the proof of Lemma~\ref{lem: keysmall}. 
\end{proof}

\begin{lemma}\label{lem: keysmall2}
	Let $\mathcal{O}$ be the tree defined in Corollary~\ref{cor: tree}. Assume that 
	\begin{equation}
	\|Ef\|_{BL^p(O)} \lesssim \|Ef^{\nsim}\|_{BL^p(O)} 
	\end{equation}
	holds for a $\geq(1-R^{-\delta})$ fraction of the leaves $O\in \mathcal{O}$ and \eqref{eq: 55} in Lemma~\ref{lem: findscales} and   \eqref{eq: twin} in Lemma~\ref{lem: find scale} hold   for some $t$ with $R^{\delta} \leq R_{j_t-1}\leq R^{1/2+\delta}$, then inequality~\eqref{eq: induction2} in Theorem~\ref{thm: induction2} holds for the function $f$ when $p\geq 3+1/5$.
\end{lemma}
\begin{proof}
	To ease the notation, write  $r=R_{j_t-1}$,\,  $S=S_t$,  \,$\mathcal{S}=\mathcal{S}_t$ and $D=d^{j_t}$. 

	If $D\leq r^{1/2}$, we apply  \eqref{eq: est11} in the proof of Lemma~\ref{lem: keylarge2} to show that 
	\begin{equation}
	\|Ef\|_{BL^p(B_R)}^p \lesssim R^{O(\epsilon_0)} D r^{\frac{5}{2}-\frac{3p}{4}} r^{-\frac{p-2}{4}} \|f\|_{L^2}^2 \cdot \max_{S\in \mathcal{S}(\kappa)} \|f^{\nsim}_{\Pi_S, \tau(S)} \|_{L^2_{avg}(\tau(S))}^{p-2}. 
	\end{equation}
Since for each $S\in \mathcal{S}(\kappa)$ \eqref{eq: choosekappa} holds,  when $D\leq r^{1/2}$ and $p\geq 3+1/5$, we apply  Lemma~\ref{lem: keysmall}  to obtain
	\begin{align*}
\|Ef\|_{BL^p(B_R)}^p &\lesssim R^{O(\epsilon_0)} D r^{\frac{5}{2}-\frac{3p}{4}} r^{-\frac{p-2}{2}} \|f\|_{L^2}^2 \cdot \max_{d(\tau)=r^{-1/2}} \|f_{\tau}\| \|_{L^2_{avg}(\tau)}^{p-2} \\
&\lesssim R^{O(\epsilon_0)} r^{4-\frac{5p}{4}} \|f\|_{L^2}^2 \cdot \max_{d(\tau)=r^{-1/2}} \|f_{\tau}\| \|_{L^2_{avg}(\tau)}^{p-2} \\
&\lesssim R^{O(\epsilon_0)} \|f\|_{L^2}^2 \cdot \max_{d(\tau)=r^{-1/2}} \|f_{\tau}\| \|_{L^2_{avg}(\tau)}^{p-2}.
\end{align*}

	Otherwise,  $D\geq r^{1/2}$.  We apply  \eqref{eq: est2} in the proof of Lemma~\ref{lem: keylarge2} to derive  that 
\begin{equation}
\|Ef\|_{BL^p(B_R)}^p \lesssim R^{O(\delta)} D^{3-p} r^{\frac{5}{2}-\frac{3p}{4}}\|f\|_{L^2}^p \lesssim R^{O(\delta)} \|f\|_{L^2}^p
\end{equation}
when $p\geq 3+ 1/5$. 
\end{proof}

\section{Proof of Theorem~\ref{thm: induction2}}\label{section: proofmain}
Theorem~\ref{thm: induction2} now follows from the lemmas proved in previous sections.  First we apply Lemma~\ref{lem: recursion} and Corollary~\ref{cor: tree} to obtain  that 
$$\|Ef\|_{BL^p(B_R)}^p \lesssim R^{\delta} \sum_{O\in \mathcal{O}} \|Ef\|_{BL^p(O)}^p, $$
and 
$$\|Ef\|_{BL^p(O)} \sim \lambda_0 \|Ef\|_{BL^p(B_R)} $$ for all leaves $O \in \mathcal{O}.$
Now we apply the decomposition in Subsection~\ref{subsection: te} with the relation $B_k\sim T_{\theta,v}$ defined  in Section~\ref{section: sim}.  Then the triangle inequality implies that 
$$\|Ef\|_{BL^p(O)} \lesssim \|Ef^{\sim}\|_{BL^p(O)} + \|Ef^{\nsim}\|_{BL^p(O)} +\mathit{RapDec}(R)\|f\|_{L^2}.$$

If there exist a $\geq R^{-\delta}$--fraction of the leaves $O \in \mathcal{O}$ such that $\|Ef\|_{BL^p(O)} \lesssim R^{\delta} \|Ef^{\sim}\|_{BL^p(O)}$, then we apply Lemma~\ref{lem: telocal} to finish the proof. Otherwise,  we apply Lemma~\ref{lem: findscales}. If the first case \eqref{case1} of Lemma~\ref{lem: findscales} happens, then we apply Lemma~\ref{lem: small} to conclude the proof. Otherwise, the second case \eqref{case2} of Lemma~\ref{lem: findscales} happens. For the integer $t$ in  \eqref{case2}, we apply
\begin{itemize}
	\item Lemma~\ref{lem: large} if $R_{j_t-1} \geq R^{1-\epsilon_0}$, 
	\item Lemma~\ref{lem: keylarge2} if $R^{1/2+\delta} \leq R_{j_t-1}\leq R^{1-\epsilon_0}$, 
	\item Lemma~\ref{lem: keysmall2} if $R^{\delta} \leq R_{j_t-1} \leq R^{1/2+\delta}$. 
\end{itemize}

The proof of Theorem~\ref{thm: induction2} is now completed.

\section{Planes}\label{section: plane}

In this section, we prove Proposition~\ref{lem: plane}, which we recall here.
 \begin{prop}\label{lem: plane1}
	Let $S$ be a fat $r$--surface.  For a fixed cap $\tau$ of radius $r^{-1/2}$, let $\mathbb{T}_{S,\tau}$ be the set of tubes $T_{\tau, w}\in  \mathbb{T}_{S, \mathit{ess}}$. Then there exists a set $\Omega(S, \tau)$ of  fat $r$--planes $\Sigma$ such that
	\begin{enumerate}
		\item $|\Omega(S, \tau)| \lesssim_d r^{O(\delta)}$,
		\item 
		$\mathbb{T}_{S, \tau} \subset \bigcup_{\Sigma \in \Omega(S, \tau)} \mathbb{T}_{\Sigma}$
		where $\mathbb{T}_{\Sigma}$ is defined as in Definition~\ref{def: fatplane}.
	\end{enumerate}
\end{prop}

Recall that by Definition~\ref{def: fatsurf}, the fat $r$--surface $S$ is the intersection of $N_{r^{1/2+\delta}}Z_S$ with a ball $B_{\rho}$ of radius $\rho\coloneqq r^{1-\delta}$, where $Z_S$ is a union of smooth algebraic surfaces with $\mathit{deg} Z_S\leq d$. 
By decomposing $Z_S$ into $\lesssim \mathit{Poly}(d)$ irreducible components, if necessary, we can assume that $Z_S$  is the zero set of an irreducible polynomial $P$. 

Let $e_i$ denote the unit vector in the direction of $x_i$--axis for $i=1, 2, 3$. After a change of variables, we can assume that $G(\omega_{\tau})= e_3$.  To ease the notation, we write $\mathbb{T}_{\tau}=\mathbb{T}_{S, \tau}$. For $i=1, 2$, we define 
\begin{equation}\label{eq: criterion}
\mathbb{T}_{\tau, i}= \{ T_{\tau,w}\in \mathbb{T}_{\tau}: \exists T_{\tau', w'}\in \mathbb{T}_S \text{ such that  } T_{\tau,w} \cap T_{\tau', w'}\cap Z_S\neq \emptyset, \text{ and } |G(\omega_{\tau})\wedge G(\omega_{\tau'}) \wedge e_i | \gtrsim 1/K .\}
\end{equation}
By the definition of $\mathbb{T}_{S, ess}$ (see Definition~\ref{def: ess}),
\begin{equation}\label{eq: decom}
\mathbb{T}_{\tau}\subset \mathbb{T}_{\tau, 1} \cup \mathbb{T}_{\tau, 2}.
\end{equation} 

We are going to prove  that Proposition~\ref{lem: plane1} holds for the collection of tubes $\mathbb{T}_{\tau, i}$ in the place of $\mathbb{T}_{S, \tau}$ for $i=1,2$.  In what follows, we focus on the collection $\mathbb{T}_{\tau,1}$, since the other case will follow in exactly the same way.   The  proof of Proposition~\ref{lem: plane1} for $\mathbb{T}_{\tau,1}$  is separated into several lemmas:   Lemma~\ref{lem: perturb} $\sim$  Lemma~\ref{lem: last}. 

Lemma~\ref{lem: perturb}  says that one can perturb $P$ slightly such that the zero set of the perturbed polynomial  is close enough to $Z(P)$ in $B_{\rho}$.  This perturbation will be used  in the proof of Lemma~\ref{lem: morse}.

\begin{lemma}\label{lem: perturb}
	Given an irreducible polynomial $P$ of degree $\leq d$ with $\nabla P\neq 0$ on $Z(P)$ and a fixed ball $B_{\rho}$ of radius $\rho$, for any $0<\lambda< r^{-1/2}$,  there exist some  constants  $0<\lambda_1, \lambda_2\leq \lambda$  depending on $P$ and $B_{\rho}$ such that  the polynomial
	\begin{equation}
	Q \coloneqq P+\lambda_1 (x_1^2+ x_3^2)^d +\lambda_2
	\end{equation}
	satisfies
	\begin{enumerate}
		\item the zero set $Z(Q)$ is a  smooth algebraic surface,
		\item $\underset{x\in Z_Q \cap 2 B_{\rho}}{\max}  \mathit{dist}(x, Z(P)), \underset{x\in Z(P)\cap B_{\rho}}{\max} \mathit{dist}(x, Z(Q)) <\lambda$.
	\end{enumerate}
\end{lemma}
\begin{proof}
	Since $|\nabla P |\neq 0$ on the compact set $Z(P)\cap 3B_{\rho}$ (here we think of $3B_{\rho}$ as a closed ball), there exists $\mu_1>0$ such that $|\nabla P|\geq \mu_1$ on $Z(P)\cap 3B_{\rho}$.  Since $|\nabla P|$ is a continuous function, there exists $0<\mu_2<\mu_1$ such that $ |\nabla P| \geq \mu_2$ on $ N_{\mu_2} Z(P) \cap 3B_{\rho}$.  Furthermore, since the closure of  $3B_{\rho} \setminus N_{\mu_2/2} Z(P)$  is compact, there exists $\mu_3>0$ such that $|P(x)|\geq \mu_3$ if $x\in 3B_{\rho} \setminus N_{\mu_2/2} Z(P)$.
	Let $\mu=\min(\mu_1, \mu_2, \mu_3)$. The remainder of the proof is divided into two steps. 
	
    \textbf{Step 1.} We show that there exists $\lambda_0>0$ such that for any $\lambda_1, \lambda_2<\lambda_0$, $$\max_{x\in Z(P)\cap B_{\rho}} \mathit{dist}(x, Z(Q))<\lambda.$$
	For any $x=(x_1, x_2, x_3)\in Z(P)\cap 3B_{\rho}$, without loss of generality, we can assume that $|\partial_{x_1} P(x)|\gtrsim \mu$. We claim that  there exists $\lambda_0>0$ such that for any $\lambda_1, \lambda_2\leq \lambda_0$, there exists $x_1'$ with  $|x_1'-x_1|\leq \lambda $ and $Q(x_1', x_2, x_3)=0$.
	To see this,  suppose that  $Q(x_1, x_2, x_3)\neq 0$. For $\mu_0$ sufficiently small, on the ball of radius $\mu_0$ centered at $x$,  we have 
	\begin{equation}|\partial_{x_1} P|\gtrsim \mu, \, |\partial_{x_1} Q|\gtrsim \mu. 
	\end{equation}
	By choosing $\lambda_0$ sufficiently small,  we do Taylor expansion on $Q(x)$ at $(x_1, x_2, x_3)$ with respect to $x_1$, 
	\begin{align*}
	Q(x_1'', x_2, x_3) &= Q(x_1, x_2, x_3)+ \partial_{x_1} Q(x) (x_1''-x_1) + O((x_1''-x_1)^2). \\
	&= O(\lambda_0) + \partial_{x_1} Q(x) (x_1''-x_1) + O((x_1''-x_1)^2)
	\end{align*}  
	Then there exists $x_1''$ such that $Q(x_1'', x_2, x_3)$ has 
	the opposite sign as $Q(x_1, x_2, x_3)$ for some $|x_1''-x_1| \lesssim \frac{\lambda_0}{\mu}< \lambda $. By the mean value theorem, there exists $x_1'$ between $x_1$ and $x_1''$ such that $Q(x_1', x_2, x_3)=0$. 
	By  the compactness of $Z(P)\cap B_{\rho}$ (here we regard $B_{\rho}$ as a closed ball), there exists a uniform $\lambda_0>0$ for all $x\in Z(P)\cap B_{\rho}$.
	
	\textbf{Step 2.} We show that by choosing $\lambda_0>0$ sufficiently small, 
	$$\max_{x\in Z(Q)\cap 2B_{\rho}} \mathit{dist}(x, Z(P))<\lambda.$$
	For any $x\in Z(Q)\cap 2B_{\rho}$, we have $P(x)\lesssim \lambda_0$, where the implicit constant in $\lesssim$ depends on the ball $B_{\rho}$ and the coefficients of the polynomial $P$. 
	By choosing $\lambda_0$ sufficiently small, we can assume that 
	$$|P(x)|\leq \mu/100 \text{ for any } x\in Z(Q)\cap 2B_{\rho}.$$
Then $x\in 3B_{\rho}\setminus N_{\mu_2/2} Z(P)$, which implies that on the ball of radius $\mu/2$ centered at $x$, $|\nabla P|\geq \mu$.  The same argument as above shows that when $\lambda_0>0$ is sufficiently small, for any $x\in Z(Q)\cap 2B_{\rho}$, there exists $x'$ with $|x'-x|<\lambda$  such that $P(x')=0$. 
	
	Finally, by Sard's theorem, for a generic $\lambda_2\leq \lambda_0$,  $Z(Q)$ is a smooth algebraic surface.
\end{proof}
\begin{cor}\label{cor: perturbt}
	Let  $0<\lambda<r^{-1/2}$ be a sufficiently small  number,   $Q$ as in Lemma~\ref{lem: perturb} and $S$ as in Proposition~\ref{lem: plane1}. For any $T_{\tau,w}\in \mathbb{T}_S$, we have 
	\begin{enumerate}
		\item $1.1 T_{\tau,w} \cap Z(Q)\cap 1.1B_{\rho}\neq \emptyset; $
		\item for any point $z\in 9.9 T_{\tau,w }\cap 1.9 B_{\rho}\cap Z(Q)$, 
		$$\mathit{Angle}(T_z Z(Q), G(\omega_{\tau})) \leq 1.1 r^{{-1/2}+2\delta}.$$
	\end{enumerate}
\end{cor}
\begin{proof}
	 Recall that $\mathbb{T}_S$ consists of tubes $T_{\tau,w}$ such that 
	 \begin{itemize}
	 	\item $T_{\tau,w}\cap Z(P)\cap B_{\rho}\neq \emptyset$.
	 	\item For any smooth point $z\in Z(P)\cap 2B_{\rho}\cap 10 T_{\tau,w}$, $$\mathit{Angle}(T_z Z(P), G(\omega_{\tau})) \leq r^{-1/2+2\delta}.$$
	 \end{itemize}
 If $z_0\in T_{\tau,w}\cap Z(P)\cap B_{\rho}\neq \emptyset$, then by Lemma~\ref{lem: perturb}, there exists $z_0'\in Z(Q)$ such that $|z_0'-z_0|\leq \lambda <r^{-1/2}$, so 
 $$z_0' \in 1.1T_{\tau,w} \cap Z(Q)\cap 1.1  B_{\rho}.$$
 
 For any $z\in 9.9 T_{\tau,w}\cap 1.9 B_r\cap Z(Q)$, by Lemma~\ref{lem: perturb}, there exists $ z' \in Z(P)$ such that $|z'-z|\leq \lambda$. So 
 $z'\in 10T_{\tau,w}\cap 2 B_r\cap Z(P)$, and we have 
 $$\mathit{Angle}(T_{z'} Z(P), G(\omega_{\tau})) \leq r^{-1/2+2\delta}.$$
 
 When $\lambda>0$ is small enough,  $|\nabla P(z')-\nabla Q(z')| \leq r^{-1/2}$ for any $z'\in 2B_{\rho}$ and $|\nabla Q (z')-\nabla Q (z)|\leq r^{-1/2}$ for any $|z-z'|\leq \lambda $. Finally we obtain
  $$\mathit{Angle}(T_{z} Z(Q), G(\omega_{\tau})) \leq 1.1r^{-1/2+2\delta}.$$
\end{proof}

 Morally speaking, Lemma~\ref{lem: perturb} and Corollary~\ref{cor: perturbt} state that $P$ can be replaced  by the  perturbed polynomial $Q$  when defining $S$ and $\mathbb{T}_S$. 

\begin{lemma}\label{lem: zv}
	There exists a vector  $v_1$ such that $\mathit{Angle}(v_1, e_1) <\lambda $ and  
	\begin{equation}
 Z_{v_1}\coloneqq \{ x: \nabla Q(x) \cdot v_1=0\}\cap Z(Q)
\end{equation}
is either  a smooth algebraic curve or an empty set.  Here $\lambda$ is a small number  as in Lemma~\ref{lem: perturb}, 
\end{lemma}
To prove Lemma~\ref{lem: zv}, we recall the parametric transversality theorem,  which is a version of Sard's theorem. 

\begin{theorem} \cite[Theorem 2.7, Page 79]{Hirsch}\label{thm: Hirsch}
	Let $V, M, N$ be smooth manifolds without boundary and $A\subset N$ a smooth submanifold.  Let $F: V \rightarrow C^{\infty}(M, N)$ satisfy the following condition: 
	\begin{enumerate}
		\item the evaluation map $F^{ev}: V\times M \rightarrow N, (v,x) \mapsto F_v(x)$ is $C^{\infty}$;
		\item $F^{ev}$ is transverse to $A$. 
	\end{enumerate}
	Then the set 
	\[  \cap\kern -0.45em  |\kern0.7em (F; A)\coloneqq \{v\in V: F_v \cap\kern-0.7em|\kern0.7em A\}\] is dense. 
\end{theorem} 
If $f: M\rightarrow N$ and $A \subset N$ a submanifold, then $f$ is said to be  transverse to $A$, which we denote by $f\cap\kern-0.7em|\kern0.7em A$, if for any  $x\in M$ with $f(x) =y \in A$, the tangent space $T_yN$ is spanned by $T_y A$ and the image $ Df_x(T_xM)$. The importance of transversality is highlighted in \cite[Theorem 3.3, Page 22]{Hirsch}, which says that if $f$ is transverse to $A$ and $f^{-1}(A)\neq \emptyset$, then $f^{-1}(A)$ is a submanifold of $M$, and the codimension of $f^{-1}(A)$ in $M$ is the same as the codimension of $A$ in $N$.

Now we are ready to prove Lemma~\ref{lem: zv}.
\begin{proof}
	Without loss of generality, we can assume that $Z(Q)$ is irreducible. Otherwise, we can   decompose  it into a union of irreducible components and work with each component separately. 
	Since $Z(Q)$ is irreducible, by B\'{e}zout's theorem, either $Z(Q)\subset \{\nabla Q\cdot v_1=0\}$ or  $\dim Z_{v_1} \leq 1$.   If $Z(Q) \subset  \{\nabla Q\cdot v_1=0\}$, then for any point $p\in Z(Q)$, $\mathit{Angle}(v_1, T_p Z(Q)) =0$, and 
	the set of $\mathbb{T}_{\tau, 1}$ is empty because of \eqref{eq: criterion}.
	So we can focus on the second case when $\dim Z_{v_1} \leq 1$.  We define the map $F: S^2 \times \mathbb{R}^3\rightarrow \mathbb{R}^2$ as $$F(v_1, x) = (Q(x), \nabla Q(x)\cdot v_1).$$  Suppose that there exists $(v_1, x)$ with $\mathit{Angle}(v_1, e_1)<\lambda$ and $F(v_1, x)=0$; otherwise, the set $Z_{v_1}$ is empty.  Since $Z(Q)$ is smooth and $ T_{(v_1, x)} (S^2 \times \mathbb{R}^3) =T_{v_1} S^2  \oplus \mathbb{R}^3 $, the map 
	\begin{align*}
	 F_{*}: T_{v_1} S^2  \oplus \mathbb{R}^3  &\rightarrow T_{(0,0) }\mathbb{R}^2\\
	  (w_1, w_2)  & \mapsto (  \nabla Q(x) \cdot w_1 + \langle  H(Q) v_1, w_2\rangle  ,  \nabla Q(x) \cdot w_2)
	\end{align*}
	is surjective. Here $H(Q)$ is the Hessian of $Q$.

	Now we apply Theorem~\ref{thm: Hirsch} with $V=S^2$, $M=\mathbb{R}^3$ and $N=\mathbb{R}^2$ and $A =\{(0,0)\}\subset \mathbb{R}^2$.  Since $F_{*}$ is surjective,  the evaluation map $F^{ev}$ is transverse to $\{ 0 \}$. By Theorem~\ref{thm: Hirsch}, the set of $v_1$ such that  the map $F_{v_1}(x) = (Q(x), \nabla Q(x) \cdot v_1)$ is transverse to $\{ 0 \}$ is dense. So such a vector  $v_1$ exists inside  any small cap of radius $\lambda$. 
	Hence,  if $F^{-1}_{v_1}(0)$ is nonempty, then  $F^{-1}_{v_1}(0)$ is a smooth curve in $\mathbb{R}^3$, and $Z_{v_1}$ is a smooth algebraic curve of degree $\lesssim \mathit{Poly}(d)$. 
\end{proof}

\begin{remark}\label{rem: cov}
We can do a change of variables $x_j\mapsto  a_{j1} x_1+ a_{j2} x_2 + a_{j3}x_3$ with $|a_{jj}-1|\lesssim \lambda $ and $ |a_{jk}|\lesssim \lambda $ for $1\leq j\neq k\leq 3$ such that $v_1 $ is changed into $e_1$, and the leading order term of $Q$ as a polynomial of $x_1$ and $x_3$ with coefficients in $\mathbb{R}[x_2]$ becomes 
$$\lambda_1 ((a_{11}x_1+ a_{13}x_3)^2 + (a_{31}x_1+a_{33}x_3)^2)^d, $$
whose zero set in  $ \mathbb{R}^2$ is $\{(0,0)\}$. 

From now on, we  assume that $v_1=e_1$ and 
\begin{equation}
\mathit{Angle}(e_3, G(\omega_{\tau})) \leq r^{-1/2}
\end{equation}
 since $\lambda<r^{-1/2}$.

\end{remark}

For any $t\in \mathbb{R}$, let us  define the plane perpendicular to $e_2$: $$\Sigma_t\coloneqq \{ (x_1, x_2, x_3): x_2=t\}.$$

We consider $Z(Q)\cap \Sigma_t$. If $Z_{v_1}$ is empty, then since $v_1$ is parallel to $\Sigma_t$, the tangent plane $T_p Z(Q)$ at any point $p\in Z(Q)$ is not parallel to $\Sigma_t$.  As a result, $Z(Q)\cap \Sigma_t$ is a smooth curve. For a fixed $t$, the polynomial $Q$ restricted to the plane $\Sigma_t$ as a polynomial of $x_1$ and $x_3$ has leading order term as defined in Remark~\ref{rem: cov}, so $Z(Q)\cap \Sigma_t$ is a compact smooth algebraic curve on the plane $\Sigma_t$.  We find a point $x\in Z(Q)\cap \Sigma_t$ with the largest (or smallest) $x_3$--coordinate. Then  $v_1=e_1\in T_x (Z(Q)\cap \Sigma_t)$.  
This contradicts the assumption that $Z_{v_1}$ is empty. 

From now on, we know that $Z_{v_1}$ is a smooth algebraic curve. 

\begin{lemma}
	There exists  $\eta \in [r^{-1/2+5\delta}, r^{-1/2+15\eta}]$ such that 
$$Z_{v_1, \eta} \coloneqq \{ x\in Z_{v_1}\cap 2B_{\rho}, \,\, \mathit{Angle}(T_x Z_{v_1}, \Sigma_t)\leq \eta\}$$
	 can be decomposed into a union of $\lesssim \mathit{Poly}(d)$ connected components $Z_{v_1, \eta}=\sqcup_j Z_{v_1, \eta, j}$.  Furthermore, each component $Z_{v_1, \eta, j}\subset N_{r^{1/2+2C\delta}} \Sigma_{t_j}$ for some $t_j \in \mathbb{R}$. 
\end{lemma}
\begin{proof}

	We decompose $Z_{v_1, \eta}$ into a union of  connected components $Z_{v_1, \eta, j}$. The end points of $Z_{v_1, \eta, j}$ are the points $p\in Z_{v_1}$ such that either 
	\begin{equation}
	\mathit{Angle} (T_p Z_{v_1}, \Sigma_t )=\eta, 
	\end{equation}
	or $p \in \partial (2 B_{\rho})$. Since $\mathit{Angle}(T_p Z_{v_1}, \Sigma_t )=\eta$ can be expressed as an algebraic equation of degree $\lesssim d$, for a generic $\eta \in [ r^{-1/2+5\delta}, r^{-1/2+15\delta}]$, the set $\{ x\in Z_{v_1}, \mathit{Angle}(T_{x}Z_{v_1}, \Sigma_t)\leq \eta\}$ has $\lesssim \mathit{Poly}(d)$ connected components.
	So
	$Z_{v_1, \eta}$ has $\lesssim \mathit{Poly}(d)$ connected components $Z_{v_1, \eta,j}$. 

 For a fixed  component $Z_{v_1, \eta, j}$ with $Z_{v_1, \eta, j} \cap B_r\neq \emptyset$, we can find a plane $\Sigma_t$ such that $\Sigma_t\cap Z_{\eta,j}\cap 2B_{\rho}\neq \emptyset$. By the definition of $Z_{v_1, \eta}$, for any point  $p\in Z_{v_1, \eta}$, the angle between $T_p Z_{v_1}$ and $\Sigma_t$ is $\leq r^{-1/2+15\delta}$ and $Z_{v_1, \eta,j}\subset  2 B_{\rho}$. Hence, the component $Z_{v_1, \eta,j} \subset N_{r^{1/2+2C\delta}} \Sigma_t$ for some absolute constant $C$. Since there are $\lesssim \mathit{Poly}(d) $ connected components $Z_{v_1, \eta,j}$ of $Z_{v_1, \eta}$, we can find $\lesssim \mathit{Poly}(d)$ planes $\Sigma_{t_1}, \dots, \Sigma_{t_M}$ with $t_1<\cdots < t_M$, such that the union of their  $r^{1/2+2C\delta}$--neighborhoods covers $Z_{v_1, \eta}$.
\end{proof}

\begin{remark}\label{rem: remove}
Note that we can remove all tubes $T_{\tau, w}$ from $ \mathbb{T}_{\tau}$ with  $3T_{\tau, w}\cap \cup_{j} N_{r^{1/2+2C\delta}}\Sigma_{t_j}\neq \emptyset$, and  the removed tubes can be  covered  by  the union $\cup_{j} N_{r^{1/2+4C\delta}} \Sigma_{t_j}$.    Indeed, since  $\mathit{Angle}(e_3,  G(\omega_\tau)) \leq r^{-1/2}$, the tubes $T_{\tau,w}$ are parallel to $\Sigma_t$ up to an angle difference of $\lesssim r^{-1/2}$.  It follows that there exist  $\lesssim r^{O(\delta)}$ fat $r$--planes $\Sigma$, such that the removed tubes lie in $\cup_{\Sigma} \mathbb{T}_{\Sigma}$. 

\end{remark}

\begin{lemma}\label{lem: noboundary}
	For any $T_{\tau,w}\in \mathbb{T}_{\tau}$, we have  $$1.9B_{\rho}\cap 9.9 T_{\tau,w} \cap(  Z_{v_1}\setminus Z_{v_1, \eta})  = \emptyset.$$
\end{lemma}

\begin{proof}
	If $z\in 1.9B_{\rho}\cap 9.9 T_{\tau,w} \cap (Z_{v_1}\setminus Z_{v_1, \eta} ), $ then
	\begin{itemize}
		\item $T_{\tau,w}\in \mathbb{T}_S$ and Corollary~\ref{cor: perturbt} imply that $\mathit{Angle}(T_z Z(Q), G(\omega_{\tau})) \leq 1.1 r^{-1/2+2\delta}$, 
		\item $z\notin Z_{v_1, \eta}$ implies that $\mathit{Angle}(T_{z} Z_{v_1}, \Sigma_t) \geq r^{-1/2+5\delta}$,
		\item $z\in Z_{v_1}$ implies that $\mathit{Angle}( T_z Z(Q), v_1) =0$ .
	\end{itemize}
But the above three angle conditions can not hold simultaneously because $v_1 =e_1$, $\Sigma_t\perp e_2$ and $\mathit{Angle}(e_3, G(\omega_{\tau})) \leq r^{-1/2}$. 
\end{proof}

For any $t_1<t_2$, we define the slab $\Sigma_{t_1, t_2}$ as 
$$\Sigma_{t_1, t_2}\coloneqq \{ (x_1, x_2, x_3): t_1\leq x_2 \leq t_2\}.$$
\begin{lemma}\label{lem: morse}
Define $$Z_{v_1, e_2^{\perp}}\coloneqq \{ x\in Z_{v_1}:  T_xZ_{v_1}\perp e_2\}.$$ Then there exist $s_0\leq s_1\leq \cdots \leq s_N\leq s_{N+1}$, $N\lesssim r^{O(\delta)} $ such that 
\begin{itemize}
	\item$s_{j+1}-s_j \leq r^{1-10\delta}$;
	\item $2B_{\rho}\subset \Sigma_{s_0, s_{N+1}}$, and  $B_{\rho}$ contains the fat $r$--surface $S$ in Lemma~\ref{lem: plane1};
	\item $Z_{v_1, e_2^{\perp}}\subset \cup_{j=1}^{N} \Sigma_{s_j}$;
	\item the set $$\big(Z(Q)\setminus ( Z_{v_1}\bigcup_{j=1}^N N_{r^{1/2+5C\delta}}\Sigma_{s_j})\big)  \bigcap \Sigma_{s_0, s_{N+1}}= \bigsqcup_l Z_l$$ can be decomposed
	into a union of connected components $Z_l$. In addition, the number of connected components $Z_l $ is $\lesssim r^{O(\delta)}$,  and the orthogonal projection $$\Pi_{v_1}: \mathbb{R}^3 \rightarrow \mathbb{R}^2, \, \, (x_1, x_2, x_3) \mapsto (x_2, x_3), $$ is injective on each $Z_l$. 
\end{itemize}
\end{lemma}
\begin{remark}
	Here $Z_{v_1, e_2^{\perp}}$ plays a similar role as $Z_{v_1, \eta}$, except that the points in  $Z_{v_1, e_2^{\perp}}$ are not required to lie in $2B_{\rho}$.  In particular, $Z_{v_1, e_2^{\perp}} \cap 2B_{\rho}\subset Z_{v_1, \eta}$ and 
	\begin{equation}\label{eq: contain}
	Z_{v_1, \eta} \subset \cup_{j=1}^N N_{r^{1/2+5C\delta}} \Sigma_{s_j}. 
	\end{equation}
	We use $Z_{v_1, e_2^{\perp}}$ to get rid of the potential issues associated to $\partial (2B_{\rho})$ in the the proof of Lemma~\ref{lem: morse}. 
\end{remark}

\begin{proof}[Proof of Lemma~\ref{lem: morse}]
	
	We claim that there exists a set of planes  $\{\Sigma_{s_k} \}$, $1 \leq k \lesssim\mathit{Poly}(d) $, such that 
	$$Z_{v_1, e_2^{\perp}}\subset \bigcup_{k}\Sigma_{s_k}.$$
	 To see this,  recall that  $Z_{v_1} = \{\nabla Q\cdot v_1=0\} \cap \{Q=0\}$.  Without loss of generality, we can assume that $Z_{v_1}$ is irreducible. Otherwise, we  decompose $Z_{v_1}$  into $\lesssim \mathit{Poly}(d)$ irreducible components and work with each component.  Since $Z_{v_1}$ is smooth, we can rewrite $Z_{v_1, e_2^{\perp}}$ as 
	$$ Z_{v_1, e_2^{\perp}}  = \{ ( \nabla (\nabla Q \cdot v_1) \times \nabla Q ) \cdot  e_2=0\}  \cap Z_{v_1}.$$
	 By B\'{e}zout's theorem,  since $Z_{v_1}$ is irreducible, $Z_{v_1, e_2^{\perp}}$ is either a union of $\lesssim \mathit{Poly}(d)$ points or $Z_{v_1, e_2^{\perp}}=Z_{v_1}$.  
	 
	 When $Z_{v_1, e_2^{\perp}}$ is a union of $\lesssim \mathit{Poly}(d)$ points, we choose $\lesssim \mathit{Poly}(d)$  planes $\Sigma_{s_k}$ to cover all those points.

	 When $Z_{v_1, e_2^{\perp}}=Z_{v_1}$, we can decompose $Z_{v_1, e_2^{\perp}}$ into $\lesssim \mathit{Poly}(d)$ connected components $Z_{v_1, e_2^{\perp}}= \cup Z_{v_1, e_2^{\perp}, k}$.
	 We show that each connected component $Z_{v_1, e_2^{\perp}, k}$  lies in a plane of the form $\Sigma_{s_k}$  by the definition of $Z_{v_1, e_2^{\perp}}$. Indeed, take two points $p_1, p_2 \in Z_{v_1, e_2^{\perp}, k}$. If $(p_1 - p_2)\cdot e_2 \neq 0$, then there exists a point $p \in Z_{v_1, e_2^{\perp}, k}$ such that  $v\cdot e_2 \neq 0$ for  any nonzero vector  $v \in T_p Z_{v_1, e_2^{\perp}}$. The tangent space $T_p Z_{v_1, e_2^{\perp}}$ is equal to $T_p Z_{v_1}$ in the case of $Z_{v_1, e_2^{\perp}}= Z_{v_1}$. This contradicts the definition of $Z_{v_1, e_2^{\perp}}$ because $\nabla(\nabla Q \cdot v_1)\times \nabla Q \in T_p Z_{v_1}$ is a nonzero vector as $Z_{v_1}$ is smooth.

	To form our collection $s_1\leq \dots \leq s_N$, we take all $s_k$ that occur in $\{\Sigma_{s_k}\}$  and add more $s\in \mathbb{R}$ to this collection, if necessary,  to ensure that  $s_{j+1}-s_j \leq r^{1-10\delta}$. Since $d\leq r^{\delta}$, $ N\lesssim r^{O(\delta)}$. Finally, we choose $s_0, s_{N+1}$ such that $2B_{\rho}\subset \Sigma_{s_0, s_N}$. 
	 
	 \begin{figure}
	 	\centering
	 	\begin{overpic}[scale=0.12 ]{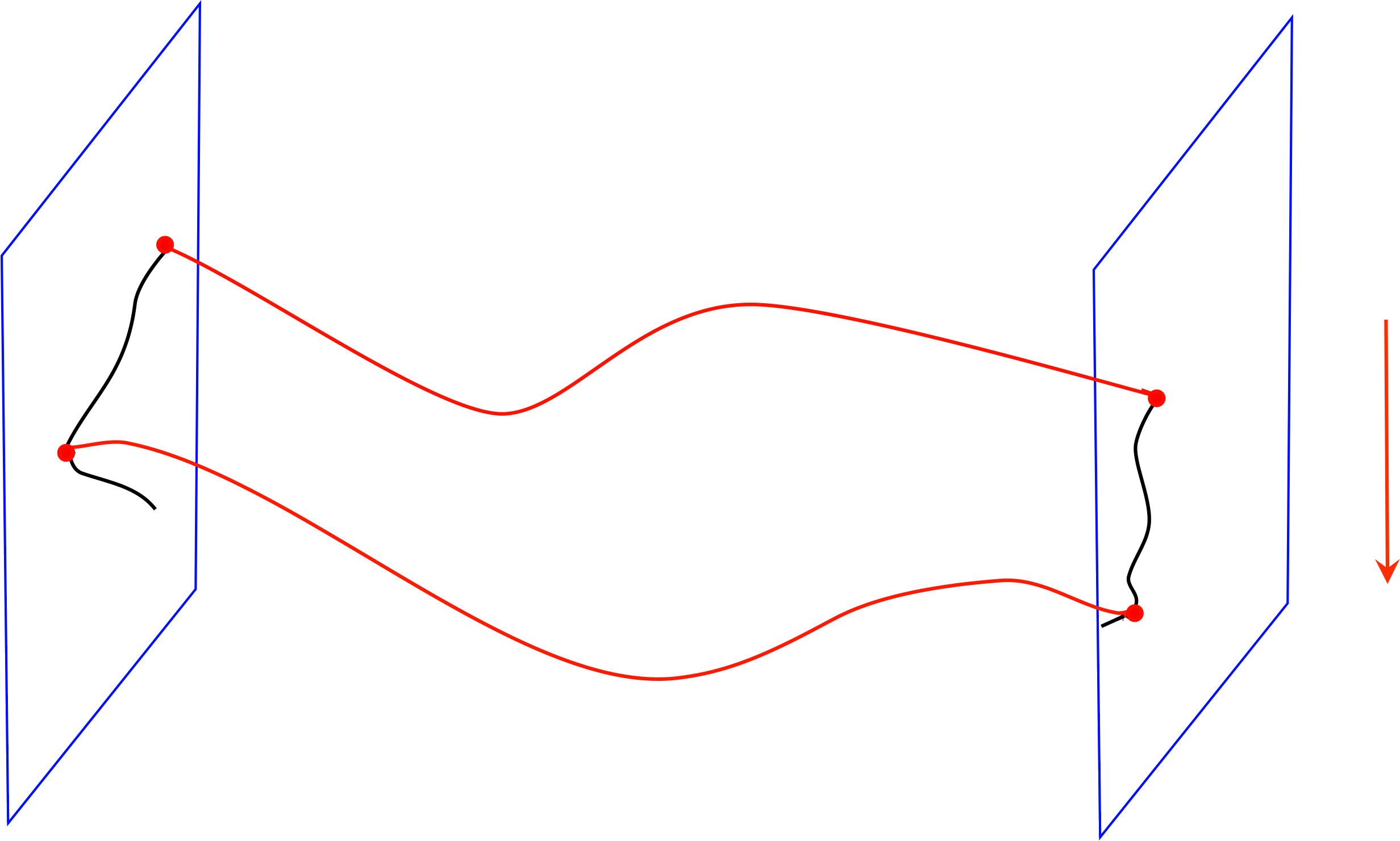}
	 		\put(40, 8){ \color{red} $Z_{v_1}$}
	 		\put(0, 55){\color{blue}$ \Sigma_a$}
	 		\put(90, 60){\color{blue}$\Sigma_b$}
	 		\put(50, 30){$Z_l$}
	 		\put(-5, 35){ $Z_l\cap \Sigma_a$}
	 		\put(95, 30){\color{red} $v_1$}
	 	\end{overpic}
 	\end{figure}
	
	For the rest of the proof, we show that when restricted on each component $Z_l$, the orthogonal projection $\Pi_{v_1}: \mathbb{R}^3\rightarrow \mathbb{R}^2$ is injective.  
	We are going to construct a diffeomorphism $$F: (Z(Q)\cap \Sigma_a)\times [a,b]\rightarrow Z(Q)\cap \Sigma_{a,b}$$ such that 
	\begin{enumerate}
		\item $F((Z_{v_1}\cap \Sigma_a)\times [a,b]) = Z_{v_1}\cap \Sigma_{a,b}$;
		\item\label{en: 2} for any $s\in [a,b]$, $F(( Z(Q)\cap \Sigma_a) \times \{s\}) = Z(Q)\cap \Sigma_s$. 
	\end{enumerate}
	Assuming that we have found such a diffeomorphism, we verify that the restriction of $\Pi_{v_1}$  on each component  $Z_l$ is injective. Note that when restricted on  $Z_l$, the map   $F$ also defines a diffeomorphism: $$(Z_l\cap \Sigma_a) \times (a,b) \rightarrow Z_l.$$ Since $Z_l$ is connected, $Z_l\cap \Sigma_s$ is also connected for each $s\in (a,b)$ because of \eqref{en: 2}.  Using this property, we  will prove that the restriction of  $\Pi_{v_1}$ on $Z_l\cap \Sigma_s$ is injective for each $s\in (a,b)$; it follows that   $\Pi_{v_1}$  is also injective  on $Z_l$. 
	To see this, we claim that 
	\begin{equation}\label{eq: conditioneta}\tag{$\star$}
	\mathit{Angle}(\nabla Q(x), e_2)\geq \eta_0>0  \text{ 	for any } x\in Z(Q)\cap \Sigma_{a,b}, 
	\end{equation}
whose proof will be provided later when we construct $F$. Given \eqref{eq: conditioneta}, $Z_l \cap \Sigma_s$ is a connected smooth curve segment on the plane $\Sigma_s$.  If $p_1, p_2\in Z_l\cap \Sigma_s$ and $\Pi_{v_1}(p_1)=\Pi_{v_1}(p_2)$, then we can find a point $p$ between $p_1$ and $p_2$ on  $Z_l\cap \Sigma_s$ such that $v_1\in T_p(Z_l\cap \Sigma_s)$, and therefore  $p\in Z_{v_1}$. This is impossible because $Z_l \cap Z_{v_1}=\emptyset$. 

\vspace{5pt}
	
	For the rest of the proof, we focus on the proof of  \eqref{eq: conditioneta} and the construction of  the diffeomorphism $F$.  To do so, we shall use a Morse theory type argument, based on the proof of  \cite[Theorem 2.2, Page 153]{Hirsch}.

	Now we regard  $Q$ as a polynomial of $x_1$ and $x_3$ with coefficients in $\mathbb{R}[x_2]$, then the leading order term of $Q$ is (after the change of variables in Remark~\ref{rem: cov})	$$\lambda_1 ((a_{11}x_1+ a_{13}x_3)^2 + (a_{31}x_1+a_{33}x_3)^2)^d, $$
	 with $|a_{11} -1|, |a_{33}-1|, |a_{13}|, |a_{31}|\leq \lambda \leq 0.1$. So  the set
	$$Z(Q)\cap \Sigma_{s_0, s_{N+1}}$$ is compact.

		For each connected component $Z_l$, its boundary $$\partial Z_l\subset \Sigma_{a}\cup \Sigma_b\cup Z_{v_1}$$ for some $a$ and $b$ with  $0<b-a<r^{1-10\delta}$.  
	  We use the function $$f: Z(Q)\cap \Sigma_{a, b} \rightarrow \mathbb{R},\quad  (x_1, x_2, x_3)\mapsto x_2 $$ as the  Morse function on $Z(Q)\cap \Sigma_{a,b}$. 
	We are going to construct a smooth vector field $X$ on $Z(Q)\cap \Sigma_{a,b}$ such that
	the flow generated by $X$ sends a level set of $f$ to another level set of $f$ and preserves the curve $Z_{v_1}$. 
	The vector field $X$ is required to satisfy the following conditions:
	\begin{enumerate}
		\item for any $x\in Z(Q)\cap \Sigma_{a, b}$,  $$X(x)\cdot \nabla^{Z(Q)} f(x) =1,$$
		where $\nabla^{Z(Q)}f$ is the gradient of $f$ on the surface $Z(Q)$, which is viewed as a Riemannian submanifold of $\mathbb{R}^3$; 
		\item \label{preserve} on  $Z_{v_1}\cap \Sigma_{a, b}$, $$X(x)= \frac{\nabla Q\times \nabla (\nabla Q \cdot e_1)}{ \nabla^{Z(Q)} f \cdot (\nabla Q \times \nabla (\nabla Q \cdot e_1))};$$ 
		\item outside a neighborhood of $Z_{v_1}\cap \Sigma_{a, b}$, 
		$$X(x)= \frac{\nabla^{Z(Q)} f (x) }{| \nabla^{Z(Q)} f(x) |^2}.$$
	\end{enumerate}

To see such a vector field $X$ exists, it suffices to show that  there exists $\eta>0$ such that
\begin{enumerate}[label=(\roman*)]
	\item \label{i} $|\nabla^{Z(Q)} f \cdot (\nabla Q \times \nabla (\nabla Q \cdot e_1))| \geq \eta $ on $Z_{v_1}\cap \Sigma_{a, b}$; 
	\item \label{ii}$| \nabla^{Z(Q)} f(x) |\geq \eta$ on $Z(Q)\cap \Sigma_{a,b}$. 
\end{enumerate}
Then  we take a smooth cutoff function $\phi(x)$ on $Z(Q)\cap \Sigma_{a,b}$, such that $\phi=1$ on $Z_{v_1}\cap Z(Q)\cap \Sigma_{a,b}$ and $\phi=0$ outside a sufficiently  small neighborhood of $Z_{v_1}\cap Z(Q)\cap \Sigma_{a,b}$,  and we define 
\begin{equation}\label{eq: X}\tag{$\star \star $}
X(x)\coloneqq \phi(x) \frac{\nabla Q\times \nabla (\nabla Q \cdot e_1)}{ \nabla^{Z(Q)} f \cdot (\nabla Q \times \nabla (\nabla Q \cdot e_1))}+ (1-\phi(x)) \frac{\nabla^{Z(Q)} f (x) }{| \nabla^{Z(Q)} f(x) |^2}.
\end{equation}

Let us verify \ref{ii} first. 
Recall that  $\nabla^{Z(Q)}f(x)$ is the projection of $e_2$ to $T_x Z(Q)$.  We first show a weaker version of \eqref{eq: conditioneta}: 
 \[
 \mathit{Angle}(\nabla Q(x), e_2)>0 \text{ for all } x\in Z(Q)\cap \Sigma_{a,b}.
 \]
 To see this, if $\nabla Q(x)$ is parallel to $e_2$, then   $x\in Z_{v_1}$ because $v_1=e_1$ and $e_1 \perp e_2$.  Furthermore,  $x$ also lies in $ Z_{v_1, e_2^{\perp}}$.  Indeed, since $\nabla Q(x) $ is parallel to $e_2$, we have $e_2 \perp T_x Z(Q)$, which implies that $e_2\perp T_xZ_{v_1}$ because $T_x Z_{v_1}\subset T_x Z(Q)$.  This is impossible because $Z(Q)\cap \Sigma_{a,b} \cap Z_{v_1, e_2^{\perp}}=\emptyset$. 
Since $Z(Q)\cap \Sigma_{a,b}$ is compact, there exists a constant $\eta_0>0$ such that \eqref{eq: conditioneta} holds. By choosing $\eta>0$ sufficiently small  depending on $\eta_0$, we have $$|\nabla^{Z(Q)}f(x)|\geq \eta$$ on $Z(Q)\cap \Sigma_{a,b}$.

Now we verify \ref{i}. Since $\nabla Q \times \nabla (\nabla Q \cdot e_1) (x) \in  T_x Z(Q)$, by the definition of $\nabla^{Z(Q)} f$, 
$$|\nabla^{Z(Q)} f(x) \cdot \big(\nabla Q(x) \times\nabla (\nabla Q \cdot e_1) (x) \big)| = |e_2 \cdot  \big(\nabla Q(x) \times\nabla (\nabla Q \cdot e_1) (x)\big)|.$$
Again, since $\mathit{Angle}(\nabla Q(x), e_2) \geq \eta_0$ for all $x\in Z(Q) \cap \Sigma_{a, b}$, when $\eta$ is sufficiently small depending on $\eta_0$, we have  $$ |e_2 \cdot \big( \nabla Q(x) \times\nabla (\nabla Q \cdot e_1) (x) \big)| \geq \eta. $$

We have shown that the vector field $X$ defined by \eqref{eq: X}  is a bounded vector field on $Z(Q)\cap \Sigma_{a,b}$. By property \eqref{preserve}, $X$ is tangential to the curve  $Z_{v_1}\cap \Sigma_{a,b}$ because $\nabla Q(x) \times \nabla (\nabla Q\cdot e_1)(x) \in T_x Z_{v_1}$ for any $x\in Z_{v_1}\cap \Sigma_{a,b}$.  Now for each $x\in Z(Q)\cap \Sigma_{a,b}$, let $J(x)$ be the maximal interval $J\subset \mathbb{R}$ such that there exists a smooth curve $\gamma: J \rightarrow Z(Q)\cap \Sigma_{a,b}$ with $\gamma'(t) = X(\gamma(t))$ and $\gamma (0)=x$.

For any $s\in [a,b]$, take $x\in f^{-1}(s)$. Since $Z(Q)\cap \Sigma_{a,b}$ is compact, $J(x)$ is closed. And the derivative of the map $$J(x) \rightarrow \mathbb{R}, \, t\mapsto f(\gamma(t))$$ is identically $1$. This means that $f(\gamma(t_1)) - f(\gamma(t_0)) =t_1-t_0$ for any $t_1, t_0 \in J(x)$. So $J(x) = [a-s, b-s]$. 

We know that $f^{-1}(a) = Z(Q)\cap \Sigma_a$ is a union of boundary components of $Z(Q)\cap \Sigma_{a, b}$. Define a map
$$F: f^{-1}(a) \times [a, b] \rightarrow Z(Q)\cap \Sigma_{a, b}, \quad  F(x, t) = \gamma(t-a, x).$$
Since $f$ increases along the trajectories of $X$,  the map $F$ is injective. $F$ is also an immersion because the $X$--trajectories are transverse to level surfaces $\Sigma_t$.  Thus $F$ is an embedding.  $F$ is also  onto because $J(x) = [a-s, b-s]$ is  maximal. This finishes the construction of $F$.
%
\end{proof}

\begin{remark}
	We can remove all tubes $T_{\tau,w}$ from $ \mathbb{T}_{\tau}$ with $$3 T_{\tau,w}\cap  \big(\cup_{j=0}^{N+1} N_{r^{1/2+5C\delta}} \Sigma_{s_j} \big)\neq \emptyset.$$
Same argument as in  Remark~\ref{rem: remove}, there exists $\lesssim r^{O(\delta)}$ fat $r$--planes $\Sigma$ such that  those removed tubes lie in  $\cup_{\Sigma} \mathbb{T}_{\Sigma}$, which is acceptable for Proposition~\ref{lem: plane1}. 
\end{remark}
\begin{lemma}\label{lem: last}
	Let $\mathbb{T}_{\tau, 1, l}$ be the set of tubes 
	$$\mathbb{T}_{\tau, 1, l}\coloneqq \{ T_{\tau,w}\in \mathbb{T}_{\tau}: \exists  T_{\tau', w'} \in \mathbb{T}_S\text{ such that } 1.1 B_{\rho}\cap 1.1 T_{\tau,w}\cap 1.1 T_{\tau', w'} \cap Z_l \neq \emptyset \text{ and } G(\omega_{\tau})\wedge G(\omega_{\tau'})\wedge v_1 \geq r^{-\delta}\}.$$
	Then there exist $\lesssim r^{O(\delta)}$ fat $r$--planes $\Sigma$ such that $ \mathbb{T}_{\tau, 1, l} \subset \cup \mathbb{T}_{\Sigma}$. 
\end{lemma}
\begin{proof}

	After a suitable translation in $\mathbb{R}^3$, we can assume that the ball $B_{\rho}$ containing $S$ is centered at the origin. Recall that $\mathit{Angle}(e_3,  G(\omega_{\tau})) \leq r^{-1/2}$ and $v_1=e_1$.  Consider  the lines on $\mathbb{R}^2 = \Pi_{v_1}(\mathbb{R}^3)$: $$l_{\pm} \coloneqq \{ (x_2, x_3): x_3= \pm 1.3 \rho\}, \, \, l_{a} \coloneqq \Pi_{v_1} (\Sigma_a), \, \, l_b\coloneqq \Pi_{v_1} (\Sigma_b). $$ Let $\mathcal{R}$ be the rectangle defined by 
	$$\mathcal{R}\coloneqq \{ (x_2, x_3): -1.3\rho < x_3< 1.3 \rho, \, \, a < x_2 <b\}. $$ Then the boundary of $\mathcal{R}$ lies in the union of $l_{\pm}$ and $l_a, l_b$. 
	
	 If $\Pi_{v_1} (Z_l) \cap \mathcal{R}=\emptyset$, then $\Pi_{v_1}(Z_l) \cap \Pi_{v_1}(B_{\rho})=\emptyset$ and $\mathbb{T}_{\tau, 1, l}$ is empty. 
	So we can assume that  $\Pi_{v_1} (Z_l) \cap \mathcal{R}\neq \emptyset$. In this case, if 
	\begin{equation}\label{eq: entcontain}
		\Pi_{v_1}(  \partial Z_l \cap Z_{v_1}    )\cap l_{\pm} \cap \mathcal{R} =\emptyset, 
	\end{equation} then $\mathcal{R}$ is entirely contained in $\Pi_{v_1}(Z_l)$.   Otherwise, we can decompose the interval $[a,b]$ into $\lesssim \mathit{Poly}(d)$ subintervals and sort the subintervals into two sets $\mathcal{J}$ and $\mathcal{I}$ such that: \begin{itemize}
		\item for any $J\in \mathcal{J}$ and any $a_2\in J$, $$\Pi_{v_1}(\partial Z_l \cap Z_{v_1}) \cap \{ (x_2, x_3): x_2=a_2,\, -1.3\rho \leq x_3 \leq 1.3 \rho \} =\emptyset.$$
		\item for any $I\in \mathcal{I}$ and any $b_2 \in I$, 
		$$\Pi_{v_1}(\partial Z_l \cap Z_{v_1}) \cap \{ (x_2, x_3): x_2=b_2,\, -1.3\rho \leq x_3 \leq 1.3 \rho \} \neq \emptyset.$$
	\end{itemize}
Indeed, such a decomposition exists.  Without loss of generality, we can assume that $Z_{v_1}$ is irreducible. Otherwise, we decompose $Z_{v_1}$ into a union of irreducible components and work with each component separately. 
Now that $Z_{v_1}$ is irreducible, if $Z_{v_1}\not\subset \Pi_{v_1}^{-1}(l_{+}\cup l_{-})$,  then $Z_{v_1}$ intersects $\Pi_{v_1}^{-1}(l_{+}\cup l_{-})$ at $\lesssim \mathit{Poly}(d)$ points. The projection of those points under $\Pi_{v_1}$ becomes the end points of the intervals in $\mathcal{I}$ and $\mathcal{J}$.  If $Z_{v_1}\subset \Pi^{-1}_{v_1}(l_{+}\cup l_{-})$, then $\mathcal{J}$ is empty and $\mathcal{I}$ contains the interval $[a, b]$.

By Lemma~\ref{lem: noboundary},  for  any $I \in \mathcal{I}$ and  any tube $T_{\tau,w}\in \mathbb{T}_{\tau}$,  $$\Pi_{v_1}(T_{\tau,w})\cap \Pi_{v_1}(Z_l) \cap \{ (x_2, x_3): x_2\in I\} =\emptyset. $$
By changing the interval $[a,b]$ into one of the subintervals $J\in \mathcal{J}$ and $Z_l$ into $Z_l\cap \{ (x_1, x_2, x_3): x_2\in J\}$ and removing tubes $T_{\tau,w}$ as in Remark~\ref{rem: remove},  we can assume that \eqref{eq: entcontain} holds. 
		
	Take a tube $T_{\tau, w^*} \in \mathbb{T}_{\tau, 1, l}$. By the definition of $\mathbb{T}_{\tau, 1, l}$, there exists another tube $T_{\tau', w'}\in \mathbb{T}_S$ such that $ 1.1 B_{\rho} \cap 1.1 T_{\tau, w^*} \cap 1.1 T_{\tau', w'}\cap Z_l \neq \emptyset$ and $|G(\omega_{\tau})\wedge G(\omega_{\tau'}) \wedge v_1|\geq r^{-\delta}$.  
	We are going to show that for any  $T_{\tau, w}\in \mathbb{T}_{\tau, 1,l}, $ the intersection  $3 T_{\tau,w}\cap 3 T_{\tau', w'}\neq \emptyset$. 
	
	Since $|G(\omega_{\tau})\wedge G(\omega_{\tau'}) \wedge v_1|\geq r^{-\delta}$, we have $\mathit{Angle}(\Pi_{v_1} (G(\omega_{\tau'})), e_3) \geq r^{-\delta}$.  In addition, since $0<b-a\leq r^{1-10\delta}$ and  $\Pi_{v_1}(T_{\tau', w'})\cap B_{\rho} \neq \emptyset$, the image $\Pi_{v_1}(3 T_{\tau', w'})$ intersects the boundary of $\mathcal{R}$ at both $l_a$ and $l_b$, but $$\Pi_{v_1} (3T_{\tau', w'}) \cap l_{\pm}=\emptyset .$$ 
Now we apply an argument in the proof of  \cite[Lemma 4.9]{Guth1}.  We draw a curve segment $l(T_{\tau', w'})$ in $Z(Q)$ starting at $z_1 \in 1.1 T_{\tau', w'} \cap 1.1 T_{\tau, w^*} \cap Z_l$ and trying to stay as close as possible to the core line of $T_{\tau', w'}$. The segment has  length $\geq 3 \rho$  and the tangent directions of $l(T_{\tau', w'})$ form an angle of $\leq r^{-1/2+2\delta}$ with $G(\omega_{\tau'})$.    In particular, $l(T_{\tau', w'}) \cap 1.5 B_{\rho}\subset 3 T_{\tau', w'}$. 
	
	For the tube $T_{\tau,w}$, we also draw a curve segment $l(T_{\tau, w})$ in $Z(Q)$ starting at $z_2\in T_{\tau,w}\cap Z_l$ and trying to stay as close as possible to the core line of $T_{\tau,w}$. The segment has length $\geq 3 \rho$, and the tangent directions of $l(T_{\tau,w})$   form an angle of $\leq r^{-1/2+2\delta}$ with $G(\omega_{\tau})$. Note that $\Pi_{v_1}(3T_{\tau,w })$ intersects the boundary of  $\mathcal{R}$ at both $l_{+}$ and $l_{-}$, and $\Pi_{v_1}(3T_{\tau,w })\cap (l_a\cup l_b)=\emptyset$. 
By Lemma~\ref{lem: noboundary}, $$1.9 B_{\rho} \cap 9.9 T_{\tau,w} \cap (Z_{v_1}\setminus Z_{v_1, \eta})=\emptyset. $$ By \eqref{eq: contain},  $l(T_{\tau, w}) \cap 1.5 B_{\rho}$ can not intersect $\partial Z_l \cap Z_{v_1}$.  As a result,  $l(T_{\tau, w}) \cap 1.5 B_{\rho}\subset Z_l$.   Finally, $$\Pi_{v_1}(l(T_{\tau,w}) \cap Z_l )\cap \Pi_{v_1}(l(T_{\tau', w'})\cap Z_l )  \neq \emptyset.$$
By Lemma~\ref{lem: morse}, the projection  $\Pi_{v_1}$ is injective when  restricted on $Z_l$.  Hence,  $l(T_{\tau, w}) \cap l(T_{\tau', w'}) \neq \emptyset$, which  implies that $3T_{\tau,w}\cap 3T_{\tau', w'} \neq \emptyset$. 

Now we have shown that any $T_{\tau,w}\in \mathbb{T}_{\tau, 1, l}$, $3T_{\tau,w}\cap 3T_{\tau', w'} \neq \emptyset$. It follows that the tubes $T_{\tau,w}\in \mathbb{T}_{\tau, 1, l}$ lie in the $r^{O(\delta)}$--neighborhood of the plane spanned by the core line of $T_{\tau', w'}$ and the direction $G(\omega_{\tau})$. 
\end{proof}

Proposition~\ref{lem: plane1} now follows from  the lemmas proved  above, as we explain now. 
\begin{proof}[Proof of Proposition~\ref{lem: plane1}]
	By \eqref{eq: decom},  it suffices to prove Proposition~\ref{lem: plane1} for the set $\mathbb{T}_{\tau,1}$.  Using Lemma~\ref{lem: perturb} and Corollary~\ref{cor: perturbt}, we can replace the zero set $Z_S$ by $Z(Q)$  when defining the fat $r$--surface $S$ up to a small change of constants. We decompose $Z(Q)$ into three parts: 
	\begin{equation}
	Z(Q)= \big(	Z(Q)\cap (\cup_{j=1}^N N_{r^{1/2+5C\delta}}\Sigma_{s_j}) \big) \bigsqcup \big(  Z(Q)\setminus (Z_{v_1}\cup_{j=1}^N N_{r^{1/2+5C\delta}}\Sigma_{s_j}) \big) \bigcup (Z_{v_1} \setminus Z_{v_1, \eta}).
	\end{equation}
	For the first part, we apply Remark~\ref{rem: remove}.  For the third part, a tube $T_{\tau,w}\in \mathbb{T}_{\tau,1}$  can not intersect $Z_{v_1}\setminus Z_{v_1, \eta}$ by  Lemma~\ref{lem: noboundary}.
	 For the second part, using $B_{2\rho}\subset \Sigma_{s_0, s_{N+1}}$ and by  Lemma~\ref{lem: morse}, 
	$$\big(Z(Q)\setminus ( Z_{v_1}\bigcup_{j=1}^N N_{r^{1/2+5C\delta}}\Sigma_{s_j})\big)  \bigcap \Sigma_{s_0, s_{N+1}}= \bigsqcup_l Z_l,$$ 
	  and the orthogonal projection $\Pi_{v_1}$ is injective on  each  connected components $Z_l$. Finally, we apply Lemma~\ref{lem: last} to conclude the proof. 
\end{proof}

\section{Proof of Lemma~\ref{lem: find scale}}\label{section: 38}
In this section, we prove Lemma~\ref{lem: find scale}, which we recall here. 

\begin{lemma}
	Let $\mathcal{O}$ be the tree defined in Corollary~\ref{cor: tree}. If there is a leaf  $O\in \mathcal{O}$ satisfying both \eqref{eq: 55} and \eqref{eq: tang}, 
	then 
	\begin{equation}\label{eq: twin2}
	\|Ef^{\nsim}_{S_t}\|_{BL^p(O)}\sim \|Ef^{\nsim}_{\Pi_{S_t}}\|_{BL^p(O)}.
	\end{equation}
\end{lemma}
\begin{proof}
	For any leaf $O$ satisfying  \eqref{eq: 55} and \eqref{eq: tang} and $1\leq l<t$, we have 
	$$\|Ef^{\nsim}_{S_l}\|_{BL^p(O)} \leq R_{j_{l-1}}^{-\delta} R_{j_{t-1}}^{\delta} \|Ef^{\nsim}_{S_t}\|_{BL^p(O)}.$$
	By \eqref{eq: relpi} in Lemma~\ref{lem: properties}, 
	 $$f^{\nsim}_{\Pi_{S_t}} = f^{\nsim}_{S_t} + \sum_{l=1}^{t-1} f^{\nsim}_{S_l, S_t} +\mathit{RapDec}(R)\|f\|_{L^2}.$$
	By the triangle inequality and the fact that $t\leq \delta^{-2}$, we are left to show  $$\|Ef^{\nsim}_{S_l, S_t}\|_{BL^p(O)} \leq R_{j_{l-1}}^{-\delta} R_{j_{t-1}}^{\delta} \|Ef^{\nsim}_{S_t}\|_{BL^p(O)} \leq R^{-\delta^3}\|Ef^{\nsim}_{S_t}\|_{BL^p(O)}. $$ The last inequality is due to the fact that $R_{j_{l-1}}^{1-\delta} \geq R_{j_{t-1}} \geq R^{\delta}$, while the first inequality follows from Lemma~\ref{lem: negli}:  $$\|Ef^{\nsim}_{S_l, S_t}\|_{BL^p(O)} \leq \|Ef^{\nsim}_{S_l}\|_{BL^p(O)} \leq R_{j_{l-1}}^{-\delta} R_{j_{t-1}}^{\delta} \|Ef^{\nsim}_{S_t}\|_{BL^p(O)}.$$
\end{proof}

\begin{lemma}\label{lem: negli}
	Consider a fat $r_1$--surface $S_1$ and a fat $r_2$--surface $S_2$ with $r_2^{1-\delta}\geq r_1$. Assume that $f$ is concentrated on wave packets from $\mathbb{T}_{S_2}$ and define $$f_{S_1}:= \sum_{T_{\tau_1, w_1}\in \mathbb{T}_{S_1}} f_{\tau_1, w_1}, \quad  f_{S_1, \mathit{trans}}:= \sum_{T_{\tau_1, w_1}\in \mathbb{T}_{S_1}, \mathit{trans} } f_{\tau_1, w_1}. $$
	The following estimates hold for any ball $B\subset S_1$ of radius $K$, where $K$ is a large constant defined in Subsection~\ref{subsection: bl},
	\begin{enumerate}
		\item $\|Ef_{S_1}\|_{BL^p_A(B)} \leq \|Ef\|_{BL^p_{A-2}(B)}$, 
		\item $\|Ef_{S_1, \mathit{trans}}\|_{BL^p_A(B)} \leq \|Ef\|_{BL^p_{A-2}(B)}$.
	\end{enumerate}
\end{lemma}
\begin{proof}
By the definition of a fat $r$--surface, we know that $S_j$ is the $r_j^{1/2+\delta}$--neighborhood of a degree $d$ algebraic surface $Z_{S_j}$ intersecting a ball $B_{r_j^{1-\delta}}$ for  $j=1, 2$.  For any smooth point $z_j\in B+ B_{r_j^{1/2+\delta}}(0) \cap Z_{S_j}$, let $Z_j$ be the tangent plane of $Z_{S_j}$ at $z_j $. By the definition of $f_{S_1}$, for any $T_{\tau_1, w_1}\in \mathbb{T}_{S_1}$ and $T_{\tau_1, w_1}\cap B\neq\emptyset$, we have 
	\[ \mathit{Angle}(G ( \omega_{\tau_j}), Z_j) \leq r_j^{-1/2+2\delta} \leq r_1^{-1/2+2\delta}.\]
	
	Given a cap $\tau \subset \mathbb{R}^2$ of radius $K^{-1}$, the set $G(\tau)\coloneqq \{ G(\omega): \omega \in \tau\}$  gives rise to   a cap on the unit sphere of approximately the same radius in $\mathbb{R}^3$.    For each $j=1, 2$, the directions parallel to $\Sigma_j$ can be represented as a unit circle $\mathcal{C}_j$ in the unit sphere in $\mathbb{R}^3$.  We decompose $\mathcal{C}_2$ into the tangential part $\mathcal{C}_{2, \mathit{tang}} \coloneqq  \mathcal{C}_2\cap N_{r_1^{-1/2+2\delta}} \mathcal{C}_1$, and the transversal part $\mathcal{C}_{2, \mathit{trans}}\coloneqq \mathcal{C}_2\setminus \mathcal{C}_{2, \mathit{tang}}.$
	
	The tangential part $\mathcal{C}_{2, \mathit{tang}} $  contains the directions of tubes in   $\mathbb{T}_{S_1}$ that pass through $B$. We first consider the case when $r_1^{-1/2+2\delta} \leq \mathit{Angle}(Z_1, Z_2) \leq K r_1^{-1/2+2\delta}$.  Since the boundary $\partial \mathcal{C}_{2, \mathit{tang}}$ is zero dimensional, for any cap $G(\tau)$ of radius $K^{-1}$, either $G(\tau)\cap \mathcal{C}_{2, \mathit{tang}}=\emptyset $ or $G(\tau)\cap \mathcal{C}_{2, \mathit{trans}}=\emptyset$,  except for the  two caps intersecting $\partial \mathcal{C}_{2, \mathit{tang}}$.  Note that this property is not necessarily true in higher dimensions: $\partial \mathcal{C}_{2, \mathit{tang}}$ might have dimension $\geq 1$, so there could be many $G(\tau)$ intersecting both the tangential part and the transversal part.
	
	By the definition of the $BL^p$--norm, 
	\begin{equation}\label{eq: muef}
	\|Ef\|_{BL^p_A(B)}^p= \mu_{Ef}(B) = \min_{V_1, \dots, V_A: \text{ lines in } \mathbb{R}^3} \big( \max_{\tau: \mathit{Angle}(G(\omega_{\tau}), V_{a}) \geq K^{-1} \text{ for all } 1\leq a \leq A} \int_B |Ef_{\tau}|^p \big),
	\end{equation}
	where $A$ is a large integer. 
	The quantity $\mu_{Ef}(B)$ takes the $(A+1)$th largest value of $\int_B |Ef_{\tau}|^p$ over all $\tau$ and therefore is equal to zero if $\supp f$ lies in a union of $\leq A$ caps $\tau$. As a result, except for  the two caps intersecting $\partial \mathcal{C}_{2, \mathit{tang}}$, any  cap $G(\tau)$  belongs to  either the tangential part or the transversal part, and 
	\begin{equation}\label{eq: 58}
 \max ( \|Ef_{S_1}\|_{BL^p_A(B)}^p , \|Ef_{S_1, \mathit{trans}}\|_{BL^p_A(B)}^p) \leq 	\|Ef\|_{BL^p_{A-2}(B)}^p. 
	\end{equation}
	When $\mathit{Angle}(Z_1, Z_2) \geq K r_1^{-1/2+2\delta}$, we can find $\lesssim 1$ caps $\tau$ covering $\mathcal{C}_{2, \mathit{tang}}$, so \eqref{eq: 58} holds by the definition of $\mu_{Ef}$. 
\end{proof}
Recall that in \eqref{eq: muef}, we have an underlying  integer $A$ in the definition of  $BL^p_A$--norm. The  integer $A$ can change from line to line because of the subadditive property: 
$$\|Ef +Eg\|_{BL^p_{A_1+A_2}(B)}^p \lesssim \|Ef\|_{BL^p_{A_1}(B)} + \|Eg\|_{BL^p_{A_2}(B)}.$$
Every time we use the triangle inequality of the $BL^p$--norm or apply \eqref{eq: 58}, the value of $A$ is reduced. 

As a consequence,  the inequalities  \eqref{eq: 55} and  \eqref{eq: tang} actually mean
\begin{equation*}
\|Ef^{\nsim}\|_{BL^p_A(O)} \leq R_{j_{t-1}}^{\delta} \|Ef^{\nsim}_{S_t}\|_{BL^p_{A_t}(O)}
\end{equation*}
 where $1\leq t\leq n$, and 
\begin{equation*}
\|Ef^{\nsim}\|_{BL^p_A(O)} \geq R_{j_{l-1}}^{\delta} \|Ef^{\nsim}_{S_l}\|_{BL^p_{A_l}(O)}
\end{equation*}
with $A_l \geq 2^l A_{l-1}$ and $A \geq 2^n A_n$ for all $1\leq l \leq n$.


\vskip.25in

\begin{thebibliography}{3}
\bibitem{Bourgian} J. Bourgain, Besicovitch type maximal operators and applications to fourier analysis, Geometric  and Functional analysis 1 (1991), no. 2, 147-187.

\bibitem{BD} J. Bourgain and C. Demeter, A study guide for the $l^2$ decoupling theorem, Chinese Annals of Mathematics, Series B 38 (2017), no. 1, 173-200.

\bibitem{BG} J. Bourgain and L. Guth, Bounds on oscillatory integral operators based on multilinear estimates, Geometric and Functional analysis 21 (2011), no. 6, 1239-1295.

\bibitem{Fefferman} C. Fefferman, Inequalities for strongly singular convolution operators, Acta Mathematica 124 (1970), no. 1, 9-36.

\bibitem{Guth} L. Guth, Distinct distance estimates and low degree polynomial partitioning, Discrete $\&$ Computational Geometry 53 (2015), no. 2, 428-444.

\bibitem{Guth1} L. Guth, A restriction estimate using polynomial partitioning, Journal of American Mathematical Society 29 (2016), no. 2, 371-413.

\bibitem{Guth2} L. Guth, Restriction estimates using polynomial partitioning II, Acta Mathematica 221 (2018), no. 1, 81-142.

\bibitem{GWZ} L. Guth, H. Wang, R. Zhang, A sharp square function estimate for the cone in $\mathbb{R}^3$, preprint arXiv: 1909.10693(2019).

\bibitem{HR} J. Hickman, K. Rogers, Improved Fourier restriction estimates in higher dimensions, Cambridge Journal of Mathematics, vol. 7, no. 3, pp. 219-282, (2019).

\bibitem{Hirsch} M. Hirsch, Differential topology,  Graduate Texts in Mathematics, No. 33, Springer-Verlag, New York- Heidelberg, 1976.

\bibitem{Milnor} J. Milnor, On the betti numbers of real varieties, Proceedings of the American Mathematical Society 15 (1964), no. 2, 275-280.

\bibitem{Stein1} Elias M Stein, On limits of sequences of operators, Annals of Mathematics (1961), 140-170.

\bibitem{Stein2} Elias M Stein, Some problems in harmonic analysis, Harmonic analysis in Euclidean spaces, 1979, 3-20.

\bibitem{Tao1} T. Tao, The Bochner-Riesz conjecture implies the restriction conjecture, Duke Math Journal 96 (1999), no. 2, 363-375.

\bibitem{Tao2} T. Tao, Endpoint bilinear restriction theorems for the cone, and some sharp null form estimates, Mathematische Zeitschrift 238 (2001), 215-268.

\bibitem{Tao3} T. Tao, A sharp bilinear restriction estimate for paraboloids, Geometric and Functional Analysis 13(2003), no. 6, 1359-1384.

\bibitem{Tao4} T. Tao, Some recent progress on the restriction conjecture, Fourier analysis and convexity, Springer, 2004, 217-243.

\bibitem{Thom} R. Thom, sur l'homologie des vari\'{e}t\'{e} alg\'{e}briques, Differential and Combinatorial Topology (A symposium in Honor of Marston Morse), 1965, pp. 255-265.

\bibitem{Wolff} T. Wolff, A sharp bilinear cone restriction estimate, Annals of Mathematics (2) 153 (2001), no. 3, 661-698.

\bibitem{Zygmund} A. Zygmund, On Fourier coefficients and transforms of functions of two variables, Studia Mathematica 50 (1974), no. 2, 189-201. 


\end{thebibliography}
\end{document}